%
%
%
\magnification 1200
\abovedisplayskip=6pt plus3pt minus3pt  
\belowdisplayskip=6pt plus3pt minus3pt  
\mathsurround=1pt      
\baselineskip = 14pt   
\lineskiplimit = 2pt   
\lineskip = 3pt        
%
%
\def\diagdisplayskip{12pt plus3pt minus3pt}
%
\def\diagbaselines{{\baselineskip20pt \lineskip3pt \lineskiplimit3pt}}
%
%
%
\font\tenmsb=msbm10
\font\sevenmsb=msbm7
\font\fivemsb=msbm5
%
%
\newfam\msbfam                   \textfont\msbfam=\tenmsb       
\scriptfont\msbfam=\sevenmsb     \scriptscriptfont\msbfam=\fivemsb
\def\msb{\fam\msbfam\tenmsb}     
\def\Bbb#1{{\msb #1}}  
\def\scr#1{{\cal #1}}  
%
%
\font\sc=cmcsc10
%
%
\def\domedskip{\par \ifdim\lastskip<\medskipamount
  \removelastskip \medskip\fi}
%
\def\sect#1{\bigbreak \leftline{\bf#1.}\nobreak\smallskip}
%
\def\proc#1{\domedskip\noindent{\bf #1.}\quad\sl}
\def\endproc{\domedskip\rm}
%
\def\prf{\domedskip\noindent{\bf Proof.}\quad}
%
\def\sqr#1#2{{\vcenter{\vbox{\hrule height.#2pt
  \hbox{\vrule width.#2pt height#1pt \kern#1pt \vrule width.#2pt}
  \hrule height.#2pt}}}}
%
\def\endpfsign{\sqr55}
%
\def\qed{{\unskip\nobreak\hfil\penalty50
  \hskip5pt\hbox{}\nobreak\hfil$\endpfsign$
  \parfillskip=0pt\finalhyphendemerits=0\domedskip}\rm}
%
\def\rk#1{\domedskip\noindent{\bf #1.}\quad}
\def\endrk{\domedskip}
%
%
\newcount\sectionnumber  \sectionnumber = 0
\newcount\resultnumber   \resultnumber = 1
%
\def\newsection{\global\advance\sectionnumber by 1 \number\sectionnumber
  \global\resultnumber = 1}
%
\def\section#1{\sect{\newsection.\quad #1}}
%
\def\thelabel{\number\sectionnumber.\number\resultnumber
  \global\advance\resultnumber by 1}
%
\def\proclaim#1{\proc{#1 \thelabel}}

%
\def\label#1{\xdef#1{\number\sectionnumber.\number\resultnumber}%
  \global\advance\resultnumber by 1 #1}
%
\def\labelplus#1#2{{\advance\resultnumber by #1
  \xdef#2{\number\sectionnumber.\number\resultnumber}}#2}
%
\def\checklabel#1{%
  \edef\xxxx{\number\sectionnumber.\number\resultnumber}%
  \global\advance\resultnumber by 1 \xxxx
  \ifx\xxxx#1\else \errmessage{
  The macro \noexpand#1 (= #1) does not correspond 
  to the label \xxxx}\fi}
%
%
\newcount\fignumber  \fignumber = 1
\def\fignum#1{\xdef#1{\number\fignumber}%
  \global\advance \fignumber by 1}
%
%
\newif\ifwithfigs
\withfigstrue
\ifwithfigs

  \ifx\LabelFigloaded\MYundefined\relax
  \else
    \message{ !!! labelfig.tex ALREADY loaded !!!}
   \fi

  \def\LabelFigloaded{\relax}


  \chardef\LabelFigCatAt\the\catcode`\@
  \catcode`\@=11

 \let\LabelFigwlog@ld\wlog
 \def\wlog#1{\relax}

 \ifx\\\MYundefined@
    \let\\\relax
 \fi


  \def\ms@g{\immediate\write16}

 \def\N@wif{\csname newif\endcsname }
 \def\Temp@ {\N@wif\ifIN@}
 \ifx\INN@\MYundefined@
    \else \let\Temp@\relax
 \fi
 \Temp@

  \def\IN@{\expandafter\INN@\expandafter}
  \long\def\INN@0#1@#2@{\long\def\NI@##1#1##2##3\ENDNI@
    {\ifx\m@rker##2\IN@false\else\IN@true\fi}%
     \expandafter\NI@#2@@#1\m@rker\ENDNI@}
  \def\m@rker{\m@@rker}
 
  \newtoks\Initialtoks@  \newtoks\Terminaltoks@
  \def\SPLIT@{\expandafter\SPLITT@\expandafter}
  \def\SPLITT@0#1@#2@{\def\TTILPS@##1#1##2@{%
     \Initialtoks@{##1}\Terminaltoks@{##2}}\expandafter\TTILPS@#2@}

 \def\Shifted@@#1#2#3{\setbox0=\hbox{#3}%
   \raise -\dp0\vbox {\kern-#2%
       \hbox {\kern#1\unhbox0\kern-#1}%
           \kern#2}}

 \newcount\gridcount
 \newbox\auxGridbox@ \newbox\hGridbox@ \newbox\vGridbox@
 \newbox\Labelbox@ \newbox\auxLabelbox@
 \newbox\Coordinatebox@
 \newtoks\Labeltoks@
 \newdimen\Wdd@ \newdimen\Htt@
 \newdimen\Wddd@ \newdimen\Httt@
 
 \def\Wr@{\immediate\write16}

 \newdimen\GL@wd
 \GL@wd=.02pt
 \def\GridLineWidth#1{\GL@wd=#1}

 \def\gobble#1{}
 \def\EdgeErr@{\Wr@{}%
      \Wr@{\string\Edges\space argument
      1, 10, 100 or 1000 please\string!}%
      }

 \newcount\Edgect@

 \def\Sweepup#1\endSweepup{}

 \def\SetEdges@{%
    \edef\Zr@@s{\expandafter\gobble\number\Edgect@\empty}%
        \count255=0\Zr@@s\relax
        \ifnum\count255=\z@\else\EdgeErr@\show\tailtest\fi
        \count255=1\Zr@@s\relax
        \ifnum\count255=\Edgect@\relax\else\EdgeErr@\show\leadtest\fi
    \EdgGl@b\edef\Zr@s{\expandafter\gobble\Zr@@s\empty}
    \ifnum\Edgect@>\@ne\relax\EdgGl@b\let\L@Dc\empty
        \else\EdgGl@b\edef\L@Dc{\string.}\fi
    \ifnum\Edgect@>\@ne\relax
        \EdgGl@b\edef\Edgescale@##1{\divide##1 by \Edgect@}%
        \else\EdgGl@b\edef\Edgescale@##1{}\fi
    }

 \def\Edges#1{\Edgect@=#1\relax
     \let\EdgGl@b\global \SetEdges@}

 \Edges{1}

 \def\hhrule{\hrule height \GL@wd\vskip-.\GL@wd}

 \def\hRule@{%
   \advance\gridcount -2%
   \vfil\hhrule\vfil
   \llap{\smash{\raise -2.5pt
     \hbox{\L@Dc\number\gridcount\Zr@s\kern2pt}}}%
   \hhrule
   }

\def\vvrule{\vrule width \GL@wd \kern-\GL@wd}

 \def\vRule@{\advance\gridcount 2%
   \hfil\vvrule\hfil
   \setbox\auxGridbox@=\vbox to 0pt
      {\vskip \Htt@\vskip 2pt
        \hbox to 0pt{\hss\L@Dc\number\gridcount\Zr@s\hss}\vss}%
      \wd\auxGridbox@=0pt \box\auxGridbox@
   \vvrule
   }

 \def\PlaceGrid@@{\gridcount=10 
  \setbox\hGridbox@=\hbox{%
        \hbox{%
             \hskip-.4pt\vrule
             \vbox to \Htt@{%
               \offinterlineskip\parindent=\z@\relax
               \hbox to \Wdd@{\hfil}
               \hRule@\hRule@\hRule@\hRule@
               \vfil\hhrule\vfil}%
             \vrule\hskip-.4pt}
    }%
  \gridcount=0%
  \setbox\vGridbox@=\hbox{%
      \vbox{\offinterlineskip\parindent=0pt\hsize=0pt
         \vskip-.4pt\hrule%
         \hbox to \Wdd@{%
                 \vtop to \Htt@{\vfil}%
                 \vRule@\vRule@\vRule@\vRule@
                 \hfil\vvrule\hfil}%
         \hrule\vskip-.4pt}}%
  \wd\hGridbox@=0pt\ht\hGridbox@=0pt
  \wd\vGridbox@=0pt\ht\vGridbox@=0pt
  \hbox{\box\hGridbox@\box\vGridbox@}%
  }

 \def\LabelsGlobal{\def\LabGl@b{\global}}
 \def\LabelsLocal{\def\LabGl@b{}}
 \LabelsGlobal 

 \def\SetLabels#1\endSetLabels{%
   \LabGl@b\Labeltoks@={#1()\\}%
   }

 \LabGl@b\Labeltoks@={()\\}

 \def\ShowGrid{\LabGl@b\let\PlaceGrid@\PlaceGrid@@}
 \def\HideGrid{\LabGl@b\let\PlaceGrid@\relax}
 \def\Grids{\ShowGrid\LabGl@b\let\GridSwitch@\ShowGrid}
 \def\noGrids{\HideGrid\LabGl@b\let\GridSwitch@\HideGrid}

 \noGrids

 \def\bAdjust@@{%
     \setbox\auxLabelbox@=\hbox{\raise \dp\auxLabelbox@
            \box\auxLabelbox@}}
 \def\bAdjust@{\let\vAdjust@\bAdjust@@}

 \def\eAdjust@@{\dimen0=-.5\ht\auxLabelbox@
     \advance\dimen0 by .5\dp\auxLabelbox@
     \setbox\auxLabelbox@=
            \hbox{\raise\dimen0\box\auxLabelbox@}}
 \def\eAdjust@{\let\vAdjust@\eAdjust@@}

 \def\tAdjust@@{%
     \setbox\auxLabelbox@=\hbox{\raise-\ht\auxLabelbox@
            \box\auxLabelbox@}}
 \def\tAdjust@{\let\vAdjust@\tAdjust@@}

 \let\vAdjust@\relax

 \def\lAdjust@{\let\hAdjust@\rlap}
 \def\rAdjust@{\let\hAdjust@\llap}

 \let\hAdjust@\relax\let\vAdjust@\relax

 \def\FetchLabel@#1(#2)#3\\{%
     \IN@0#2@@\ifIN@
        \setbox0=\hbox{\ignorespaces#1#3\unskip}%
        \ifdim\wd0>0pt
           \ms@g{}%
           \ms@g{ !!! Bad label(s)? !!!}%
           \message{ #1(#2)#3}%
        \fi
        \def\LabelMole@##1\endFetchLabel@{%
            \IN@0()\\@##1@%
            \ifIN@\def\Temp@{\FetchLabel@##1\endFetchLabel@}%
            \else\def\Temp@{}%
            \fi
            \Temp@
           }%
     \else
       \ignorespaces#1\unskip
       \setbox\auxLabelbox@=%
         \hbox to 0pt{\hss\ignorespaces\hAdjust@
          {\ignorespaces#3\unskip}\hss}%
       \vAdjust@
       \let\hAdjust@\relax\let\vAdjust@\relax
       \AugmentLabelBox@@{#2}%
       \ht\Labelbox@=0pt\dp\Labelbox@=0pt
       \let\LabelMole@\FetchLabel@%
     \fi\LabelMole@}

 \newtoks\XYSep@ 
 \def\SetXYSeparator#1{%
     \IN@0#1@@\ifIN@\XYSep@{*}%
     \else
     \XYSep@{#1}%
     \fi
     }

 \SetXYSeparator*

 \def\AugmentLabelBox@@#1{%
     \IN@0\the\XYSep@ @#1@\ifIN@
       \SPLIT@0\the\XYSep@ @#1@%
       \setbox\Labelbox@=\hbox to 0pt{%
         \unhbox\Labelbox@
         \Shifted@@{\the\Initialtoks@\Wddd@}%
         {\the\Terminaltoks@\Httt@}%
         {\box\auxLabelbox@}}%
     \else
         \ms@g{}%
         \ms@g{ !!! Bad insertion point. !!!}%
         \message{ (#1\ this point was rejected.)}%
     \fi
    }

 \def\FetchOption@#1[#2]#3\endFetchOption@{%
    \def\temp{#1}
    \ifx\temp\empty
       \Edgect@=#2\relax
       \let\EdgGl@b\relax
       \SetEdges@
       \Cleaner@#3%
    \fi}

 \def\Cleaner@#1[@]{\Labeltoks@{#1}}
     
 \def\PlaceLabels@@{\mathsurround=0pt
     \def\Cr@{\\}%
     \let\L\lAdjust@\let\R\rAdjust@
     \let\B\bAdjust@\let\E\eAdjust@\let\T\tAdjust@
     \expandafter\FetchOption@\the\Labeltoks@[@]\endFetchOption@
     \Wddd@=\Wdd@ \Edgescale@\Wddd@ 
     \Httt@=\Htt@ \Edgescale@\Httt@
     \expandafter\FetchLabel@\the\Labeltoks@\endFetchLabel@
     \box\Labelbox@
     }%

 \let \PlaceLabels@\PlaceLabels@@

 \def\AffixLabels#1{\setbox\Coordinatebox@=\hbox{#1}%
      \Wdd@=\wd\Coordinatebox@ \Htt@=\ht\Coordinatebox@
      \advance\Htt@ \dp\Coordinatebox@
      \hbox{\copy\Coordinatebox@\kern-\Wdd@ 
           \Shifted@@{0pt}{-\dp\Coordinatebox@}%
           {\PlaceLabels@\PlaceGrid@}%
           \kern\Wdd@}%
      \GridSwitch@ 
      \LabGl@b\Labeltoks@{()\\}%
      }
 
   \let\wlog\LabelFigwlog@ld   
   \catcode`\@=\LabelFigCatAt  


 
                                By

              Raymond S\'eroul <A18645@FRCCSC21.BITNET>
                                and 
              Laurent Siebenmann <lcs@topo.math.u-psud.fr>
    
              VERSIONS: July 1991, Oct 1991, Jan 1992, July 1992

INTRODUCTION

      This labelling package is intended for TeX users who
rely on non-TeX sources for for their graphics inserts.  It
provides means for adding TeX labels to such inserts with a
minimum of fuss. 

       For most labels, TeX users have in the past found it
reasonably convenient to rely on non-TeX sources. Typical
occasions when an inescapable need for TeX labels seemed to
arise are

 (a) when the graphics program lacks certain exotic or complex
mathematical symbols

 (b) when the very highest typographical quality is wanted for the
labels

 (c) when labels included with the graphics fail to print, 
 and you cannot figure out why (cf. boxedeps.doc).  The labels
 provided by labelfig.tex are 100

       Since this package first appeared, many users, who in the
past scarcely dreamed of using TeX labels, have come to use
nothing but.  So it is now appropriate to add

Intoxication Warning:  TeX labels may be addictive and expensive. 

     If you have a fast preview you may disagree, and even find
that this package provides an agreeable paste-up environment; see
extra applications at end.

     Note to publishers: It is possible and convenient to ultimately
export the TeX labels produced by labelfig.tex to become an integral
part of the EPS file. This is often desired by a publisher who typically
uses an "upmarket" graphics or page layout program, with which the
staff is skilled in perfecting figures.  See Appendix I for
a recipe.

     The authors are grateful to Patrick Ion of Math Reviews for
helpful comments and encouragement.

BASIC INSTRUCTIONS

    After reading in the macro file using

preview or proof your figure with a coordinate grid printed on
top, by typing the following:

    \ShowGrid  
    \AffixLabels{<the graphics insertion>}

Here <the graphics insertion> is what you would type to insert
the graphics object alone without the grid.  This must provide
for the space around it. For example <the graphics insertion>
might well be \BoxedEPSF{MyFigure scaled 700} using the
boxedeps.tex macro package (from same source); this provides a
TeX box containing the encapsulated PostScript insert specified by
the file MyFigure. \AffixLabels{...} provides the grid (supposing
\ShowGrid is present) and later, once you have specified labels
using the grid, it will "tack on" the labels.

     The grid is a sort of (usually elongated) checkerboard of
ten rows and ten columns and its (internal) partitions are by
default numbered  .1, ... ,.9  both horizontally (X-coordinate
running left to right) and vertically (Y-coordinate running bottom
to top).  Thus the points enclosed by the grid correspond to the
points of the unit square in the cartesian "X-Y" plane, the lower
left corner corresponding to the origin (0,0).  By extrapolation,
the full page corresponds to a larger rectangle in the plane.

     These coordinates serve to position labels as follows.
Before the \AffixLabels{...} command type label specifications:

  \SetLabels
   (<X-coordinate>*<Y-coordinate>) <first label> \\
   .
   .
   .
   (<X-coordinate>*<Y-coordinate>)  <last label> \\
  \endSetLabels

Each row specifies one label and is terminated by \\.  In each
row, the position indicator comes first; it is written as a
standard cartesian point except that the X- and Y- coordinates
are separated by * rather than a comma because TeX allows a
comma as decimal point. There are no dimension units to specify
as the unit is the grid itself.

     By default, this cartesian point specifies where the middle
of the baseline of the label will be located.  However if you precede
the point by \L [or \R] the left [or right] edge of the baseline will
be located there. Similarly you may also precede the point by \T, \E,
or \B to vertically align the top equator or bottom of the label box
at the specified point.  This gives nine standard positions of
the label with respect to the insertion point --- corresponding to
the eight principle points of the compas and the center

                     \L\T     \T      \R\T

                     \L\E     \E      \R\E

                     \L\B     \B      \R\B

But this neglects the default "baseline" level of TeX,
giving potentially three more positions

                     \L    <no tag>   \R

For text, the baseline level is often the preferred. Its relation to
the others is variable. It will often coincide with the bottom level,
as happens for "X".  But it is often distinct, as for "g", in which
case you have in all 12 distinct positions rather than 9.

     It is convenient to think of this specification of label
position as attaching the label by a thumb-tack to the coordinate
grid. There are up to twelve positions of the thumb-tack on the
label, while the position of the thumb-tack on the coordinate grid is
arbitrary.  Normally, one choses the position of the thumb-tack on
the label to be the one that is the closest to the item being
labeled.  There are good reasons for this "rule of thumb":

   (a)  It facilitates correct positioning at first try.

   (b)  If the scale of the figure must be altered after labels
have been affixed, the labels have a good chance of remaining well
positioned.

   (c)  The visible grid need not extend beyond the "bounding box"
for the figure, because the best preferred position is always
(at least almost) within the bounding box .

The second reason is particularly important. Indeed it often
happens that scale has to be altered after labelling begins, in
order to either provide space for the labels, or to adjust
proportions between the labels and the figure.  (The size of labels
is unaffected by scaling.)

     Here is an artificial but self-contained test which uses
TeX rules to make a graphics object.

TEST

    Do not skip this!



 \def\FrameIt#1{\hbox{\vrule$\vcenter {\hrule\kern3pt%
             \hbox {\kern3pt #1\kern3pt}%
               \kern3pt\hrule}$\relax\vrule}}

 \def\Caption#1#2{\FrameIt{%
       \vtop {\hsize=#1\relax \parindent=0pt
         \leftskip=0pt \rightskip=0pt plus15pt
         \parfillskip=0pt
         \lineskip=1pt\baselineskip=0pt
         #2}}}

 \def\FirstQuadrant{\hbox to 100pt{\vrule\vbox to 100pt{%
        \hbox to 100pt{\hfil}\vfil\hrule}\hss}}


  \SetLabels
    \R(.5*.2) $\zeta\,\cdot$\\
    (.9*-.10) $\xi$\\
    \R(-.03*.9) $\eta$\\
    \T(.5*.9) \Caption{70pt}{%
          \it The norm of
          $g(\xi+i\eta)$ is indicated on
          contours of this invisible surface.}\\
  \endSetLabels

  \AffixLabels{\FirstQuadrant}

  \end

  Note that the coordinates to use for labels are indicated on the
edges of the grid (when visible) corresponding to the conventional
x- and y- axes of the Cartesian plane. By default the grid is
1-by-1. However, by the command \Edges{100}, you can change this
to 100-by-100 and many users find this alternative most
convenient. Place the command \Edges{...} in your style file (or
header) since its effect is is global. Other possible edge values
are 10 and 1000.

  If you use the command \Edges{...} at all, do so with care.  For
if you accidentally delete an \Edges{...} command your labels will
abruptly be badly misplaced and may logically but mysteriously
generate "dimension too big" errors under TeX and "off page" errors
under your driver.  

  You can dictate the edgescale for an individual figure by giving
the scale in brackets immediately after \SetLabels.  Thus, to
import into an article using say \Edge{100} a figure labelled using
another edgescale, say the original 1-by-1 default, you can use
\SetLabels[1]...\endSetLabels.


GETTING IT DOWN PAT

     Complicated labeling deserves the same respect as
complicated mathematics.  Do not expect it to come out perfect the
first time!  What is needed in either case is a mechanism to
repeatedly typeset troublesome pieces.

     One mechanism is always available.  One does complicated
labelling in a separate "test" file involving just the figure being
labelled;  a texpert will know how to \dump TeX's current state as
a temporary format that restarts rapidly at each retry.  Usually,
one then pastes the completed labelled figure back into the main
TeX file, but, of course, one can also \input it as an auxiliary
file.

     If you do not have a TeXpert at handy, here is a first
approximation to an efficient setup. By deletions reduce a copy
of your article to just a few lines before and after the figure.
Now label the figure, and finally, copy and paste the labelled
figure to the original article. Then copy the next figure to label
into this testbed and repeat. The TeXpert can improve the  speed
at which TeX starts up, by compiling a format specifically for
your article; just one caution: best NOT include in the format
ephemeral details of setup like \Set<mydriver>ArtSpecials (from
boxedeps.tex because this reads  figure dimensions which you may
change during your work session.

     An improved mechanism to repeatedly typeset troublesome
pieces is now available on the Macintosh; it is called LinoTeX;
see the same ftp sources.  It could be set up on many types
of computer.

     Before using labelfig.tex to attach labels to a graphics
object inserted using boxedeps.tex or BoxedArt.tex, make it a
firm rule to carefully adjust the bounding box using the trimming
commands of these packages, and also at least tentatively scale
and position the object. Beware of changing the grid inadvertently
after the labels have been positioned.  For example, correcting
the bounding box of a PostScript graphics object can foul up the
labels by changing the coordinate grid to which the labels are
attached. This is particularly true for the trimming  commands of
boxedeps.tex and BoxedArt.tex. However, as noted already, change
of scale is much less disruptive, and modest adjustments should be
well tolerated.

     Sometimes the labels protrude so far from the bounding box
of a figure that the figure has to be repositioned.  Best do this
by ad hoc spacing, say using \hglue and \vglue; altering the
bounding box would create a vicious circle.

     Remember that you are responsible for preventing labels
from overlapping. You are responsible for all label typography
including size and style. A label is really just about anything
that can be put in a TeX box. Note that spaces at the beginning
and end of labels will normally be suppressed; if you really want
them you must protect them with TeX braces.

     This package temporarily sets the \mathsurround parameter
of TeX to zero  while the labels are being affixed. This is done
because nonzero \mathsurround space would influence the position
of left and right aligned labels; then, when a texpert or printer
modifies mathsurround, diagram labeling might be disastrously
altered. There is a small price to pay involving labels that are
formatted as caption boxes including mathematics: you  may want or
need to specify an explicit mathsurround space within the caption
box; it will not influence anything outside.

     Those hostile to the use of * as separator between
the X and Y coordinates of label insertion points, are free to
impose another using \SetXYSeparator{<the new separator>}.  
Americans may prefer "," to "*" since they never use a 
comma as a decimal point; on the other hand, * may be more visible.

APPENDIX (I)  MERGING labelfig.tex LABELS INTO AN EPSF GRAPHICS OBJECT.

     As promised in the introduction, here is a recipe useful for
publishers. It works at least on Macintosh and at least for vectorized
graphics and Adobe type1 fonts.  (There is surely a similar recipe for
PCs under MSWindows.)

 (a)  Use boxedeps.tex utility to integrate the figure given by the eps
file, "x.eps" say, with a visible frame around it.  See
\ShowDisplacementBoxes command in boxedeps.tex.  To get precise results
automatically it is important to use the \Trim... commands of
boxedeps.tex making the "DisplacementBox" neatly fit the figure.

 (b)  Use the TeX printer driver and LaserWriter (versions >= 8.1.1) to
export to an EPSF the DVI page containing the integrated, labelled
figure. You now have an EPS file  "xx.eps"  that contains too much, and at
the wrong scale, and at wrong position.

 (c)  Convert the EPSF to an Adode Illustrator format EPSF using
the shareware utility called epsConvert by Sam Weiss
1993-- (currently $25).

 (d)  In Illustrator (or a compatible program), group the labels and the
"DisplacementBox"; copy them to the clipboard and paste them into "x.ps".
This step requires that all the label fonts be "visible to the Macintosh.

 (e)  Translate and scale the pasted group consisting of the labels plus
the "DisplacementBox" so as to make the "DisplacementBox" the bounding
box of (labelless) figure represented by "x.eps".  At this point the
labels will be correctly placed on the figure "x.eps".

 (f)  Ungroup and delete the "DisplacementBox".  The result is the
desired single EPS file, "x+.eps" say, It contains the original figure
plus its labels.  

     Using grouping and ungrouping appropriately in "x+.eps", a
publisher's staff can very efficiently improve label positions etc.

APPENDIX II)  SOME EXOTIC APPLICATIONS

     The grid of labelfig.tex is analogous to a light-table in
classical page makeup with wax or latex glue.  In principle, you
can use it to compose any page from its indivisible parts.  This
even has some of the artisanal charm of classical paste-up
provided you have a fast screen preview to make the process
"interactive".

     In practice labelfig.tex is a tool for nonstandard jobs.
Here are a few going beyond the labelling already discussed.

(I)  GRAPHICS INTEGRATION.

     This is accomplished by treating the imported graphics
objects as labels.  The underlying graphics object is then
typically an empty  \vbox to <dimension>{\vfill} in a TeX
\midinsert...\endinsert construction.  A label line
might be of the form

   (.1*.1) \special{... MyFigure ...}\\

The exact form of the special command varies from driver to
driver.  However, in the case of encapsulated PostScript graphics
(EPSF norm), by relying on boxedeps.tex, one can have the
following standard syntax (independant of driver  (see
boxedeps.doc for details.
  
  (.1*.1) \BoxedEPSF{MyFigure scaled <scale in mils>}\\

This may be slow since it requires TeX to read the PostScript
file to read bounding box using many complex macros.  So you
may want to try

  (.1*.1) \EPSFSpecial{MyFigure}{<scale in mils>}\\

which is fast and driver independant, but it squashes the
bounding box, normally to its lower left corner.

     Similarly for graphics of the Macintosh PICT norm ---
using BoxedArt.tex (same sources) in place of boxedeps.tex.

     This approach to integration is to be recommended when
one is assembling a composite graphics object.

 (II)  COMMUTATIVE DIAGRAM ENHANCEMENT

     Commutative diagrams or arrays of mathematical objects
connected by arrows of various sorts are common in mathematics.
The mathematical objects require the use of TeX.  Recently TeX
acquired a good collection of arrows of all slopes --- that of
LamSTeX --- plus pwerful macros to build the diagrams.

     However, even the LamSTeX collection is often
inadequate; it lacks for example double shafted arrows, dotted
arrows and curved arrows. Fortunately it is possible to produce
such arrows on an individual basis using sophisticated graphics
programs such as Illustrator and AldusFreehand (both serving
the EPSF norm) or using Metafont (with its public domain norm).
Since the creation of each new arrow is a work of love, you
probably want to limit the number of arrows by using LamSTeX
for most arrows. The 40K commutative diagram module of LamSTeX
has been adapted to work with AmSTeX and a copy may be posted
with LabelFig and related files. Unfortunately no one has yet
offered a version that works with Plain TeX or LaTeX.

       Suffice it here to say that when the exotic arrow has
been somehow imported into TeX, labelfig.tex treats it as a
label that one affixes to the commutative diagram.  Two other
steps will be treated in separate notes, namely the matter of
extracting the dimension specifications for the arrow and the
construction of the arrow --- for these steps are far from
unique and often depend intimately on your computer environment. 
Notes for the Macintosh-Textures-Illustrator combination are
found in the file ExoticArrows.doc.

 (III) NESTING 

Ingenuity pays off in exploiting labelfig.tex. One can
mix graphics and typography quite freely.  labelfig.tex is good
for freeform or overlapping arrangements, while boxedeps.tex (or
BoxedArt.tex) is best for regimented non-overlapping
arrangements --- and the two can be combined.

     The default behavior of labelfig.tex is not ideal 
for nesting objects, because to prevent trouble for beginners
the register for labels is globally cleared when \AffixLabels
concludes.  But there are switches available

      \LabelsGlobal      \LabelsLocal

which change this.  To understand this, extend the above test 
by something like:


 \LabelsLocal

 \SetLabels
    (.5*.5) AAA\\
 \endSetLabels

 {
 \SetLabels
    (.5*.5) ZZZ\\
 \endSetLabels
   \AffixLabels{\FirstQuadrant}
 }

   \AffixLabels{\FirstQuadrant}


     There are however potential pitfalls.  Neither
labelfig.tex nor boxedeps.tex has been tested under extreme
conditions. Problems may occur if their procedures are
indiscriminately nested. For boxedeps.tex (not labelfig.tex)
there is a precise cause for worry, namely many of its
variables are "global", which means that TeX braces will not
provide the protection one might expect.

COMMAND SUMMARY FOR labelfig.tex

  Here [...] means optional (one or zero)
       [...]* means any number of such constructs

  \SetLabels
    [[<P>](<X><Sep><Y>) <label> \\]*
  \endSetLabels
  \ShowGrid  
  \AffixLabels{<the figure>}

   --- <P> is tack position, one of eleven or empty
              order irrelevant

                   \L\T      \T      \R\T

                   \L\E      \E      \R\E

                     \L               \R

                   \L\B      \B      \R\B

   --- (<X><Sep><Y>) insertion point;
  <Sep> is separator, = * by default;
  \SetXYSeparator{<Sep>} changes it.
   <X> and <Y> are real numbers

  --- <label> a label to attach 

  --- <the figure> the figure to label 

  \GlobalLabels (default)     
  \LocalLabels  setting for nested constructs.

 \Grids makes ALL grids appear; \HideGrid then makes just next disappear.
 \noGrids returns to default.  The commands are always global.

 \GridLineWidth{<dimension>} adjusts width of grid lines. Default is very
small, to give "hairline" effect. If your grid lines are missing try
setting \GridLineWidth{1pt}.

 \Edges#1 globally changes the edge size of all grids to the numerical 
value #1, which must be 1, 10, 100, or 1000.  The default is 1.

VERSION HISTORY.
 --- Jan 1993: \Edges#1 and [??] option after \SetLabels
 --- July 1992: \Grids, \noGrids, \HideGrid;
       Gridlines become hairlines; \GridLineWidth{<dimension>}.
 --- Oct 1991, Jan 1992: \SetXYSeparator{<Sep>},  \LabelsGlobal,
       \LabelsLocal.
 --- July 1991: first release

Address for bugs and other feedback:

        Raymond S\'eroul
        IREM and Lab. de Typographie Informatise
        Univ. Rene Descartes
        Strasbourg

    Tel 33-88-41-63-45
    Email:  A18645@FRCCSC21.BITNET

        Laurent Siebenmann
        Mathematique, Bat. 425,
        Univ de Paris-Sud,
        91405-Orsay,
        France

    Tel 33-1-6941-7949; 
    Email: lcs@topo.math.u-psud.fr  

\input epsf.tex
\chardef\newinsCatAt\the\catcode `\@
\catcode `\@=11
%
%
%
\newskip\insertskipamount\newskip\inserthardskipamount
\insertskipamount 12pt plus2pt     
\inserthardskipamount 4pt          
\def\insertskip{\vskip\insertskipamount}
%
%
\newskip\LastSkip
\def\SaveLastSkip{\LastSkip\lastskip}
\def\RestoreLastSkip{\nobreak\vskip-\LastSkip\vskip\LastSkip}
%
%
\newcount\SplitTest
\def\SetSplitTest{\SplitTest\insertpenalties
  \insert\topins{\floatingpenalty1}%
  \advance\SplitTest-\insertpenalties}
%
%
\def\midinsert{\par
 \SaveLastSkip\penalty-150\SetSplitTest\RestoreLastSkip
 \ifnum\SplitTest=-1
  \@midfalse\p@gefalse\else\@midtrue\fi\@ins}
\def\@ins{\par\begingroup\setbox\z@\vbox\bgroup%
  \vglue\inserthardskipamount}
\def\endinsert{\egroup 
  \if@mid \dimen@\ht\z@ \advance\dimen@\dp\z@
    \advance\dimen@\insertskipamount
    \advance\dimen@\pagetotal\advance\dimen@-\pageshrink
    \ifdim\dimen@>\pagegoal\@midfalse\p@gefalse\fi\fi
  \if@mid%
    \ifdim\lastskip<\insertskipamount\removelastskip\insertskip\fi
    \nointerlineskip\box\z@\penalty-200\insertskip
  \else%
    \SaveLastSkip
    \insert\topins{\penalty100 
    \splittopskip\z@skip
    \splitmaxdepth\maxdimen \floatingpenalty\z@
    \ifp@ge \dimen@\dp\z@
    \vbox to\vsize{\unvbox\z@\kern-\dimen@}
    \else \box\z@\nobreak\insertskip\fi}
    \RestoreLastSkip
   \fi\endgroup}
%
\catcode `\@=\newinsCatAt

\epsfverbosetrue
\def\epsfsizedbox#1#2{\epsfysize=#2 \epsfbox{#1}}
\def\cl{\centerline}

\font\ninemsb=msbm9
\font\sixmsb=msbm6
\font\ninerm=cmr9                   \font\sixrm=cmr6   
\font\ninei=cmmi9                   \font\sixi=cmmi6  
\font\ninesy=cmsy9                  \font\sixsy=cmsy6  
\font\ninebf=cmbx9                  \font\sixbf=cmbx6  
\font\nineex=cmex9                  
\font\nineit=cmti9                  
\font\ninesl=cmsl9                  
\font\ninett=cmtt9                  
\def\small{%
%
%
\textfont0=\ninerm \scriptfont0=\sixrm \scriptscriptfont0=\fiverm
\def\rm{\fam0\ninerm}
%
%
\textfont1=\ninei \scriptfont1=\sixi \scriptscriptfont1=\fivei
%
%
\textfont2=\ninesy \scriptfont2=\sixsy \scriptscriptfont2=\fivesy
%
%
\textfont3=\nineex \scriptfont3=\nineex \scriptscriptfont3=\nineex
%
%
\textfont\bffam=\ninebf \scriptfont\bffam=\sixbf
\scriptscriptfont\bffam=\fivebf \def\bf{\fam\bffam\ninebf}%
%
%
\textfont\itfam=\nineit \def\it{\fam\itfam\nineit}%
\textfont\slfam=\ninesl \def\sl{\fam\slfam\ninesl}%
\textfont\ttfam=\ninett \def\tt{\fam\ttfam\ninett}%
%
%
\textfont\msbfam=\ninemsb \scriptfont\msbfam=\sixmsb
\scriptscriptfont\msbfam=\fivemsb \def\msb{\fam\msbfam\ninemsb}%
%
%
\normalbaselineskip=11pt%
\setbox\strutbox=\hbox{\vrule height8pt depth3pt width0pt}%
%
%
\normalbaselines\rm}    

\def\fig#1#2\endfig{%
\midinsert\cl{#2}\vglue6pt\cl{\small Figure #1}\endinsert}
\def\figwspace#1#2#3#4\endfigwspace{%
\midinsert\vglue#2\cl{#4}\vglue#3\vglue6pt\cl{\small Figure #1}\endinsert}

\def\theFigWyeDelta#1{
\figwspace{#1}{0pt}{14pt}
\small
\SetLabels
\L\E(.16*.7)$x$\\
(.07*.22)$y$\\
(.23*.22)$z$\\
\L\E(.86*.89)$x$\\
(.72*.12)$y$\\
(.98*.12)$z$\\
(.85*.17)$x$\\
\R(.80*.4)$z$\\
\L(.91*.4)$y$\\
\T\E(.5*-.15)$x,y,z\in G$, $xyz=1$.\\
\endSetLabels
\AffixLabels{\epsfsizedbox{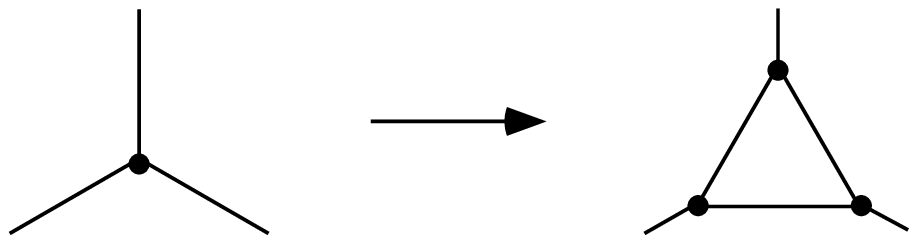}{1truein}}
\endfigwspace
}

\else 
  \def\theFigWyeDelta#1{\relax}
\fi
%
%
%
\def\z{\Bbb Z}               
\def\zmod#1{\z_{#1}}         
\def\zmodb#1{\zmod{#1}}      
\def\zt{\zmod{2}}            
\def\ztb{\zmodb{2}}          
\def\maps#1#2#3{#1\colon #2\to #3} 
\def\cover{\tilde}           
\def\indexset{\scr{C}}       
\def\achar{\delta}           
\def\mchar{\varepsilon}      
\def\shriek{^{\,!}}          
\def\id{\mathop{{\rm id}}\nolimits}         
\def\ker{\mathop{{\rm Ker}}\nolimits}       
\def\im{\mathop{{\rm Im}}\nolimits}         
\def\rank{\mathop{{\rm rank}}\nolimits}     
\def\hom{\mathop{{\rm Hom}}\nolimits}       
\def\lk{\mathop{{\rm Lk}}\nolimits}         
\def\graded{\mathop{{\rm G}}\nolimits}     
\def\veci{\vec\imath}        
\def\vecj{\vec\jmath}        
\def\restrict{\,|\,}         
\def\sep{\,|\,}              
\def\genby#1{\langle #1 \rangle}	
%
\def\smallfrac#1#2{{\textstyle {#1\over #2}}}
\def\half{\smallfrac{1}{2}}
%
\def\smallchoose#1#2{{\textstyle {#1 \choose #2}}}
%
\def\mapright#1{\buildrel #1 \over \longrightarrow}
\def\hookmapright#1{\buildrel #1 \over \hookrightarrow}
\def\mapdown#1{\big\downarrow\rlap
  {$\vcenter{\hbox{$\scriptstyle#1$}}$}}
%
\def\cite#1{[#1]}
\def\Flapan{1}
\def\LeeWein{2}
\def\Massey{3}
\def\Sakuma{4}
\def\Watkins{5}
\def\ZS{6}
%
\def\ThmDTwo{8.1}
\def\ThmDThreeMCirc{8.2}
\def\ThmDThreeMob{8.3}
\def\ThmDFour{8.7}
\def\ThmDFive{8.8}
\def\ResultSeqTaut{6.6}
\long\def\ignore#1\endignore{}
\hbox{}
\bigskip
\centerline{\sc The Homology of Abelian Covers of Knotted Graphs}
\smallskip
\centerline{\sc R.~A.~Litherland}
\bigskip
\rk{Abstract}
Let $\cover M$ be a regular branched cover of a homology
3--sphere $M$ with deck group $G\cong \ztb^d$ 
and branch set a trivalent
graph $\Gamma$; such a cover is determined by a coloring
of the edges of $\Gamma$ with elements of $G$.
For each  index--2 subgroup $H$ of $G$,
$M_H = \cover M/H$ is a double branched cover of $M$. Sakuma has
proved that $H_1(\cover M)$ is isomorphic, modulo 2--torsion,
to $\bigoplus_H H_1(M_H)$, and has 
shown that $H_1(\cover M)$
is determined up to isomorphism by $\bigoplus_H H_1(M_H)$
in certain cases; specifically, when
$d=2$ and the coloring is such that
the branch set of each cover $M_H\to M$ is connected,
and when $d=3$ and $\Gamma$ is the complete graph $K_4$. 
We prove this for
a larger class of coverings: when $d=2$,
for any coloring of a connected graph;
when $d=3$ or $4$, for an infinite class of colored graphs;
and when $d=5$, for a single coloring of the Petersen graph.
\endrk
\rk{AMS Subject Classification} Primary: 57M12.
Secondary: 57M25, 57M15.
\endrk

\section{Introduction} 

For our purposes, a graph is a 1--dimensional
polyhedron $\Gamma$. A vertex of $\Gamma$ is
a point at which $\Gamma$ is not a 1--manifold,
and an edge is the closure of a component of the complement
of the set of vertices. A component of $\Gamma$ that contains a vertex
is naturally a graph in the combinatorial sense (possibly with loops
or multiple edges). A component without vertices is a single edge homeomorphic
to $S^1$, which we call a circular edge (as opposed to a
loop, which is homeomorphic to $S^1$, but contains a vertex). 

\rk{Remark} None of our theorems apply to graphs with circular 
edges, but they are needed for some lemmas.
\endrk

All the graphs we consider are trivalent;
 this does not exclude circular edges. If $\Gamma$
is a trivalent graph the number $V$ of vertices
and the Euler characteristic $\chi(\Gamma)$ 
are related by $V=-2\chi(\Gamma)$, and the number of non-circular edges
is $-3\chi(\Gamma)$. By a cycle in a trivalent graph we mean a 
(possibly empty) subgraph homeomorphic to a disjoint union of
circles; these are in one-to-one correspondence with the elements 
of $H_1(\Gamma;\zt)$. A cycle with one component is called
a circuit.
If $\Gamma'$ is a subgraph
of $\Gamma$, we use $\Gamma\setminus\Gamma'$ to denote the closure
of the set-theoretic complement $\Gamma-\Gamma'$.
We call a graph simple if it has no loops
or multiple or circular edges.

Let $d$ be an integer greater than 1, and let $G$ be a
(multiplicative) group isomorphic to $\ztb^d$.
Let $M$ be a homology 3--sphere 
and let $\maps{\pi}{\cover M}{M}$
be a regular branched cover with deck group $G$
and branch set a graph $\Gamma\subset M$. 
Then $\cover M$ is a manifold iff $\Gamma$ is trivalent;
we assume that this is the case.
For each edge $e$ of $\Gamma$,
the stabilizer $G_e$ of a lift of $e$ to $\cover M$ is a subgroup 
of $G$ of order 2 (and is independent of the lift since $G$ is abelian).
We color $e$ with the non-trivial element of $G_e$. 
The colors $g_1$, $g_2$ and $g_3$ of the edges
at a vertex $v$ are the non-trivial elements of the stabilizer
of a lift of $v$ (a group isomorphic to $\zt\oplus\zt$),
so they satisfy the relation $g_1 g_2 g_3 =1$. 
Conversely,
any coloring of the edges of $\Gamma$ by non-trivial
elements of $G$ satisfying this relation 
at each vertex
defines a homomorphism
$H_1(M-\Gamma;\zt)\to G$ sending each meridian of an edge 
to the corresponding color. The corresponding branched covering
is connected iff the colors
of the edges generate $G$; we shall always assume this is so,
and call $\Gamma$ a $G$--colored graph. 
We regard two $G$--colorings of $\Gamma$ as identical
if they differ only by automorphisms of $G$ and $\Gamma$.
We sometimes write $G(d)$ for $G$ to indicate 
the value of $d$ under consideration.
When we refer to a basis of $G$, or to independent elements
of $G$, we are considering $G$ as a $\zt$ vector space.

\rk{Remark} Any coloring of the edges by elements of $G$
satisfying the above relations defines a homomorphism
from $H^1(\Gamma;\zt)$ to $G$, and vice-versa, so we can 
choose such a coloring with the colors generating $G$
iff the first Betti number $b_1(\Gamma)$ of $\Gamma$ is at least
$d$. However, this may fail to be a $G$--coloring as just defined
since some of the colors may be the identity. If $G$ has a bridge $e$,
the color of $e$ must be $1$ since $e$ represents zero in
$H^1(\Gamma;\zt)$. If $\Gamma$ does not have a bridge,
the existence of a $G$--coloring is not guaranteed;
when $d=2$, a $G$--coloring is just a Tait coloring,
and the question of which bridgeless trivalent graphs
have a Tait coloring has a long history.  
\endrk

Let $\indexset^\star=\indexset^\star(G)$ be the set 
of all subgroups of $G$ of index 2,
and let 
$\indexset = \indexset(G)=\indexset^\star \cup\{G\}$. 
If $\Gamma$ is a $G$-colored graph, for $H\in\indexset$
we let $\Gamma_H$ be the union of the edges of $\Gamma$
whose colors are not in $H$; this is a cycle in $\Gamma$.
If $\Gamma$ is embedded in a homology 3--sphere $M$
with branched cover $\maps\pi{\cover M}M$,
let $M_H=\cover M/H$ for $H\in\indexset$.
If $H\in\indexset^\star$, there is a 2--fold branched covering
$\maps{\rho_H}{M_H}{M}$ whose branch set 
is the link $\Gamma_H$. There is also 
a branched covering $\maps{\pi_H}{\cover M}{M_H}$
with group $H$, whose branch set $\Delta_H$ is the inverse
image of $\Gamma\setminus\Gamma_H$. When $H=G$, $M_G=M$
and we let $\pi_G=\pi$ and $\rho_G=\id$.
Sakuma showed that 
$H_1(\cover M)$ and
$\bigoplus_{H\in\indexset^\star}H_1(M_H)$ 
are isomorphic  modulo 2--torsion
\cite{\Sakuma, Theorem 14.1},
and determined the 2--torsion
of $H_1(\cover M)$ when $d=2$ and each $\Gamma_H$ is connected,
and when $d=3$ and $\Gamma=K_4$ \cite{\Sakuma, Theorem 14.2}.
Our first theorem generalizes part (1) of 
\cite{\Sakuma, Theorem 14.2}, because 
$\bigoplus_{H\in\indexset^\star}H_1(M_H)$
has odd order when all the $\Gamma_H$ are connected,
so the exact sequence of the theorem is split.

\proc{Theorem \ThmDTwo} If $d=2$ and $\Gamma$ is connected, then
there is a short exact sequence
$$0 \to \bigoplus_{H\in\indexset^\star}H_1(M_H)
\mapright\beta H_1(\cover M) \to \ztb^{b_1(\Gamma)-2}\to 0,$$
and 
$\beta\bigl(\bigoplus_{H\in\indexset^\star}H_1(M_H)\bigr)
=2H_1(\cover M)$.
\endproc

There are infinitely many $G(2)$--colorings of connected graphs
for which the $\Gamma_H$ are not all
connected; see Example \labelplus1\ExMobDTwo. In this case,
the above sequence does not split; nevertheless, 
$H_1(\cover M)$ is determined up to isomorphism by
$\bigoplus_{H\in\indexset^\star}H_1(M_H)$. 
This is a consequence of the case $p=2$ and 
$e=1$ of the following proposition, whose proof is a simple
application of the structure theorem for finitely 
generated abelian groups, and is omitted.

\proc{Proposition \label\Algebra} Let $A$ and $B$ be finitely
generated abelian groups, $p$ a prime, and $e$ a positive integer.
If $p^eA \cong p^eB$ and $A/p^eA \cong B/p^eB$, then $A \cong B$.
\qed

\rk{Example \checklabel\ExMobDTwo} 
Let $\Gamma$ be an $n$--rung M\"obius ladder.
Recall that this graph consists of
a $2n$--circuit (the 
rim) together with its diameters (the rungs).
(It is usual to require $n\geq 3$, but the cases $n=1$ or $2$
make sense; when $n=1$ we have the theta-curve, and when
$n=2$ we have $K_4$.) When $n\geq 2$, $\Gamma$ is simple,
and we take the vertices to be
$v_0,\ldots,v_{2n-1}$ and the edges to be
$\sigma_i=\{v_i,v_{i+1}\}$ and 
$\tau_i=\{v_i,v_{i+n}\}$, the subscripts 
being taken modulo $2n$. The $\sigma_i$ form the rim, and the $\tau_i$
are the rungs.
Let $\Gamma'$ be a non-empty cycle in $\Gamma$ that
contains $k$ rungs. If $k=0$, $\Gamma'$ is the rim;
otherwise, $\Gamma'$ is connected if $k$ is odd,
and has $k\over 2$ components if $k$ is even.

Now take $d=2$, and let the non-trivial elements of $G$
be $g_1$, $g_2$ and $g_3$. Give all the rungs the color
$g_1$, and give the edges of the rim the colors
$g_2$ and $g_3$ alternately. If $H=\genby{g_1}$
then $\Gamma_H$ is the rim, while if $H=\genby{g_2}$
or $\genby{g_3}$ then $\Gamma_H$ contains all $n$ rungs.
Thus every $\Gamma_H$ is connected iff $n$ is odd or $n=2$.
\endrk

We say that a $G$--coloring of a graph $\Gamma$ is unsplittable
if, for any $g\in G$, deleting the edges of $\Gamma$
with color $g$ leaves a connected graph. If $\Gamma$ has 
an unsplittable coloring, then either $\Gamma$ is the theta-curve
(in which case $d=2$), or $\Gamma$ is connected and simple.
First, taking $g=1$ shows that $\Gamma$ is connected,
and in particular has no circular edges. Since $\Gamma$
has no bridges, it has no loops. If $\Gamma$ is not the theta-curve
and has a pair of multiple 
edges, these are adjacent to two distinct edges with the same color.
Deleting these edges disconnects $\Gamma$, contrary to the definition.

A circuit $C$ in a $G$--colored graph $\Gamma$ will be called special
if there is some $H\in\indexset^\star$ such that
$\Gamma_H=C$ and $\Gamma\setminus C$ is connected.
Note that if this is so then $\Gamma$ is unsplittable iff the
result of deleting from $\Gamma\setminus C$ all edges
with color $h$ is a forest whenever $1\neq h\in H$.

\proc{Theorem \ThmDThreeMCirc} 
Let $d=3$ and let $\Gamma$ have an unsplittable coloring with 
a special $m$--circuit. Then $3\leq m\leq b_1(\Gamma)$,
there is a short exact sequence
$$0 \to \bigoplus_{H\in\indexset^\star}H_1(M_H)
\mapright\beta H_1(\cover M) \to 
\zmodb{4}^{m-3}\oplus\ztb^{2(b_1(\Gamma)-m)}\to 0,$$
and 
$\beta\bigl(\bigoplus_{H\in\indexset^\star}H_1(M_H)\bigr)
=4H_1(\cover M)$.
\endproc

This implies part (2) of 
\cite{\Sakuma, Theorem 14.2}, since $K_4$
has a unique 
$G(3)$--coloring, which is
unsplittable and
has a special 3--circuit. Once again,
when Theorem \ThmDThreeMCirc\ applies,
Proposition \Algebra\ shows that
$H_1(\cover M)$ is determined up to isomorphism by
$\bigoplus_{H\in\indexset^\star}H_1(M_H)$. We now show
that Theorem \ThmDThreeMCirc\ applies to infinitely many colored
graphs.

\proc{Proposition \label\MCircExists} Let $m$ and $b$ be integers
with $3\leq m\leq b$. Then there is 
a graph $\Gamma$ with $b_1(\Gamma)=b$ 
and an unsplittable $G(3)$--coloring of $\Gamma$ which has a special
$m$--circuit.
\endproc

\prf First we show that for $m\geq 3$ there is
a graph $\Gamma$ with $b_1(\Gamma)=m$ 
and an unsplittable $G(3)$--coloring of $\Gamma$ which has a special
$m$--circuit. Let $T$ be a tree with $m$ vertices
of valence 1 (its leaves) and $m-2$ vertices of valence 3
(its forks); such
trees exist for any $m\geq 2$. Form $\Gamma$ by adding
an $m$--circuit $C$ through the leaves of $T$.
Pick $H_0 \in \indexset^\star$ and $g_0\in G-H_0$. It is easy
to color the edges of $T$ with non-trivial elements of $H_0$
so that the required relation holds at each fork.
Further pick an edge $e_0$ and a vertex $v_0$ of $C$.
Give $e_0$ the color $g_0$. There is then a unique way to
color the other edges of $C$ so that the required relation
holds at every vertex except perhaps $v_0$. If we take the product
over all vertices $v$ of the product of the edge-colors at
$v$, the result is $1$, since each edge-color appears twice.
It follows that the required relation holds at $v_0$ as well.
Since $m\geq 3$, $T$ has at least one fork,
and so all non-trivial elements of $H_0$ are used to color $T$,
and the edge-colors of $\Gamma$ generate $G$. Also, all the colors
of $C$ are in $G-H_0$. It follows first that they are non-trivial,
so we do have a $G$--coloring, and second that $\Gamma_{H_0} = C$,
so that $C$ is a special $m$--circuit.
Since deleting edges from a tree always leaves a forest,
this coloring is unsplittable.

\fignum\FigWyeDelta
\theFigWyeDelta\FigWyeDelta

Now, if $\Gamma$ is any 
unsplittable $G(3)$--colored graph with a special
$m$--circuit, performing the operation of Figure \FigWyeDelta\ at
any vertex not on that circuit yields a graph
$\Gamma'$ which is unsplittable, has a special $m$--circuit,
and has $b_1(\Gamma')=b_1(\Gamma)+1$; the general case follows.
\qed

We give some specific examples of such colorings.

\rk{Example \label\ExMobDThree} Let $\Gamma$ be an $n$--rung M\"obius ladder 
($n\geq 2$). It is possible to determine all unsplittable
$G(3)$--colorings of $\Gamma$ with a special circuit; we shall
describe them but omit the verification that there are no others.
First, an $(n+1)$--circuit consisting of one rung together 
with half the rim has complementary graph a tree. By the first
part of the above proof, there is an unsplittable coloring for
which this circuit is special. Next, suppose that $n\geq 3$
and let $\{x_1,x_2,x_3\}$ be a basis of $G$. Color the rim edge
$\sigma_0$ with $x_1$, the rung $\tau_0$
with $x_2$, the rung $\tau_1$ with $x_2x_3^{n-1}$, and all
other rungs with $x_3$. There is a unique way to complete
the coloring, and there is a special 4--circuit corresponding
to the subgroup $\genby{x_1x_2, x_3}$; unsplittability is easily checked.
Finally, there is an exceptional coloring when $n=4$: color the
rung $\tau_0$ with $x_1x_2x_3$, $\tau_i$ with $x_i$ for $1\leq i\leq 3$,
and the rim edge $\sigma_0$ with $x_2$. This determines an unsplittable
coloring with a special 4--circuit corresponding to $\genby{x_1, x_2}$.
\endrk

\rk{Example \label\ExGenPetersen} 
In \cite\Watkins, the generalized Petersen graph $P(n,k)$
was defined for $1\leq k \leq n-1$ and $n\neq 2k$ as follows.
It has $2n$ vertices $u_0,\ldots u_{n-1},v_0,\ldots v_{n-1}$,
and edges of three kinds, namely $\sigma_i=\{u_i,u_{i+1}\}$,
$\tau_i=\{v_i,v_{i+k}\}$ and $\rho_i=\{u_i,v_i\}$,
where the subscripts are taken modulo $n$. The edges
$\sigma_i$ form an $n$--circuit (the outer rim);
if $k$ is coprime to $n$ (as we shall assume), 
so do the edges $\tau_i$
(the inner rim). The edges $\rho_i$ are called rungs.
Pick $H_0 \in \indexset^\star$ and $g_0\in G-H_0$.
Color the edges of the inner rim with non-trivial elements
of $H_0$ so that adjacent edges receive distinct colors and all
three elements appear. This forces colors on the rungs.
If one edge of the outer rim is given the
color $g_0$, there is a unique way to complete
the $G$--coloring. Then $\Gamma_{H_0}$ is the outer rim,
whose complementary graph is connected; it is easy to see
that this coloring is unsplittable. This example
does not arise from the construction of Proposition \MCircExists.

If $n=2m+1$ and $k=2$ there is also an unsplittable
coloring with a special $(n+1)$--circuit;
the complementary graph to the circuit
$u_1 u_2 \ldots u_{2m} v_{2m} v_1 u_1$ is a tree,
so there is an unsplittable coloring for
which this circuit is special.
\endrk

We have one other theorem in the case $d=3$.

\proc{Theorem \ThmDThreeMob} 
Let $\Gamma$ be an $n$--rung M\"obius ladder ($n\geq 2$) with
a $G(3)$-coloring, and let $g_0$ be the product of the colors on the rungs.
Suppose that $g_0\neq 1$, and let $k$ be the number of rungs with
color $g_0$. If $k=0$, there is a short exact sequence
$$0 \to \bigoplus_{H\in\indexset^\star}H_1(M_H)
\mapright\beta H_1(\cover M) \to 
\zmodb{4}^{n-2}\to 0,$$
while if $k>0$ 
there is a short exact sequence
$$0 \to \bigoplus_{H\in\indexset^\star}H_1(M_H)
\mapright\beta H_1(\cover M) \to 
\zmodb{4}^{n-k-1}\oplus\ztb^{2(k-1)}\to 0.$$
In either case,
$\beta\bigl(\bigoplus_{H\in\indexset^\star}H_1(M_H)\bigr)
=4H_1(\cover M)$.
\endproc

There is considerable overlap between Theorems \ThmDThreeMCirc\ and
\ThmDThreeMob; all the colorings of Example \ExMobDThree\ apart
from the exceptional coloring for $n=4$ satisfy the hypothesis
of Theorem \ThmDThreeMob. However, it is easy to see that there
are infinitely many colorings satisfying that hypothesis
that do not have a special circuit.

Next we consider some $G(4)$--colorings of M\"obius ladders.

\rk{Example \label\ExMobDFour} Let $d=4$,
and let $\{x_1,x_2,x_3,x_4\}$ be a basis of $G$. Let
$\Gamma$ be an $n$--rung M\"obius ladder with $n\geq 3$.
Give the colors $x_1$, $x_2$ and $x_1x_2x_3^n$ to one rung
each, and give all other rungs the color $x_3$. If we
give any rim edge the color $x_4$, there is then a unique
way to color the remaining edges with elements of $G$
so that the required relation holds at each vertex,
and this does give a $G$--coloring. 
Here every $\Gamma_H$ is connected;
this can be seen by listing all the $\Gamma_H$, but it is easier
to make use of the following lemma.
\endrk

\proc{Lemma \label\EdgeParity} Let $e_1,\ldots,e_n$ be distinct
edges of a $G$--colored graph $\Gamma$ with colors $g_1,\ldots,g_n$.
For $H\in\indexset^\star$, the number of these edges contained
in $\Gamma_H$ is even iff $g_1\cdots g_n\in H$.
\endproc

\prf Let $\achar_H$ be the homomorphism from $G$ to $\zt$
with kernel $H$, and let
$k$ of the edges $e_1,\ldots,e_n$ be contained in
$\Gamma_H$. Since $e_i$ is contained in $\Gamma_H$ iff
$\achar_H(g_i)=1$, $\achar_H(g_1\cdots g_n) = k\bmod 2$,
and the result follows.
\qed

For the colorings of Example \ExMobDFour, the product
of the colors on the rungs is $x_3$. Let $H\in \indexset^\star$.
If $x_3\notin H$  then $\Gamma_H$ contains an odd number of rungs 
by the lemma, while if $x_3\in H$ then $\Gamma_H$ contains at most 
three rungs; in either case, $\Gamma_H$ is connected.

\proc{Theorem \ThmDFour} Let $d=4$ and
let $\Gamma$ be an $n$--rung
M\"obius ladder with $n\geq 3$. Give $\Gamma$
the $G(4)$--coloring of Example \ExMobDFour. Then
$$H_1(\cover M) \cong \cases{
\bigoplus_{H\in\indexset^\star}H_1(M_H)\oplus \zt,
& if $n=3$;\cr
\bigoplus_{H\in\indexset^\star}H_1(M_H)
\oplus \zmod{8} \oplus \ztb^{4n-14},
&if $n\geq 4$.
}$$
\endproc

Our final theorem deals with a particular coloring
of the Petersen graph.

\rk{Example \label\ExPetersen} We use the notation of Example
\ExGenPetersen, and let $\Gamma$ be the Petersen graph $P(5,2)$.
Let $d=5$, and let $G$ have a basis $\{x_0,\ldots,x_4\}$.
Color the edge $\sigma_i$ with $x_i$,
the edge $\tau_i$ with $x_{i-1}x_{i+2}$,
and the edge $\rho_i$ with $x_{i-1}x_i$, all subscripts
being taken modulo 5. We leave it to the reader to check that
this is indeed a $G$--coloring. This graph has six disconnected
cycles, all of which contain an odd number of the edges $\tau_i$.
Since the product of the colors on the $\tau_i$ is $1$, it follows
from Lemma \EdgeParity\ that every $\Gamma_H$ is connected.
\endrk

\proc{Theorem \ThmDFive} Let $d=5$, and let $\Gamma$ 
be the Petersen graph with
the $G(5)$--coloring of Example \ExPetersen. Then
$$H_1(\cover M) \cong \bigoplus_{H\in\indexset^\star}H_1(M_H)
\oplus \zmod{16} \oplus \zmodb{4}^4 \oplus \ztb^{2}.$$
\endproc

The rest of this section sets out some notation. In the next section
we give the plan of the proof 
and explain the organization 
of the rest of the paper.

We deal often with direct sums $\bigoplus_{H\in\indexset'}\Lambda_H$,
where the $\Lambda_H$ are abelian groups 
indexed by a subset $\indexset'$
of $\indexset$. It is convenient to regard an element of
$\bigoplus_{H\in\indexset'}\Lambda_H$ as a formal linear combination
$\sum_{H\in\indexset'}\lambda_H H$ with $\lambda_H\in\Lambda_H$.
When all the $\Lambda_H$ are equal, we use
the notation $\Lambda^{\indexset'}$
for $\bigoplus_{H\in\indexset'}\Lambda$. As in the proof of
Lemma \EdgeParity, for $H\in\indexset$,
we let $\achar_H$ be the homomorphism $G\to \zt$ with kernel $H$;
we also let $\mchar_H$ the homomorphism with kernel 
$H$ from $G$ to the group
$\{\pm 1\}$ of units of $\z$ (a character of $G$).

If $X$ is a polyhedron, $C(X;\Lambda)$ will denote the
simplicial chain complex of some fixed but anonymous
triangulation of $X$, with coefficients in the abelian group
$\Lambda$. When the coefficient group is omitted, it is understood to 
be $\z$, except in \S7, where it is understood to be $\zt$.
We assume that the simplices of the triangulation have been oriented,
and by a simplex of $X$ we shall mean
a simplex of the triangulation with the chosen orientation;
thus the simplices of $X$ form a basis for $C(X)$. 
We let $S(X)$ be the set of all simplices of $X$, and 
$S_i(X)$ the subset of $i$--simplices.
If $\maps{f}{X}{Y}$
is a simplicial map, the induced maps on chain complexes and homology
will also be denoted by $f$ without further decoration. If
$f$ is a regular branched covering and the triangulation
of $X$ is obtained by lifting that of $Y$, we have the transfer
map $C(Y)\to C(X)$; recall that this sends a simplex $\sigma$
to $\sum_{k\in K}k\tilde\sigma$, where $K$ is the deck group and
$\tilde \sigma$ is one lift of $\sigma$. This map and the 
induced map on homology will both be denoted by $f\shriek$.
We let $b_i(X)$ be the $i$th Betti number of $X$.

\section{Outline of the proof}
Consider a regular branched covering
$\maps{\pi}{\cover M}{M}$ of a homology 3--sphere $M$,
with deck group $G$ and branch set a $G$--colored graph $\Gamma$.
Triangulate $M$ so that $\Gamma$ is triangulated by a subcomplex,
and lift this triangulation to triangulations of the $M_H$
and $\cover M$. We have various transfer maps
$\maps{\rho_H\shriek}{C(M)}{C(M_H)}$
and $\maps{\pi_H\shriek}{C(M_H)}{C(\cover M)}$. We define
chain maps
$$\eqalignno{
\alpha\colon C(M)^{\indexset^\star} &\to \bigoplus_{H\in\indexset}C(M_H)\cr
{\rm by}\quad \alpha\Bigl(\sum_{H\in\indexset^\star}c_H H\Bigr)
&= \sum_{H\in\indexset^\star}\left(\rho_H\shriek(c_H)H-c_H G\right)
\quad\hbox{for $c_H \in C(M)$, $H\in\indexset^\star$,}\cr
\noalign{\hbox{and}}
\beta\colon \bigoplus_{H\in\indexset}C(M_H) &\to C(\cover M)\cr
{\rm by}\quad \beta\Bigl(\sum_{H\in\indexset}d_H H\Bigr)
&= \sum_{H\in\indexset}\pi_H\shriek(d_H)
\quad\hbox{for $d_H \in C(M_H)$, $H\in\indexset$.}
}$$
We also let $\maps{\gamma}{C(\cover M)}{C(M;\zmod{2^{d-1}})}$
be the composite
of $\maps{\pi}{C(\cover M)}{C(M)}$ and
reduction of the coefficients modulo $2^{d-1}$. 

Consider the sequence
$$0 \to C(M)^{\indexset^\star}
\mapright\alpha \bigoplus_{H\in\indexset}C(M_H)
\mapright\beta C(\cover M) \mapright\gamma C(M;\zmod{2^{d-1}}) 
\to 0. \eqno{(\label\MainSeq)}$$
This is not exact, but we do have the following result.
\proclaim{Lemma} The chain map $\alpha$ is injective,
$\beta\alpha = 0$, $\gamma\beta = 0$, and $\gamma$ is surjective.
\endproc

\prf 
Define 
$\maps{\alpha'}{\bigoplus_{H\in\indexset}C(M_H)}
{C(M)^{\indexset^\star}}$
by
$$\alpha'\Bigl(\sum_{H\in\indexset}d_H H\Bigr) 
= \sum_{H\in\indexset^\star} \rho_H(d_H)H.$$
Then
$$\alpha'\alpha\Bigl(\sum_{H\in\indexset^\star}c_H H\Bigr)
= \sum_{H\in\indexset^\star}\rho_H\rho_H\shriek(c_H)H
=2 \sum_{H\in\indexset^\star}c_H H,$$
so $\alpha$ is injective. Next,
$$\beta\alpha\Bigl(\sum_{H\in\indexset^\star}c_H H\Bigr)
= \sum_{H\in\indexset^\star}\left(\pi_H\shriek\rho_H\shriek(c_H)
-\pi\shriek(c_H)\right) = 0,$$
so $\beta\alpha=0$. Further,
$$\pi\beta\Bigl(\sum_{H\in\indexset}d_H H\Bigr)
= \sum_{H\in\indexset}\rho_H\pi_H\pi_H\shriek(d_H)
= \sum_{H\in\indexset}|H| \rho_H(d_H),$$
so $\gamma\beta = 0$.
Finally, $\maps{\pi}{C(\cover M)}{C(M)}$ is clearly onto,
and hence so is $\gamma$.
\qed

The sequence (\MainSeq) thus decomposes into four short
exact sequences:

$$\eqalignno{
&0\to C(M)^{\indexset^\star} \mapright\alpha \ker\beta 
\to \ker\beta/\im\alpha \to 0; &(\label\ShortSeqOne)\cr
&0\to \ker\beta \hookmapright\iota \bigoplus_{H\in \indexset}C(M_H) 
\mapright\beta \im\beta \to 0; &(\label\ShortSeqTwo)\cr
&0\to \im\beta \hookrightarrow \ker\gamma
\to \ker\gamma/\im\beta \to 0; &(\label\ShortSeqThree)\cr
{\rm and}\quad &0 \to \ker\gamma\hookrightarrow C(\cover M)
\mapright\gamma C(M;\zmod{2^{d-1}}) \to 0. &(\label\ShortSeqFour)
}$$
The last of these  relates the homology groups of $\cover M$
and the complex
$\ker\gamma$; the first homology is all we need.

\proclaim{Lemma} We have 
$H_1(\ker\gamma)\cong H_1(\cover M)$.
\endproc
\prf Since $M$ is an integral homology sphere, it is also a
$\zmod{2^{d-1}}$ homology sphere,
so part of the long exact sequence of (\ShortSeqFour) 
is
$0 \to H_1(\ker\gamma) \to H_1(\cover M) \to 0$.
\qed

To extract information from the exact sequences 
(\ShortSeqOne)--(\ShortSeqThree), we need to study the complexes
$\ker\beta/\im\alpha$ and $\ker\gamma/\im\beta$. This
leads us to consider certain chain complexes associated
to a $G$--colored graph $\Gamma$. These chain complexes are defined
and studied in \S4, after some preliminary results on the graded
ring of $G$ in \S3.  
In \S5, we determine the complex $\ker\beta/\im\alpha$,
and in \S6, we determine
the quotients of a filtration of $\ker\gamma/\im\beta$.
In \S7 we prove some results
on the $\zt$ homology of 2-- and 4--fold branched covers,
and in \S8 we prove our theorems.

\section{The graded ring of $G$}

As always, $G$ is a group isomorphic to $\ztb^d$, but in this
section we do not assume that $d\geq 2$.
For $H\in\indexset$, the character $\mchar_H$ extends to a ring 
homomorphism $\maps{\mchar_H}{\z[G]}{\z}$ on the group ring
of $G$.
The fundamental ideal $I=I[G]$ of $G$ is the kernel of
$\mchar_G$; we also let $J=J(G)$ be the ideal of those
$\lambda\in \z[G]$ for which $\mchar_G(\lambda)\equiv 0 \pmod 2$.
Note that $J = I\oplus 2\z$.
We consider the associated graded rings $A=A(G) = \graded_I(\z[G])$
and $B=B(G)=\graded_J(\z[G])$. (See \cite{\ZS, p.248}.)
Consider first the ring $B$. The group of homogeneous elements
of degree $k$ is $B_k = J^k/J^{k+1}$, and $B$
is an algebra over $B_0 = \z[G]/J$,
which we identify with $\zt$. We denote the image in $B_k$
of $\lambda\in J^k$ by $[\lambda]_k$; the product is given by
$[\lambda]_k[\mu]_l = [\lambda\mu]_{k+l}$. Turning to $A$,
we have $(1-g)^2=2(1-g)$ for $g\in G$, so $2I\leq I^2$.
Let $k\geq 1$.
It follows that $2I^k\leq I^{k+1}$, and hence
$J^k = I^k \oplus 2^k\z$. Therefore we may identify $A_k$ with its image in $B_k$, and $B_k$  is the direct sum of $A_k$ and a copy
of $\zt$ generated by $[2^k]_k$. Note also that $A_kB_l=A_{k+l}$.
Of course $A_0\cong\z$; below, when we refer to $A_k$,
it is to be understood that $k\geq 1$.

We shall determine the structure of the algebra $B$, 
and hence that of $A$. (The structure of $\graded_I(\z[G])$
when $G$ is free abelian was determined by Massey in 
\cite\Massey.)
We define a function $\maps{\omega}{G}{A_1}$ by
$\omega(g)=[1-g]_1$.

\proc{Lemma \label\AOneEqG} The function 
$\omega$ is an isomorphism, and for any
$\lambda=\sum_{g\in G}\lambda_gg \in I$ we have
$\omega\bigl(\prod_{g\in G}g^{\lambda_g}\bigr) = [\lambda]_1$.
\endproc
\prf We compute
$$(1-g)+(1-h)-(1-gh)= 1-g-h+gh =(1-g)(1-h)\in I^2,$$
so
$[1-g]_1+[1-h]_1=[1-gh]_1$, and $\omega$ is a homomorphism.
The function $[\lambda]_1 \mapsto \prod_{g\in G}g^{\lambda_g}$
is a well-defined homomorphism $A_1 \to G$ sending
$\omega(g)$ to $g$. For $\lambda \in I$ we have
$\lambda = -\sum_{g\in G}\lambda_g(1-g)$,
so $[\lambda]_1 = \sum_{g\in G}\lambda_g\omega(g)
= \omega\bigl(\prod_{g\in G}g^{\lambda_g}\bigr)$,
and the result follows.
\qed

For $0\leq l\leq d$, we let $\scr I_l$ be the set
of $l$-tuples $\veci = (i_1,i_2,\ldots,i_l)$
of integers with $1\leq i_1<i_2<\cdots <i_l\leq d$.

\proc{Lemma \label\JKBasis} Let $\{x_1,\ldots,x_d\}$ be a basis
of $G$, and (for $k\geq 0$) let $\scr B_k$ be the set consisting
of the elements 
$(1-x_{i_1})\cdots (1-x_{i_l})$
for $k\leq l \leq d$ and $\veci \in \scr I_l$,
together with the elements 
$2^{k-l}(1-x_{i_1})\cdots (1-x_{i_l})$
for $0\leq l<k$, $l\leq d$ and $\veci \in \scr I_l$.
(When $l=0$, the empty product $(1-x_{i_1})\cdots (1-x_{i_l})$
is taken to be $1$.)
Then $\scr B_k$ is a basis of $J^k$ (as a $\z$--module).
Further, an element $\lambda$
of $\z[G]$ is in $J^k$ iff
$\mchar_H(\lambda) \equiv 0\pmod{2^k}$ for all 
$H\in \indexset$.
\endproc

\prf Every element $g$ of $G$ can be written uniquely
in the form $g=x_{i_1}\cdots x_{i_l}$ for $0\leq l \leq d$
and $\veci\in I_l$; call $l$ the length of $g$.
Then $g$ is the unique element of maximal length appearing
in 
$(1-x_{i_1}) \cdots (1-x_{i_l})$, 
and it follows that the 
$(1-x_{i_1}) \cdots (1-x_{i_l})$ 
are linearly independent;
therefore so are the elements of $\scr B_k$. Let
$V_k$ be the additive subgroup of $\z[G]$ spanned by $\scr B_k$,
and let $W_k$ be the subgroup of those $\lambda\in \z[G]$
such that $\mchar_H(\lambda) \equiv 0\pmod{2^k}$ for all
$H\in \indexset$.
Clearly $V_k\leq J^k$. 
Since $\mchar_H(\lambda)\equiv \mchar_G(\lambda) \pmod 2$,
we have $J= W_1$, and it follows that $J^k\leq W_k$ for all $k\geq 0$.

It remains to prove that $W_k\leq V_k$.
Let $\lambda=\sum_{g\in G}\lambda_g g$ be a non-zero
element of $\z[G]$. Let $l$
be the maximum length of those $g$ with $\lambda_g\neq 0$,
and let $n$ be the number of those $g$ of length $l$
with $\lambda_g\neq 0$. Call the pair $(l, n)$ the
weight of $\lambda$, and order weights lexicographically.
Suppose that $W_k\not\leq V_k$, and let $\lambda$ be an
element of $W_k-V_k$ of minimum weight $(l,n)$.
Let $h = x_{j_1}\cdots x_{j_l}$ ($\vecj\in\scr I_l$)
have $\lambda_h\neq 0$.
If $l \geq k$, then 
$\lambda - (-1)^l\lambda_h (1-x_{j_1}) \cdots (1-x_{j_l})$
is an element of $W_k-V_k$ of smaller weight than $\lambda$,
a contradiction.
Suppose that $l < k$. Let $G'$ be the subgroup of $G$ generated
by $x_{j_1},\ldots,x_{j_l}$, and $G''$ the subgroup
generated by the other $x_i$, so $G=G'\oplus G''$.
Since $\lambda\in W_k$, 
$$
\sum_{H'\in\indexset(G')}\mchar_{H'}(h)\mchar_{H'\oplus G''}(\lambda)
\equiv 0 \pmod{2^k}.
$$
Now 
$$
\sum_{H'\in\indexset(G')}\mchar_{H'}(h)\mchar_{H'\oplus G''}(\lambda)
= \sum_{g\in G}\Bigl(\sum_{H'\in\indexset(G')}
\mchar_{H'}(h)\mchar_{H'\oplus G''}(g)\Bigr)\lambda_g.
$$
Let $g=g'g''$, with $g'\in G'$ and $g''\in G''$. Then 
$\mchar_{H'}(h)\mchar_{H'\oplus G''}(g)=\mchar_{H'}(hg')$,
and $\sum_{H'\in\indexset(G')}\mchar_{H'}(hg')$ is $0$
if $g'\neq h$, and $2^l$ if $g'=h$. If $g'=h$ and $g''\neq 1$,
then $\lambda_g=0$ by our choice of $h$.
It follows that $2^l\lambda_h \equiv 0 \pmod{2^k}$, or
$\lambda_h\equiv 0 \pmod{2^{k-l}}$.
Now $\lambda - (-1)^l\bigl(\lambda_h/2^{k-l}\bigr)
2^{k-l}(1-x_{j_1}) \cdots (1-x_{j_l})$
is an element of $W_k-V_k$ of smaller weight than $\lambda$,
and this contradiction proves that $W_k\leq V_k$.
\qed

As an immediate consequence of this lemma, we have bases
for $A_k$ and $B_k$.
\proc{Lemma \label\BKBasis} Let $\{x_1,\ldots,x_d\}$ be a basis
of $G$. The elements 
$$[2^{k-l}(1-x_{i_1}) \cdots (1-x_{i_l})]_k
=[2]_1^{k-l}\omega(x_{i_1})\cdots\omega(x_{i_l})$$
for $0\leq l\leq\min\,\{k,d\}$ and $\veci \in \scr I_l$
form a basis of $B_k$ (as a $\zt$ vector space),
and those for $1\leq l\leq\min\,\{k,d\}$ form a basis
for $A_k$.
\qed

Note that this implies that multiplication by $[2]_1$
defines injections $B_k \to B_{k+1}$ for $k\geq 0$
and $A_k\to A_{k+1}$ for $k\geq 1$, and that these
are onto for $k\geq d$.

\proc{Lemma \label\SymmAlg} The graded algebra $B$
is the quotient of the symmetric algebra of $B_1$
by the relations $a^2 = [2]_1a$ for $a \in A_1$.
\endproc

\prf The given relations do hold in $B$: by Lemma \AOneEqG,
any element of $A_1$ equals $\omega(g)$ for some $g\in G$,
and $\omega(g)^2 = [(1-g)^2]_2 = [2(1-g)]_2 = [2]_1\omega(g)$.
Therefore, if $\hat B$  is 
the quotient of the symmetric algebra of $B_1$ by these relations,
there is an epimorphism $\hat B \to B$. But 
if $\{x_1,\ldots,x_d\}$ is a basis
of $G$, then $\hat B_k$
is generated by the elements
$[2]_1^{k-l}\omega(x_{i_1})\cdots\omega(x_{i_l})$
for $0\leq l\leq\min\,\{k,d\}$ and $\veci \in \scr I_l$,
and these map to independent elements in $B_k$ by Lemma
\BKBasis.
\qed

The $\zt$ vector space $\ztb^{\indexset^\star}$ is a commutative 
algebra under componentwise multiplication. Its identity
element $\sum_{H\in\indexset^\star}H$ will be denoted
by $1^{\indexset^\star}$.
We may define a linear map $\maps\Omega{B_1}{\ztb^{\indexset^\star}}$
by $\Omega(\omega(g)) = \sum_{H\in\indexset^\star}\achar_H(g)H$
for $g\in G$, and $\Omega([2]_1) = 1^{\indexset^\star}$.
Since $x^2=x$ for all $x \in \ztb^{\indexset^\star}$,
it follows from Lemma \SymmAlg\ that $\Omega$ extends (uniquely)
to an algebra homomorphism $\maps\Omega{B}{\ztb^{\indexset^\star}}$.

\proc{Lemma \label\OmegaInj} The map $\Omega$ restricts
to an injection on $A_k$ for $1\leq k$,
and on $B_k$ for $0\leq k \leq d-1$. Further,
$\Omega$ maps each of $A_d$ and $B_{d-1}$ onto
$\ztb^{\indexset^\star}$.
\endproc

\prf We show first that $\Omega$ maps $A_d$ isomorphically
onto $\ztb^{\indexset^\star}$. For any $g_1,\ldots,g_d\in G$,
we have 
$\Omega(\omega(g_1)\cdots\omega(g_d))
= \sum_{H\in\indexset^\star}\achar_H(g_1)\cdots\achar_H(g_d)H$.
Given $H_0\in\indexset^\star$ we may find a basis
$\{x_1,\ldots,x_d\}$ of $G$ with $x_i\notin H_0$
for $1\leq i\leq d$. Then
$\achar_{H_0}(x_1)\cdots\achar_{H_0}(x_d) =1$ and
$\achar_H(x_1)\cdots\achar_H(x_d) = 0$ for any $H\neq H_0$,
so $\Omega(\omega(x_1)\cdots\omega(x_d)) = H_0$.
Thus $\Omega$ maps $A_d$ onto $\ztb^{\indexset^\star}$.
Since $\dim A_d = 2^d-1=\dim \ztb^{\indexset^\star}$,
$\Omega$ is also injective on $A_d$.

Next, let $\{x_1,\ldots,x_d\}$ be any basis of $G$,
and consider the basis elements
$b_{\veci} = [2]_1^{d-l}\omega(x_{i_1})\cdots\omega(x_{i_l})$
($0\leq l\leq d$, $\veci \in \scr I_l$) of $B_d$.
Let $s$ be the sum of all $\Omega(b_{\veci})$.
For each $H\in\indexset^\star$, the coefficient
of $H$ in $\Omega(b_{\veci})$ is 1 if
$x_{i_1},\ldots,x_{i_l}\notin H$, and 0 otherwise.
Therefore the coefficient of $H$ in $s$ is the number
of subsets of $\{x_1,\ldots,x_d\}\cap (G-H)$, taken modulo
2. Since $\{x_1,\ldots,x_d\}\cap (G-H)$ is non-empty,
this number is even, so $s=0$. It follows that $\Omega$
maps the subspace of $B_d$ spanned by the $b_{\veci}$
for $l<d$ isomorphically onto $\ztb^{\indexset^\star}$.
Since multiplication by $[2]_1$ maps $B_{d-1}$
isomorphically onto this space and
$\Omega([2]_1b) = \Omega(b)$ for all $b\in B$,
$\Omega$ also maps $B_{d-1}$ isomorphically onto
$\ztb^{\indexset^\star}$.
Since multiplication by $[2]_1$ maps $B_k$ injectively into
$B_{k+1}$, it follows that $\Omega$ is injective on $B_k$
for $0\leq k \leq d-1$, and therefore on $A_k$ for
$1\leq k \leq d-1$. Finally, multiplication by $[2]_1$ maps $A_k$ 
isomorphically onto $A_{k+1}$ for $k\geq d$, and hence $\Omega$
is injective on $A_k$ for $k\geq d$.
\qed

There is an inner product on $\ztb^{\indexset^\star}$ given
by 
$$\bigl(\sum_{H\in\indexset^\star}a_H H\bigr) \cdot
\bigl(\sum_{H\in\indexset^\star}b_H H\bigr) = 
\sum_{H\in\indexset^\star}a_H b_H.$$
Note that for any $x$ and $y$ in $\ztb^{\indexset^\star}$,
$x\cdot y = 1^{\indexset^\star}\cdot(xy)$.

\proc{Lemma \label\Orthog} For $1\leq k \leq d-1$, we
have $\Omega(A_k) = \Omega(B_{d-k-1})^\perp$,
where $^\perp$ denotes the orthogonal complement
with respect to the above inner product.
\endproc
\prf Since $\dim A_k = \sum_{i=1}^k{d\choose i}$ and
$\dim B_{d-k-1} = \sum_{i=0}^{d-k-1}{d\choose i}
= \sum_{i=k+1}^d{d\choose i}$, we have
$\dim A_k + \dim B_{d-k-1} = 2^d-1 = \dim \ztb^{\indexset^\star}$.
Therefore it suffices to prove that $\Omega(a)\cdot\Omega(b) = 0$
for $a \in A_k$ and $b\in B_{d-k-1}$. 
Since $\Omega(a)\cdot\Omega(b) 
= 1^{\indexset^\star}\cdot\Omega(ab)$
and $ab \in A_{d-1}$, it is enough to show
that $1^{\indexset^\star}\cdot\Omega(A_{d-1}) = 0$. For $g_1,\ldots,g_{d-1}\in G$,
$1^{\indexset^\star}\cdot\Omega(\omega(g_1)\cdots\omega(g_{d-1}))$ 
is the
number of $H\in\indexset^\star$ that contain none of
$g_1,\ldots,g_{d-1}$, taken modulo 2. Since 
$g_1,\ldots,g_{d-1}$ do not generate $G$, this number
is even, and we are done.
\qed

Now let $G'$ be a subgroup of $G$, and set $G''=G/G'$.
We have an epimorphism $\z[G]\to\z[G'']$ inducing
epimorphisms $I[G]^k \to I[G'']^k$ and
$A_k(G) \to A_k(G'')$ for all $k\geq 1$.
We denote the kernels of these maps by
$\z[G,G']$, $I^k[G,G']$, and $A_k(G,G')$. 
For $k=1$, $A_1(G,G')$
is just the image of $G'$ under the isomorphism
$\maps{\omega}{G}{A_1(G)}$. 

\proc{Lemma \label\OmegaAKRel} Let $G'\leq G$
and $a\in A_k(G)$ ($1\leq k \leq d$). Let
$\Omega(a) = \sum_{H\in\indexset^\star(G)} a_H H$.
Then $a\in A_k(G,G')$ iff $a_H=0$ whenever $H\geq G'$.
\endproc
\prf Let $G''=G/G'$.
There is a linear map 
$\ztb^{\indexset^\star(G)}\to \ztb^{\indexset^\star(G'')}$
sending $H\in \indexset^\star(G)$ to $H/G'$ if $H\geq G'$,
and to zero otherwise; its kernel consists of all
$\sum_{H\in\indexset^\star(G)} a_H H$
such that $a_H=0$ whenever $H\geq G'$. We also have the algebra 
homomorphism $\maps{\Omega''}{B(G'')}{\ztb^{\indexset^\star(G'')}}$.
Restricting to $A_k(G)$, we have a commutative diagram
{\abovedisplayskip=\diagdisplayskip  
\belowdisplayskip=\diagdisplayskip   
$$
\def\normalbaselines{\diagbaselines} 
\matrix{
A_k(G) & \mapright{\Omega} & \ztb^{\indexset^\star(G)}\cr
\mapdown{} && \mapdown{}\cr
A_k(G'') & \mapright{\Omega''} & \ztb^{\indexset^\star(G'')}\cr
}$$}%
By Lemma \OmegaInj, $\Omega''$ is injective, and the result follows.
\qed

\section{Homology groups of colored graphs}

Let $\Gamma$ be a $G(d)$--colored graph ($d\geq 2$),
and fix a triangulation of $\Gamma$.
In this section we study chain complexes
$C(\Gamma\sep k)$
(for $k=1,2,\ldots d$)
of $\zt$ vector spaces associated to
this triangulation. (Why we should want to do this
will emerge in later sections.) Let $\sigma$ be a simplex of
$\Gamma$. 
If $\sigma$ is a vertex of $\Gamma$, we let $G_\sigma$
be the subgroup of $G$ generated by the colors of
the edges of $\Gamma$ incident to $\sigma$, which is isomorphic
to $\zt\oplus\zt$.
If $\sigma$ is any other simplex of $\Gamma$, then $\sigma$
is contained in a unique edge of $\Gamma$,
and we let $g_\sigma$ be the color of this edge,
and $G_\sigma$ the subgroup of $G$ generated by 
$g_\sigma$. We also set $A_k^\sigma = A_k(G, G_\sigma)$.
We let $C(\Gamma\sep k)$ be the subcomplex 
of $C(\Gamma; A_k)$
generated by all chains of the form $a\sigma$
where $\sigma$ is a simplex of $\Gamma$ and 
$a\in A_k^\sigma$.
(This is a subcomplex because if $\tau$ is a face of $\sigma$
then $A_k^\sigma\leq A_k^\tau$.) We let $b_i(\Gamma\sep k)$
be the dimension of the $i^{\rm th}$ homology group $H_i(\Gamma\sep k)$
of $C(\Gamma\sep k)$. Of course, the homology groups
are zero except in dimensions 0 and 1, and $H_1(\Gamma\sep k)$
is equal to the space $Z_1(\Gamma\sep k)$ of 1--cycles.
We let $\chi(\Gamma\sep k)= b_0(\Gamma\sep k)-b_1(\Gamma\sep k)$
be the $\zt$ Euler characteristic of $C(\Gamma\sep k)$.
It is clear that the homology of $C(\Gamma\sep k)$ is unchanged by
subdivision, and therefore independent of the triangulation.

\proc{Lemma \label\ChiK} We have 
$\chi(\Gamma\sep k) = -{d-2 \choose k-1}\chi(\Gamma)$.
(In the case $k= d$ we are using the convention that
${n\choose r}= 0$ for $r>n$.)
\endproc
\prf By Lemma \BKBasis, the dimension of $A_k^\sigma$
is $a = \sum_{i=1}^k \bigl({d\choose i}-{d-2\choose i}\bigr)$
if $\sigma$ is a vertex of $\Gamma$, and 
$b = \sum_{i=1}^k \bigl({d\choose i}-{d-1\choose i}\bigr)$
otherwise. Therefore
$\chi(\Gamma\sep k) = b\chi(\Gamma) +(a-b)V$,
where $V$ is the number of vertices of $\Gamma$.
Since $\Gamma$ is
trivalent, $V = -2\chi(\Gamma)$, so
$\chi(\Gamma\sep k) = -(2a-3b)\chi(\Gamma)$,
and it is easy to compute that $2a-3b={d-2 \choose k-1}$.
\qed

\proc{Lemma \label\HKOne} We have
$b_0(\Gamma\sep 1)=b_1(\Gamma)$ and
$b_1(\Gamma\sep 1)=b_0(\Gamma)$.
\endproc

\prf Let $a = \sum_{\sigma\in S_1(\Gamma)} a_\sigma\sigma$ 
($a_\sigma \in A_1^\sigma$)
be a 1--chain of $C(\Gamma\sep 1)$. For each 1--simplex
$\sigma$ of $\Gamma$, $A_1^\sigma \cong G_\sigma \cong \zt$,
with non-trivial element $\omega(g_\sigma)$.
If $g_1$, $g_2$, and $g_3$ are the colors of three
edges meeting at a vertex,
$\omega(g_1)+\omega(g_2)+\omega(g_3)
=\omega(g_1g_2g_3)=0$. It follows that 
$a$ is a cycle iff, for each component
$\Gamma'$ of $\Gamma$, the $a_\sigma$ for 1--simplices
$\sigma$ of $\Gamma'$ are either all zero or all
non-zero. This proves that $H_1(\Gamma\sep 1)\cong\ztb^{b_0(\Gamma)}$,
or $b_1(\Gamma\sep 1)=b_0(\Gamma)$, and it then follows from
Lemma \ChiK\ that
$b_0(\Gamma\sep 1)=b_1(\Gamma)$.
\qed

\proclaim{Lemma} For $1\leq k\leq d$,
there is an injection 
$\maps{\iota_k}{B_{k-1}}{Z_1(\Gamma\sep k)}$ defined by
$\iota_k(b) = \sum_{\sigma\in S_1(\Gamma)}\omega(g_\sigma)b\sigma$.
\endproc

\prf It is clear that the given formula defines a linear map
from $B_{k-1}$ to $C_1(\Gamma\sep k)$.
Let $b\in B_{k-1}$, and let $\tau$ be a 0--simplex
of $\Gamma$. If $\tau$ is not a vertex of $\Gamma$, it
is clear that the coefficient of $\tau$ in $\partial \iota_k(b)$
is zero. If $\tau$ is a vertex and the adjacent edge-colors
are $g_1$, $g_2$ and $g_3$, this coefficient is
$\sum_{i=1}^3 \omega(g_i)b = \omega(g_1g_2g_3)b = 0$.
Thus $\iota_k(b)$ is a cycle. It remains to show that
$\iota_k$ is injective; suppose that $\iota_k(b)=0$.
Then $\omega(g_\sigma)b = 0$ for every 1--simplex $\sigma$ of $\Gamma$.
Since the $g_\sigma$ generate $G$, this implies that
$ab=0$ for every $a\in A_1$, and therefore for every
$a \in A_d$. Now $\Omega(a)\cdot\Omega(b) 
= 1^{\indexset^\star}\cdot\Omega(ab) = 0$. Since $\Omega$
maps $A_d$ onto $\ztb^{\indexset^\star}$, it follows that
$\Omega(b)=0$; since $\Omega$ is injective on
$B_{k-1}$, we have $b=0$.
\qed

We call $\Gamma$ $k$--taut if $\iota_k$ is an isomorphism;
by Lemma \BKBasis, this occurs iff
$b_1(\Gamma\sep k)=\sum_{i=0}^{k-1}{d \choose i}$.
By Lemma \HKOne, $\Gamma$
is 1--taut iff it is connected. To give examples of $k$--taut
graphs for $k>1$, we use a different description of the 
chain complex $C(\Gamma\sep k)$. For $1\leq k\leq d$,
the injection $\maps\Omega{A_k}{\ztb^{\indexset^\star}}$
induces an injection
$\maps\Omega{C(\Gamma;A_k)}{C(\Gamma;\ztb^{\indexset^\star})}$,
which is onto for $k=d$.
We identify $C(\Gamma;\ztb^{\indexset^\star})$
with $C(\Gamma;\zt)^{\indexset^\star}$. Since
$\Omega([2]_1a) = \Omega(a)$ and $[2]_1A_k\leq A_{k+1}$,
we have a chain of subcomplexes
$$\Omega C(\Gamma;A_1)\leq \Omega C(\Gamma;A_2)\leq \cdots
\leq\Omega C(\Gamma;A_d) = C(\Gamma;\zt)^{\indexset^\star}.$$
A chain of $C(\Gamma;\zt)^{\indexset^\star}$ is of the form
$\sum_{\sigma\in S(\Gamma),H\in\indexset^\star} 
a_{\sigma, H}\sigma H$
with coefficients
$a_{\sigma, H}$ in $\zt$. It belongs
to $\Omega C(\Gamma;A_k)$ iff, for each simplex $\sigma$,
$\sum_{H\in\indexset^\star}a_{\sigma, H}H \in \Omega(A_k)$.
We let $C'(\Gamma \sep k)$ be the subcomplex $\Omega C(\Gamma \sep k)$
of $\Omega C(\Gamma;A_k)$. For $1\leq k\leq d-1$ and $a\in A_k$,
we have $a\in A_k^\sigma$ iff $[2]_1a\in A_{k+1}^\sigma$;
it follows that $C'(\Gamma\sep k) = \Omega C(\Gamma;A_k)\cap
C'(\Gamma\sep k+1)$. By Lemma \OmegaAKRel, a chain
$\sum_{\sigma\in S(\Gamma),H\in\indexset^\star} 
a_{\sigma, H}\sigma H\in C(\Gamma;\zt)^{\indexset^\star}$
belongs to $C'(\Gamma\sep d)$ iff
$a_{\sigma,H}=0$ whenever $H\geq G_\sigma$.
Now $H\geq G_\sigma$ iff $\sigma$
is not a simplex of $\Gamma_H$, so we may identify
$C'(\Gamma\sep d)$ with 
$\bigoplus_{H\in\indexset^\star}C(\Gamma_H;\zt)$.
It follows that a 1--chain 
$\sum_{\sigma,H} a_{\sigma, H}\sigma H$ 
of $C'(\Gamma\sep k)$ is a cycle iff, for each $H\in\indexset^\star$,
$a_{\sigma,H}$ is constant on each component of $\Gamma_H$.
We let $W(\Gamma\sep k)$ be the subspace of $Z'_1(\Gamma\sep k)$
consisting of all 1--chains  of $C'(\Gamma\sep k)$
such that, for each $H\in\indexset^\star$,
$a_{\sigma,H}$ is constant on all of $\Gamma_H$.

\proclaim{Lemma}
For $1\leq k\leq d$, $W(\Gamma\sep k) = \Omega\iota_k(B_{k-1})$.
\endproc

\prf We first prove the case $k = d$. $\Omega$ maps $B_{d-1}$
isomorphically onto $\ztb^{\indexset^\star}$, and
there is an isomorphism $\ztb^{\indexset^\star}\to W(\Gamma\sep d)$ 
sending
$\sum_{H\in\indexset^\star}b_H H$ to 
$\sum_{H\in\indexset^\star}\sum_{\sigma\in S_1(\Gamma_H)}b_H\sigma H$.
We show that the composite is equal to $\Omega\iota_d$. 
If $b\in B_{d-1}$ and $\Omega(b)= \sum_{H\in\indexset^\star}b_H H$
then
$$\Omega\iota_d(b) 
= \sum_{\sigma\in S_1(\Gamma)}\Omega(\omega(g_\sigma)b)\sigma
= \sum_{\sigma\in S_1(\Gamma),H\in\indexset^\star}
\achar_H(g_\sigma)b_H\sigma H
= \sum_{H\in\indexset^\star}\sum_{\sigma\in S_1(\Gamma_H)}b_H\sigma H,
$$
and this case is proved. 

Now let $k< d$. Since 
$W(\Gamma\sep k) = W(\Gamma\sep d) \cap C'(\Gamma\sep k)$,
it is enough to prove that 
$\Omega\iota_k(B_{k-1}) = \Omega\iota_d(B_{d-1})\cap C'(\Gamma\sep k)$.
Suppose that $b_k\in B_{k-1}$ and $b_d\in B_{d-1}$ are such
that $\Omega(b_k)=\Omega(b_d)$. Then
$$\Omega\iota_k(b_k) 
= \sum_{\sigma\in S_1(\Gamma)}\Omega(\omega(g_\sigma))\Omega(b_k)\sigma
= \sum_{\sigma\in S_1(\Gamma)}\Omega(\omega(g_\sigma))\Omega(b_d)\sigma
= \Omega\iota_{d}(b_d).
$$
Since, for any $b\in B_{k-1}$, $[2]_1^{d-k}b\in B_{d-1}$
and $\Omega([2]_1^{d-k}b)=\Omega(b)$, it follows that
$\Omega\iota_k(B_{k-1})$ is contained in 
$\Omega\iota_d(B_{d-1})\cap C'(\Gamma\sep k)$.
Conversely, for $b\in B_{d-1}$, we have $\Omega\iota_d(b) \in
C'(\Gamma\sep k)$ iff, for each $\sigma\in S_1(\Gamma)$,
$\Omega(\omega(g_\sigma)b)\in \Omega(A_k) = \Omega(B_{d-k-1})^\perp$
(using Lemma \Orthog). For $b'\in B_{d-k-1}$, 
$\Omega(\omega(g_\sigma)b)\cdot\Omega(b') = 
\Omega(b)\cdot\Omega(\omega(g_\sigma)b')$. Since the $g_\sigma$
generate $G$, the elements $\omega(g_\sigma)b'$ generate
$A_1B_{d-k-1}=A_{d-k}$. Therefore $\Omega\iota_d(b) \in
C'(\Gamma\sep k)$ iff $\Omega(b)\in\Omega(A_{d-k})^\perp
=\Omega(B_{k-1})$, and the proof is complete.
\qed

Thus $\Gamma$ is $k$--taut iff $W(\Gamma\sep k) = Z'_1(\Gamma\sep k)$.
Since $W(\Gamma\sep k) = W(\Gamma\sep k+1)\cap C'(\Gamma\sep k)$
for $k<d$, we have:
\proc{Lemma \label\TautDescends} For $1\leq k<d$, if $\Gamma$ is $(k+1)$--taut
then it is $k$--taut.
\qed

Since $C'(\Gamma\sep d) =
\bigoplus_{H\in\indexset^\star}C(\Gamma_H;\zt)$, we have:

\proclaim{Lemma} A $G(d)$--colored graph $\Gamma$ is $d$--taut
iff $\Gamma_H$ is connected for all $H\in\indexset^\star$.
\qed

Thus the $G(d)$--colored graphs of Theorems
\ThmDFour\ and \ThmDFive\ are $d$--taut.
Those of Theorems \ThmDTwo, \ThmDThreeMCirc\ and 
\ThmDThreeMob\ are not, in general, but they are 
$(d-1)$--taut. For Theorem \ThmDTwo\ this is clear;
for the remaining cases
we need the following description of $C'(\Gamma\sep d-1)$.

\proc{Lemma \label\TautCondition} Let
$a=\sum_{\sigma\in S(\Gamma),H\in\indexset^\star} 
a_{\sigma, H}\sigma H$ be an element of 
$C'(\Gamma\sep d)$.
Then $a \in C'(\Gamma\sep d-1)$ iff 
$\sum_{H\in\indexset^\star}a_{\sigma,H} = 0$ 
for all $\sigma\in S(\Gamma)$.
\endproc

\prf We know that $a\in C'(\Gamma\sep d-1)$ iff
$\sum_{H\in\indexset^\star}a_{\sigma,H} H\in \Omega(A_{d-1})$
for all $\sigma\in S(\Gamma)$. By Lemma \Orthog,
$\Omega(A_{d-1}) = \Omega(B_0)^\perp$. Now $B_0$
is generated by $[1]_0$ and 
$$\Omega([1]_0)\cdot\sum_{H\in\indexset^\star}a_{\sigma,H} H =
\sum_{H\in\indexset^\star}a_{\sigma,H}.$$
The result follows.
\qed

If a $G(d)$--colored graph $\Gamma$ is $(d-1)$--taut,
we shall say simply that $\Gamma$ is taut.
\proc{Lemma \label\MCircIsTaut} 
If $d=3$ and $\Gamma$ has an unsplittable $G$--coloring
with a special circuit, then $\Gamma$ is taut.
\endproc

\prf Since $\Gamma$ is simple, 
we may use the
natural triangulation in which the 0--simplices are the
vertices and the 1--simplices are the edges.
Let $H_0 \in \indexset^\star$ be such that $\Gamma_0=\Gamma_{H_0}$
is a special circuit, and 
let the non-trivial
elements of $H_0$ be $h_1$, $h_2$ and $h_3$.
The remaining elements of $\indexset^\star$ fall into three pairs depending on their intersections with $H_0$; we let $H_i$
and $H_i'$ be those for which
$H_i\cap H_0 = \genby{h_i}=H_i'\cap H_0$.
We also let $\Gamma_i = \Gamma_{H_i}$ and
$\Gamma'_i = \Gamma_{H'_i}$. 
Let $\sum_{\sigma,H}a(\sigma,H)\sigma H$ be any 1--cycle 
of $C'(\Gamma\sep 2)$, the sum
being over edges $\sigma$ and $H\in\indexset^\star$.
Since $\Gamma_0$ is connected, $a(\sigma,H_0)$ is constant 
on $\Gamma_0$. For notational simplicity, we show only that
$a(\sigma,H_1)$ is constant on $\Gamma_1$.

Let $S$
be the set of all edges colored $h_3$.
If $\sigma\in S$, we have 
$a(\sigma,H_0) = a(\sigma,H_3) = a(\sigma,H'_3) = 0$ and
$$a(\sigma,H_1) + a(\sigma,H'_1) + a(\sigma,H_2) + a(\sigma,H'_2) = 0.
\eqno{(\label\tempeq)}$$
Define an equivalence relation $\sim$ on $S$ by setting
$\sigma_1 \sim \sigma_2$ if 
$$a(\sigma_1,H_1) + a(\sigma_1,H'_1) 
= a(\sigma_2,H_1) + a(\sigma_2,H'_1).$$
Suppose that $\sigma_1$ and $\sigma_2$
are in $S$ and each have a vertex in common with an edge $\tau$
of $\Gamma\setminus\Gamma_0$. If the color of $\tau$ is $h_2$ then 
$\sigma_1$ and $\sigma_2$ lie in the same component of
$\Gamma_1$, and in the same component of $\Gamma'_1$.
Therefore $a(\sigma_1,H_1) = a(\sigma_2,H_1)$ and
$a(\sigma_1,H'_1) = a(\sigma_2,H'_1)$, so $\sigma_1\sim\sigma_2$.
Now, by (\tempeq), $\sigma_1\sim\sigma_2$ iff
$$a(\sigma_1,H_2) + a(\sigma_1,H'_2) 
= a(\sigma_2,H_2) + a(\sigma_2,H'_2),$$
and it follows similarly that $\sigma_1\sim\sigma_2$ if $\tau$
has color $h_1$.
Since $\Gamma\setminus\Gamma_0$ is connected, 
it follows that $\sigma_1\sim\sigma_2$
for any $\sigma_1$ and $\sigma_2$ in $S$.

Now define an equivalence relation $\approx$ on $S$
by setting $\sigma_1 \approx \sigma_2$ if
$a(\sigma_1,H_1) = a(\sigma_2,H_1)$. If $\sigma_1$
and $\sigma_2$ belong to the same component
of $\Gamma_1$ then $\sigma_1 \approx \sigma_2$.
Since $\sigma_1\sim\sigma_2$, we have $\sigma_1\approx\sigma_2$
iff $a(\sigma_1,H_1') = a(\sigma_2,H_1')$, and so
$\sigma_1 \approx \sigma_2$ if $\sigma_1$
and $\sigma_2$ belong to the same component
of $\Gamma_1'$.
Now $\Gamma_1\cup\Gamma'_1$ is the result of
deleting all edges colored $h_1$ from $\Gamma$,
which is connected since $\Gamma$ is unsplittable.
It follows that
$\sigma_1 \approx\sigma_2$ for all $\sigma_1$ and $\sigma_2$ in $S$;
i.e., that $a(\sigma,H_1)$ is constant on $S$.
Now any component of $\Gamma_1$ contains an edge
of $S$, 
so $a(\sigma, H_1)$ is constant
on $\Gamma_1$.
\qed

\proc{Lemma \label\MobIsTaut} Let $d=3$, and let $\Gamma$
be a M\"obius ladder with a $G$--coloring in which the
product of the colors on the rungs is non-trivial. Then
$\Gamma$ is taut.
\endproc

\prf We make $\indexset$ into an (additive) abelian group
by setting $H+K = \ker(\achar_H+\achar_K)$. For $H\in\indexset$,
we let $w_H=\sum_{\sigma\in S_1(\Gamma_H)}\sigma
=\sum_{\sigma\in S_1(\Gamma)}\achar_H(g_\sigma)\sigma
\in Z_1(\Gamma;\zt)$.
Then $w_{H+K}=w_H+w_K$, and the $w_H$ form a subgroup
of $Z_1(\Gamma;\zt)$ isomorphic to $\indexset$.
A 1--cycle of $C'(\Gamma\sep d)$ may be written in the form
$z = \sum_{H\in\indexset^\star}z_H H$, with 
$z_H\in Z_1(\Gamma_H;\zt)$. Then (by Lemma \TautCondition)
$z$ is in $C'(\Gamma\sep d-1)$ iff $\sum_{H\in\indexset^\star}z_H=0$
(the sum being taken in $Z_1(\Gamma;\zt)$), while $z$ is in
$W(\Gamma\sep d)$ iff each $z_H$ is a multiple of $w_H$.
Therefore $\Gamma$ is taut iff, given $z_H\in Z_1(\Gamma_H)$
for $H\in\indexset^\star$ with $\sum_{H\in\indexset^\star}z_H=0$,
each  $z_H$ is a multiple of $w_H$; since $\Gamma_G$ is empty,
we may replace $\indexset^\star$ by $\indexset$ in this statement.

Now let $g_0$ be the product of the colors on the rungs.
If the color of a rim edge $\sigma_i$ is $g$, then the color
of the edge $\sigma_{i+n}$ (where $n$ is the number of rungs)
is $gg_0$, so $g\neq g_0$. Since $G$ is generated by the colors
on the rungs and a single edge of the rim, at least two
distinct elements appear as rung colors. Thus the edges
colored $g_0$ form a proper subset of the rungs, and deleting them
leaves a M\"obius ladder $\Gamma'$. The $G$--coloring of $\Gamma$
induces a $G'$--coloring of $\Gamma'$, where 
$G'=G/\genby{g_0}\cong \ztb^2$. Since $\Gamma'$ is connected, it is taut.
Let $\indexset'$ be the set of $H\in\indexset$ such that $g_0\in H$.
There is a bijection $\indexset'\to\indexset(G')$ given by
$H\mapsto H'=H/\genby{g_0}$, and $\Gamma'_{H'}=\Gamma_H$.
By Lemma \EdgeParity, $\Gamma_H$ contains an even number of rungs
iff $H\in\indexset'$. Suppose now that $z_H\in Z_1(\Gamma_H;\zt)$
for $H\in\indexset$ and $\sum_{H\in\indexset}z_H=0$,
and set 
$z= \sum_{H\in\indexset'}z_H = \sum_{H\in\indexset-\indexset'}z_H$.
If $H\notin \indexset'$, $\Gamma_H$ is connected, so $z_H$
is automatically a multiple of $w_H$. It follows that
$z$ is equal to $w_K$ for some $K\in\indexset$. For
$H\in\indexset'$, each component of $\Gamma_H$ contains
zero or two rungs, and so the sum of the coefficients
of the rungs in $z_H$ is zero. Therefore the same is true
of $w_K$, which implies that $K\in\indexset'$.
Now $(z_K+w_K)+\sum_{H\in\indexset'-\{K\}}z_H=0$, and it follows
from the tautness of $\Gamma'$ that $z_H$ is a multiple of $w_H$
for $H\in\indexset'$ as well. Therefore $\Gamma$ is taut.
\qed

Much of the approach outlined in \S2 goes through for
any taut $G$--colored graph, but not for non-taut
graphs. This raises the question of how extensive
the class of taut graphs is. For M\"obius ladders,
one can determine all the taut colorings. If $d=3$ then
any taut coloring satisfies the hypothesis of Theorem
\ThmDThreeMob, apart from the exceptional coloring
of the 4--rung ladder in Example \ExMobDThree. (Actually, there
is a coloring of the 3--rung ladder for which the product
of the colors on the rungs is $1$, but an automorphism of the graph
takes it to one for which the product is non-trivial.)
For $d=4$, apart from the colorings of Example \ExMobDFour,
all taut colorings are obtained as follows. Suppose that $n\geq 4$,
and and let $\{x_1,x_2,x_3,x_4\}$ be a basis of $G$. Give the colors
$x_1$, $x_2$, $x_3$, and $x_1x_2x_3x_4^{n-1}$ to one rung each,
and give all other rungs the color $x_4$. It is possible to
complete the coloring, and the result is taut (but not 4--taut).
For $d\geq 5$, there is no taut coloring of any M\"obius ladder.

Also, the operation
of Figure \FigWyeDelta\ takes taut graphs to taut graphs,
and so generates infinitely many further examples
for $d\leq 5$; I know of no taut graphs for $d\geq 6$.
Even for taut graphs, one encounters some difficulties
which will be discussed after Lemma \ResultSeqTaut,
and which I have been unable to overcome for the added examples
just mentioned.

\section{The chain complex $\ker\beta/\im\alpha$}

We now return to the consideration of a regular branched covering
$\maps{\pi}{\cover M}{M}$ of a homology 3--sphere,
with deck group $G$ and branch set a $G$--colored graph $\Gamma$,
and of the chain maps $\alpha$, $\beta$ and $\gamma$ defined
in \S2.
For each simplex $\sigma$ of $M$, choose a lift $\tilde\sigma$
of $\sigma$ to $\cover M$, 
and for $H\in\indexset$ let 
$\sigma_H=\pi_H(\tilde\sigma)$. (In particular,
$\sigma_G=\sigma$.) Let $G_\sigma$ be the stabilizer of $\tilde\sigma$,
and let $A_k^\sigma = A_k(G, G_\sigma)$. If $\sigma$ is a 
simplex of $\Gamma$, these definitions agree with those of the
previous section; otherwise, $G_\sigma = 1$ and
$A_k^\sigma = 0$.
Also let $\indexset_\sigma$ be the set of 
$H\in\indexset$ such that $H\geq G_\sigma$,
and $\indexset^\star_\sigma = \indexset_\sigma -\{G\}$.

\proc{Lemma \label\GensModImA} 
$\bigl(\bigoplus_{H\in\indexset}C(M_H)\bigr)\big/\im \alpha$
is generated by the $\sigma_H H$ for $\sigma\in S(M)$ and
$H\in\indexset$, and $\im\beta$ is generated by the 
$\pi_H\shriek(\sigma_H)$.
\endproc

\prf If $H=G$ or $H\notin \indexset_\sigma$,
then $\sigma_H$ is the unique lift of $\sigma$ to $M_H$,
while if $H\in\indexset^\star_\sigma$
there is one other lift $\sigma_H'$ of
$\sigma$ to $M_H$. In the last case, 
$\alpha(\sigma H) = (\sigma_H+\sigma_H')H - \sigma_G G$.
This gives the first statement, and the second follows
since $\im\alpha\leq \ker\beta$.
\qed

For $\sigma\in S(M)$ and $g\in G$, the simplex
$g\tilde\sigma$ of $\cover M$ depends only on
the image of $g$ in $G/G_\sigma$.
We fix once and for all a right inverse for the projection
$G\to G/G_\sigma$, and thereby identify $G/G_\sigma$
with a complement of $G_\sigma$ in $G$.
A basis for $C(\cover M)$ is given by all
$g\tilde\sigma$ for $\sigma\in S(M)$ and 
$g\in G/G_\sigma$. 
Note that there is a bijection
$\indexset_\sigma\to\indexset(G/G_\sigma)$, namely 
$H\mapsto H/G_\sigma$.

\proc{Lemma \label\Indep} For each $\sigma\in S(M)$, the elements 
$\pi_H\shriek(\sigma_H) \in C(\cover M)$ for $H\in \indexset_\sigma$
are linearly independent. For $H\in\indexset-\indexset_\sigma$,
$2\pi_H\shriek(\sigma_H)=\pi_G\shriek(\sigma_G)$.
\endproc

\prf For $H\in\indexset_\sigma$, we have
$$\pi_H\shriek(\sigma_H)
=\sum_{h\in H}h\tilde\sigma
=|G_\sigma|\sum_{h\in H/G_\sigma}h\tilde\sigma
=|G_\sigma|\sum_{g\in G/G_\sigma}
\half(\mchar_H(g)+1)g\tilde\sigma.
$$
Let $T$ be the matrix with rows indexed by $H\in\indexset_\sigma$,
columns indexed by $g\in G/G_\sigma$, and entries
$\mchar_H(g)$, and let $J$ be the matrix with all entries $1$.
To prove the first statement, we must show that 
$\det(T+J)\neq 0$. Now $T$ is just
the character table of $G/G_\sigma$, and the 
orthogonality relations show that $\det T\neq 0$
(in fact, that $\det T= \pm n^{n/2}$
where $n=|G/G_\sigma|$). Expand $\det(T+J)$ by multilinearity
in the rows. Since the row of $T$ corresponding to
$G\in\indexset_\sigma$ consists entirely of ones,
all but two of the terms are zero, and the remaining
two are equal to $\det T$, so $\det(T+J)=2\det T\neq 0$.

Now let $H\in\indexset-\indexset_\sigma$. Then
$\rho_H\shriek(\sigma_G)=2\sigma_H$, so 
$\pi_G\shriek(\sigma_G)=\pi_H\shriek\rho_H\shriek(\sigma_G)
=2\pi_H\shriek(\sigma_H)$.
\qed

\proc{Lemma \label\KerBImA}
The chain complex $\ker\beta/\im \alpha$ is isomorphic to
$C'(\Gamma\sep d-1)$.
\endproc

\prf The complex $C'(\Gamma\sep d-1)$ was defined as 
a subcomplex of $C(\Gamma;\zt)^{\indexset^\star}$,
which is in turn a subcomplex of $C(M;\zt)^{\indexset^\star}$.
As a subcomplex of $C(M;\zt)^{\indexset^\star}$,
$C'(\Gamma\sep d-1)$ consists of those chains
$\sum_{\sigma\in S(M),H\in\indexset^\star}a_{\sigma,H}\sigma H$
such that, for each $\sigma$, 
$\sum_{H\in\indexset^\star}a_{\sigma,H} = 0$ and 
$a_{\sigma,H}=0$ if $H\in\indexset_\sigma^\star$
(because these equations imply that $a_{\sigma,H}=0$
whenever $\sigma$ is not in $\Gamma$).

Let $\maps{\bar\rho_H}{C(M_H)}{C(M;\zt)}$ be the composite
of $\maps{\rho_H}{C(M_H)}{C(M)}$ and reduction of the coefficients
modulo 2. Define 
$\maps{\zeta}{\bigoplus_{H\in\indexset}C(M_H)}
{C(M;\zt)^{\indexset^\star}}$
by
$$\eqalignno{
\zeta\Bigl(\sum_{H\in\indexset}d_H H\Bigr) 
&= \sum_{H\in\indexset^\star} \bar\rho_H(d_H)H
\quad\hbox{for $d_H\in C(M_H)$, $H\in\indexset$.}\cr
\noalign{\hbox{Then}}
\zeta\alpha\Bigl(\sum_{H\in\indexset^\star}c_H H\Bigr)
&= \sum_{H\in\indexset^\star}\bar\rho_H\rho_H\shriek(c_H)H
=0\quad\hbox{for $c_H\in C(M)$, $H\in\indexset^\star$.}
}$$
Thus $\zeta$ induces a map from
$\bigl(\bigoplus_{H\in\indexset}C(M_H)\bigr)\big/\im \alpha$
to $C(M;\zt)^{\indexset^\star}$; we shall show that 
$\ker\beta/\im\alpha$
is mapped isomorphically to $C'(\Gamma\sep d-1)$. By Lemmas
\GensModImA\ and \Indep, any element
of $\ker\beta/\im\alpha$ has a representative of the form
$$c= \sum_{\sigma\in S(M)}\Bigl(a_\sigma \sigma G 
+\sum_{H\notin\indexset_\sigma}b_{\sigma, H} \sigma_H H\Bigr)
\quad{\rm for\ }a_\sigma, b_{\sigma,H}\in\z,
$$
and such an element is in $\ker\beta$ iff
$2a_\sigma+\sum_{H\notin\indexset_\sigma}b_{\sigma,H}=0$
for each $\sigma$. It follows immediately that the image
of $\ker\beta/\im\alpha$ in $C(M;\zt)^{\indexset^\star}$
is $C'(\Gamma\sep d-1)$. Further, the chain
c is in $\ker\zeta$ iff each $b_{\sigma,H}$
is even, and then
$$\eqalign{
\alpha\Bigl(\sum_{\sigma\in S(M),H\notin\indexset_\sigma}
\half b_{\sigma,H}\sigma H\Bigr)
&= \sum_{\sigma\in S(M),H\notin\indexset_\sigma}
\Bigl(b_{\sigma,H}\sigma_H H 
- \half b_{\sigma,H}\sigma G\Bigr)\cr
&=\sum_{\sigma\in S(M)}\Bigl(a_\sigma \sigma G 
+\sum_{H\notin\indexset_\sigma}b_{\sigma, H} \sigma_H H\Bigr)=c
}$$
provided $c \in\ker\beta$.
This completes the proof.
\qed

\proc{Lemma \label\HKerB} There is a short exact sequence
$$0\to \z^{\indexset^\star} \to H_0(\ker\beta)
\to H_0(\Gamma\sep d-1) \to 0,$$
and $H_1(\ker\beta) \cong H_1(\Gamma\sep d-1)$.
\endproc

\prf By Lemma \KerBImA, the sequence (\ShortSeqOne) becomes
$$0\to C(M)^{\indexset^\star} \mapright\alpha \ker\beta
\to C'(\Gamma\sep d-1) \to 0.$$
In the long exact homology sequence,
the map $H_1(\Gamma\sep d-1) \to H_0(M)^{\indexset^\star}$
is zero since $H_1(\Gamma\sep d-1)$ is torsion
and $H_0(M)\cong \z$. Therefore the long exact sequence
gives exact sequences
$$\eqalign{
&0\to H_1(\ker\beta) 
\to H_1(\Gamma\sep d-1) \to 0\cr
{\rm and}\quad &0\to \z^{\indexset^\star} \to H_0(\ker\beta)
\to H_0(\Gamma\sep d-1) \to 0.}$$
\qed

We now turn to the sequence (\ShortSeqTwo). Note that
the map induced on first homology by the map
$\maps\beta{\bigoplus_{H\in\indexset}C(M_H)}{\im\beta}$
from that sequence may be regarded as a map from
$\bigoplus_{H\in\indexset^\star}H_1(M_H)$ to $H_1(\im\beta)$
since $H_1(M_G)=H_1(M)=0$.

\proc{Lemma \label\ResultSeqOneLeft} If $\Gamma$ is taut, the map
$\beta$ from $\bigoplus_{H\in \indexset^\star}H_1(M_H)$
to $H_1(\im\beta)$ is injective.
\endproc

\prf By Lemma \HKerB,
part of the long exact sequence of (\ShortSeqTwo) becomes
$$H_1(\Gamma\sep d-1) \mapright\iota 
\bigoplus_{H\in \indexset^\star}H_1(M_H) 
\mapright\beta H_1(\im\beta).$$
We must show that the map $\iota$ in this sequence is trivial.
Any element of $H_1(\Gamma\sep d-1)=Z'_1(\Gamma\sep d-1)$ 
has the form $z=\sum_{H\in\indexset^\star}z_H H$,
where $z_H\in Z_1(\Gamma_H;\zt)$ and $\sum_{H\in\indexset^\star}z_H=0$
in $Z_1(\Gamma;\zt)$.
Let $\maps{\bar\rho_H}{C(M_H)}{C(M;\zt)}$ be 
as in the proof of Lemma \KerBImA. The inverse image of $\Gamma_H$
in $M_H$ is a link $L_H$, and we may take 
$w_H\in Z_1(L_H)\leq Z_1(M_H)$ with $\bar\rho_H(w_H)=z_H$.
Then $\sum_{H\in\indexset^\star}\bar\rho_H(w_H)=0$, so there is an
element $w_G$ of $Z_1(M)$ with
$2w_G+\sum_{H\in\indexset^\star}\rho_H(w_H) = 0$.
Let $w = \sum_{H\in\indexset}w_H H\in
\bigoplus_{H\in\indexset}Z_1(M_H)$. For $H\in\indexset^\star$,
$\rho_H\shriek\rho_H(w_H)=2w_H$ since $w_H$ is in $Z_1(L_H)$.
Therefore
$$0 = \pi\shriek\Bigl(2w_G+\sum_{H\in\indexset^\star}\rho_H(w_H)\Bigr)
= 2\pi\shriek(w_G)+\sum_{H\in\indexset^\star}\pi_H\shriek
\rho_H\shriek\rho_H(w_H)
=2\sum_{H\in\indexset}\pi_H\shriek(w_H)
=2\beta(w),$$
so $w\in\ker\beta$.
It follows from the proof of Lemma \KerBImA\ that
the element of $H_1(\ker\beta)$ represented by $w$
corresponds to $z$ under the isomorphism
$H_1(\ker\beta)\cong H_1(\Gamma\sep d-1)$ of Lemma \HKerB,
and so $\iota(z)$ is the element of 
$\bigoplus_{H\in \indexset^\star}H_1(M_H)$ represented by $w$.

Since $\Gamma$ is taut, for each $H\in\indexset^\star$, $z_H$
is a multiple of the mod 2 fundamental class of $\Gamma_H$,
and we may take $w_H$ to be a multiple of the fundamental
class of $L_H$ for some orientation of $L_H$.
Since $L_H$ bounds a lift of a Seifert surface
for $\Gamma_H$, $w_H$ represents zero in $H_1(M_H)$,
and so $\iota(z)=0$ as required.
\qed

\proc{Lemma \label\HZeroImB} We have $H_0(\im\beta)\cong \z$.
\endproc

\prf The end of the long exact sequence of (\ShortSeqTwo) shows that
$$\maps{\beta}{\bigoplus_{H\in\indexset}H_0(M_H)}{H_0(\im\beta)}$$
is surjective. In fact the restriction of $\beta$ to
$\bigoplus_{H\in\indexset^\star}H_0(M_H)$ is surjective
since $\pi_G\shriek$ factors through $\pi_H\shriek$ for any 
$H\in\indexset^\star$.
Let $H\in\indexset^\star$. For $\sigma \in S_0(M)$,
the 0--simplices $\sigma_H$ all represent the same generator
of $H_0(M_H)\cong\z$. The image $x_H$ of this generator
in $H_0(\im\beta)$ is represented by $\pi_H\shriek(\sigma_H)$
for any $\sigma\in S_0(M)$. Define an equivalence relation
on $\indexset^\star$ by setting $H\sim K$ if $x_H=x_K$. Suppose
that there is some $\sigma\in S_0(M)$ such that neither
$H$ nor $K$ is in $\indexset_\sigma$. Then, by Lemma \Indep,
$\pi_H\shriek(\sigma_H) = {1\over 2}\pi_G\shriek(\sigma_G)
=\pi_K\shriek(\sigma_K)$,
and so $H\sim K$. Now suppose $H_1$ and $H_2$ are any two elements
of $\indexset^\star$. For $i=1$ or $2$, there is some color
$g_i\notin H_i$ (since the colors generate $G$), and a 0--simplex
$\sigma_i$ of $\Gamma$ with $g_i \in G_{\sigma_i}$. Thus
$H_i \notin\indexset_{\sigma_i}$. We may find $K\in\indexset^\star$
containing neither $g_1$ nor $g_2$. 
Then $K\notin \indexset_{\sigma_i}$,
so $H_i\sim K$ for $i = 1$ or $2$. Therefore $H_1\sim H_2$,
and there is only one equivalence class. This shows that 
$H_0(\im\beta)$ is cyclic.

On the other hand, the image of $x_H$ under the 
map $H_0(\im\beta)\to H_0(\cover M)\cong\z$ induced by inclusion 
is $2^{d-1}$ times a generator, and therefore $H_0(\im\beta)\cong\z$.
\qed

\proc{Lemma \label\ResultSeqOne} If $\Gamma$ is taut,
there is a short exact sequence
$$0 \to \bigoplus_{H\in \indexset^\star}H_1(M_H)
\mapright\beta H_1(\im\beta)
\to \ztb^{b_1(\Gamma)-d} \to 0.$$
\endproc

\prf By the previous lemma,
part of the long exact sequence of (\ShortSeqTwo) is
$$\bigoplus_{H\in \indexset^\star}H_1(M_H)
\mapright\beta H_1(\im\beta)
\to H_0(\ker\beta)
\mapright\iota \z^\indexset \to \z \to 0,\eqno{(\label\tempeq)}$$
and the first map is injective by Lemma \ResultSeqOneLeft.
It remains to prove that 
$\ker\iota\cong \ztb^{b_1(\Gamma)-d}$.
We show that there is a commutative diagram
{\abovedisplayskip=\diagdisplayskip  
\belowdisplayskip=\diagdisplayskip   
$$
\def\normalbaselines{\diagbaselines} 
\matrix{
0 & \to & \z^{\indexset^\star} & \mapright\phi & H_0(\ker\beta) 
& \mapright\psi & \ztb^{b_0(\Gamma\sep d-1)} &\to& 0\cr
&& \mapdown\id && \mapdown{\iota} && \mapdown\theta &&  \cr
0 & \to & \z^{\indexset^\star} & \mapright{\iota\phi} & \z^\indexset 
& \mapright{} & \z\oplus \ztb^{|\indexset|-2} & \to & 0\cr
&&&& \mapdown{} &&\mapdown{} &&\cr
&&&& \z & \mapright\id & \z &&\cr
&&&& \downarrow && \downarrow &&\cr
&&&& 0 && 0 &&
}$$}%
in which the rows and the central column are exact.
The central column is part of (\tempeq), and the top row is
the short exact sequence of Lemma \HKerB.
The composite $\maps{\iota\phi}{\z^{\indexset^\star}}{\z^\indexset}$
is the map 
$H_0(M)^{\indexset^\star}\to\bigoplus_{H\in\indexset}H_0(M_H)$
induced by $\alpha$, so it is given by
$\iota\phi\bigl(\sum_{H\in\indexset^\star}a_H H\bigr) = 
\sum_{H\in\indexset^\star}2a_H H 
- \bigl(\sum_{H\in\indexset^\star}a_H\bigr) G$.
This is injective and 
has cokernel isomorphic to $\z\oplus \ztb^{|\indexset|-2}$,
so we obtain the exact second row.
The maps in the right-hand column may now be 
defined to make the diagram commute.
Diagram-chasing shows that the right-hand column is exact, 
so that $\im\theta \cong \ztb^{|\indexset|-2}$
and $\ker\theta \cong \ztb^{b_0(\Gamma\sep d-1)-|\indexset|+2}$.
More diagram-chasing shows  that $\psi$ maps
$\ker \iota$ isomorphically onto $\ker \theta$.
Since $\Gamma$ is taut, $b_1(\Gamma\sep d-1) = \dim B_{d-2}
= 2^d-d-1$ and $\Gamma$ is connected. By Lemma \ChiK,
$\chi(\Gamma\sep d-1)=-\chi(\Gamma)=b_1(\Gamma)-1$.
Hence
$$b_0(\Gamma\sep d-1)-|\indexset|+2
=2^d-d-1+b_1(\Gamma)-1 -2^d+2 = b_1(\Gamma)-d,$$
and we are done.
\qed

\section{The chain complex $\ker\gamma/\im\beta$}

We may regard $C(\cover M)$ as a $\z[G]$--module.
As such, it is generated (though not freely) by the 
$\tilde \sigma$ for $\sigma\in S(M)$.
For $1\leq k\leq d-1$, we let $D(k)$ be the subcomplex 
of $C(\cover M)$ consisting of chains
$c=\sum_{\sigma\in S(M)} \lambda_\sigma\tilde\sigma$
($\lambda_\sigma\in\z[G]$) satisfying, for all $\sigma$,
$$\eqalign {
\mchar_G(\lambda_\sigma) &\equiv 0 \pmod {2^{d-1}}\cr
{\rm and}\quad \mchar_H(\lambda_\sigma) &\equiv 0 \pmod {2^k}
\quad{\rm for\ }H\in\indexset^\star_\sigma.
}$$
This is well-defined because, for $\lambda\in\z[G]$
and $\sigma\in S(M)$, the chain $\lambda\cover\sigma$
determines the image $\bar\lambda$ of $\lambda$ in $\z[G/G_\sigma]$,
and hence determines 
$\mchar_H(\lambda)=\mchar_{H/G_\sigma}(\bar\lambda)$
for $H\in \indexset_\sigma$.
By Lemma \JKBasis, $I^k C(\cover M)\leq D(k)$. Recall that
we have identified $G/G_\sigma$ with a subgroup of $G$,
and hence $\z[G/G_\sigma]$ with a subring of $\z[G]$.

\proc{Lemma \label\ImBEqDK} We have $\ker\gamma= D(1)$ and
$\im\beta = D(d-1)$.
\endproc

\prf From the definition of $\gamma$, 
$\sum_{\sigma\in S(M)} \lambda_\sigma\tilde\sigma\in\ker\gamma$
iff $\mchar_G(\lambda_\sigma)\equiv 0\pmod{2^{d-1}}$ for
each $\sigma$. Since $\mchar_H(\lambda)\equiv \mchar_G(\lambda)
\pmod 2$ for all $H\in\indexset$ and $\lambda\in\z[G]$,
it follows that $\ker\gamma=D(1)$. For $\sigma\in S(M)$,
let $(\im\beta)_\sigma = \im\beta\cap \z[G]\tilde\sigma$.
To show that $\im\beta=D(d-1)$, it is enough to show
that $\lambda\tilde\sigma \in (\im\beta)_\sigma$ iff
$\mchar_H(\lambda)\equiv 0\pmod{2^{d-1}}$ for all
$H\in \indexset_\sigma$. We may assume that 
$\lambda=\sum_{g\in G/G_\sigma}\lambda_g g\in\z[G/G_\sigma]$,
and so $\mchar_H(\lambda)=\mchar_{H/G_\sigma}(\lambda)$ for 
$H\in\indexset_\sigma$.
Consider the chain 
$\sum_{H\in\indexset_\sigma}\mchar_H(\lambda)\pi_H\shriek(\sigma_H)
\in C(\cover M)$. We have
$$\eqalignno{
\sum_{H\in\indexset_\sigma}\mchar_H(\lambda)\pi_H\shriek(\sigma_H)
&=\sum_{H\in\indexset_\sigma}\mchar_{H/G_\sigma}(\lambda)
|G_\sigma|\sum_{h\in H/G_\sigma}h\tilde\sigma\cr
&=|G_\sigma|\sum_{g\in G/G_\sigma, H\in\indexset_\sigma}
\mchar_{H/G_\sigma}(\lambda)\half(\mchar_{H/G_\sigma}(g)+1)
g\tilde\sigma\cr
&=\half|G_\sigma|\sum_{g\in G/G_\sigma}\Bigl(
\sum_{H\in\indexset_\sigma}
\bigl(\mchar_{H/G_\sigma}(\lambda g)
+\mchar_{H/G_\sigma}(\lambda)\bigr)\Bigr)
g\tilde\sigma\cr
&=\half|G_\sigma|\sum_{g\in G/G_\sigma}|G/G_\sigma|
(\lambda_g+\lambda_1)g\tilde\sigma.\cr
\noalign{\hbox{That is,}}
\sum_{H\in\indexset_\sigma}\mchar_H(\lambda)\pi_H\shriek(\sigma_H)
&= 2^{d-1}\lambda\tilde\sigma + 
\half|G/G_\sigma|\lambda_1\pi_G\shriek(\sigma_G).&(\label\temp)\cr
}$$
Suppose first that $G_\sigma = 1$. Then $\indexset_\sigma=\indexset$,
and by Lemmas \GensModImA\ and \Indep, a basis for
$(\im\beta)_\sigma$ consists of the $\pi_H\shriek(\sigma_H)$
for $H\in\indexset$. In this case, (\temp) gives
$$\sum_{H\in\indexset}\mchar_H(\lambda)\pi_H\shriek(\sigma_H)
= 2^{d-1}\bigl(\lambda\tilde\sigma + 
\lambda_1\pi_G\shriek(\sigma_G)\bigr),
$$
and it follows that $\lambda\tilde\sigma\in(\im\beta)_\sigma$
iff $\mchar_H(\lambda)\equiv 0\pmod{2^{d-1}}$ for all
$H\in \indexset$. Now suppose that $G_\sigma\neq 1$.
Then a basis for
$(\im\beta)_\sigma$ consists of the $\pi_H\shriek(\sigma_H)$
for $H\in\indexset^\star_\sigma$ and $\pi_{H_0}\shriek(\sigma_{H_0})$
for any one $H_0\in\indexset-\indexset_\sigma$, and (\temp) gives
$$ 2\mchar_G(\lambda)\pi_{H_0}\shriek(\sigma_{H_0})+
\sum_{H\in\indexset_\sigma^\star}
\mchar_H(\lambda)\pi_H\shriek(\sigma_H)
= 2^{d-1}\lambda\tilde\sigma + 
|G/G_\sigma|\lambda_1\pi_{H_0}\shriek(\sigma_{H_0}).
$$
In this case, $\lambda\tilde\sigma\in(\im\beta)_\sigma$
iff $\mchar_H(\lambda)\equiv 0\pmod{2^{d-1}}$ for all
$H\in \indexset_\sigma^\star$ and 
$2\mchar_G(\lambda)\equiv |G/G_\sigma|\lambda_1\pmod{2^{d-1}}$.
But $|G/G_\sigma|\lambda_1 = 
\sum_{H\in\indexset_\sigma}\mchar_H(\lambda)$, so this is true iff
$\mchar_H(\lambda)\equiv 0\pmod{2^{d-1}}$ for all
$H\in \indexset_\sigma$,
as required.
\qed

Thus we have a filtration
$\im\beta = D(d-1)\leq\cdots\leq D(1) = \ker\gamma$,
and instead of dealing directly with the complex
$\ker\gamma/\im\beta$, we consider the quotients of this
filtration.

The following notation will be used in the proofs of the next lemma
and Lemma \labelplus2\ResultSeqTwo.
Let $\sigma\in S(M)$, and let
$\partial\sigma = \sum_{\tau\in S(M)}i_{\sigma,\tau}\tau$.
Thus $i_{\sigma,\tau}=\pm 1$ if $\tau$ is a face of $\sigma$,
and $i_{\sigma,\tau}=0$ otherwise. If $\tau$ is a face of $\sigma$,
there is a unique element $g_{\sigma,\tau}$ of $G/G_\tau\leq G$
such that $g_{\sigma,\tau}\tilde\tau$ is a face
of $\tilde\sigma$; we set $g_{\sigma,\tau}=1$ otherwise.
Then $\partial\tilde\sigma = 
\sum_{\tau\in S(M)}i_{\sigma,\tau}g_{\sigma,\tau}\tilde\tau$.

Recall that $C(\Gamma\sep k)$ was defined as a subcomplex
of $C(\Gamma;A_k)$, which is in turn a 
subcomplex of $C(M;A_k)$.
$C(\Gamma\sep k)$ is the subcomplex of $C(M;A_k)$ generated
by all chains $a\sigma$ where $\sigma$ is a simplex of
$M$ and $a \in A_k^\sigma$, because $A_k^\sigma = 0$
if $\sigma$ is not a simplex of $\Gamma$.

\proc{Lemma \label\DKQuot} For $1\leq k \leq d-2$, we
have $D(k)/D(k+1) \cong C(M;A_k)/C(\Gamma\sep k)$.
\endproc

\prf Since $C(M;\z[G])$ is the free $\z[G]$--module on the simplices
of $M$, there is a unique $\z[G]$--module homomorphism
$\eta$ from $C(M;\z[G])$ to $C(\cover M)$ sending 
$\sigma \in S(M)$ to $\tilde\sigma$; of course, 
$\eta$ is not a chain map. Nevertheless, its kernel is
a subcomplex; it is generated by $\lambda\sigma$ for
$\sigma\in S(M)$ and $\lambda\in \z[G,G_\sigma]$.
The subcomplex $C(M;I^k)$ is sent by $\eta$ to
$I^kC(\cover M)\leq D(k)$; 
the kernel of $\eta\restrict C(M;I^k)$
is the subcomplex $E(k)$
generated by $\lambda\sigma$ for
$\sigma\in S(M)$ and $\lambda\in I^k[G,G_\sigma]$.
For $1\leq k\leq d-2$, we may identify $C(M;I^k)/C(M;I^{k+1})$
with $C(M;A_k)$, and $E(k)/E(k+1)$ with $C(\Gamma\sep k)$.
Then we have an induced
map $\bar\eta_k$ from
$C(M;A_k)$ to 
$D(k)/D(k+1)$, whose kernel
contains $C(\Gamma\sep k)$.
For $\lambda\in I^k$ and $\sigma\in S(M)$, we have
$$
(\eta\partial-\partial\eta)(\lambda\sigma)
=\sum_{\tau\in S(M)}i_{\sigma,\tau}
\lambda(1-g_{\sigma,\tau})\tilde\tau
\in I^{k+1}C(\cover M)\leq D(k+1),
$$
which shows that $\bar\eta_k$ is a chain map.

Suppose that $\lambda\tilde\sigma\in D(k)$. We may assume that
$\lambda \in\z[G/G_\sigma]$. For any
$H\in \indexset$, there is some $H'\in \indexset_\sigma$
with $H\cap G/G_\sigma = H'\cap G/G_\sigma$, and so
$\mchar_H(\lambda)=\mchar_{H'}(\lambda)\equiv 0 \pmod{2^k}$.
By Lemma \JKBasis, $\lambda-\mchar_G(\lambda)\in I^k$;
also $\eta\bigl((\lambda-\mchar_G(\lambda))\sigma\bigr)
=\lambda\tilde\sigma-\mchar_G(\lambda)\tilde\sigma$.
But $\mchar_G(\lambda)\equiv 0\pmod{2^{d-1}}$,
and so $\mchar_G(\lambda)\tilde\sigma\in D(k+1)$. Therefore
$\bar\eta_k$ maps $C(M;A_k)$ onto
$D(k)/D(k+1)$.

Next, suppose $\lambda\in I^k$
and $\sigma\in S(M)$ are such that $[\lambda]_k\sigma$
is in the kernel of $\bar\eta_k$; that is,
$\lambda\tilde\sigma\in D(k+1)$. Take $\mu\in \z[G/G_\sigma]$
so that $\mu\tilde\sigma=\lambda\tilde\sigma$.
As before, for 
$H\in \indexset$, there is some $H'\in \indexset_\sigma$
with $H\cap G/G_\sigma = H'\cap G/G_\sigma$, and so
$\mchar_H(\mu)=\mchar_{H'}(\mu)=\mchar_{H'}(\lambda)
\equiv 0 \pmod{2^{k+1}}$.
Also $\mchar_G(\mu)=\mchar_G(\lambda)=0$, so it follows from
Lemma \JKBasis\ that $\mu\in I^{k+1}$. Since 
$\lambda\tilde\sigma = \mu\tilde\sigma$, 
$\lambda-\mu \in I^k[G,G_\sigma]$, so
$[\lambda]_k= [\lambda-\mu]_k$ is in
$A_k^\sigma$.
It follows that the kernel of $\bar\eta_k$ is equal
to $C(\Gamma\sep k)$, and so $\bar\eta_k$ induces
the desired isomorphism of chain complexes from 
$C(M;A_k)/C(\Gamma\sep k)$ to $D(k)/D(k+1)$.
\qed

\proc{Lemma \label\HDModD} For $1\leq k \leq d-2$,
we have $H_0(D(k)/D(k+1))=0$,
$H_1(D(k)/D(k+1))\cong \ztb^{b_0(\Gamma\sep k)-\dim A_k}$,
and $H_2(D(k)/D(k+1))\cong H_1(\Gamma\sep k)$.
\endproc

\prf Lemma \DKQuot\ gives a short exact sequence
$$0\to C(\Gamma\sep k) \to C(M;A_k) \to D(k)/D(k+1)\to 0.$$
The long homology sequence 
gives exact sequences
$$\eqalign{
&0 \to H_2(D(k)/D(k+1))\to H_1(\Gamma\sep k) \to 0\cr
{\rm and}\quad&0\to H_1(D(k)/D(k+1))
\to H_0(\Gamma\sep k) \to A_k \to H_0(D(k)/D(k+1))\to 0.
}$$
The map $H_0(\Gamma\sep k) \to A_k$ in the second of
these has image containing $\omega(g)b$ for any
$b\in B_{k-1}$ and any $g\in G$ that appears as an
edge color. Since the colors generate $G$, this map
is onto, and the result follows.
\qed

\proc{Lemma \checklabel\ResultSeqTwo} If $1\leq k \leq d-2$
and $\Gamma$ is $k$--taut, there is a short exact sequence 
$$0\to H_1(D(k+1)) \to H_1(D(k))
\to \ztb^{{d-2 \choose k-1}(b_1(\Gamma)-1)-{d\choose k}+1}\to 0.$$
\endproc
\prf By Lemma \HDModD, part of the long exact sequence of 
$$0\to D(k+1)\to D(k) \to D(k)/D(k+1) \to 0$$ is
$$\eqalign{
H_2(D(k)) \mapright{\phi_k} H_1(\Gamma\sep k) \to H_1(D(k+1)) 
&\to H_1(D(k)) \mapright{\psi_k} 
\ztb^{b_0(\Gamma\sep k)-\dim A_k}\cr
&\to H_0(D(k+1)) \to H_0(D(k)) \to 0,
}$$
whether or not $\Gamma$ is $k$--taut.
Suppose that $H_0(D(k+1))\cong \z$. It follows that
$\psi_k$ is onto, and that $H_0(D(k))\cong\z$.
Since $H_0(D(d-1))\cong\z$ by Lemmas \ImBEqDK\ and \HZeroImB,
a downward induction on $k$ shows that $\psi_k$
is onto for $1\leq k\leq d-2$.
Now, since $\Gamma$ is $k$--taut, $b_1(\Gamma\sep k)= \dim B_{k-1}$
and $$\chi(\Gamma\sep k)= -\smallchoose{d-2}{k-1}\chi(\Gamma)
= \smallchoose{d-2}{k-1}(b_1(\Gamma)-1),$$ and so
$$\eqalign{
b_0(\Gamma\sep k)-\dim A_k 
&= \smallchoose{d-2}{k-1}(b_1(\Gamma)-1) +\dim B_{k-1}-\dim A_k \cr
&= \smallchoose{d-2}{k-1}(b_1(\Gamma)-1) - \smallchoose{d}{k} + 1.
}$$
It only remains to prove that
$\phi_k$ is onto.

In the rest of the proof, $\sigma$ always denotes a 3--simplex of $M$,
$\tau$ a 2--simplex, and $\upsilon$ a 1--simplex, so that, for
example, $\sum_{\sigma,\tau}$ indicates a sum over
$\sigma\in S_3(M)$ and $\tau\in S_2(M)$. 
We assume that the orientations of the 3--simplices of $M$ are
induced by an orientation of $M$, so that $c=\sum_\sigma\sigma$
represents a generator of $H_3(M)$. Consider the chain 
$\tilde c=\sum_\sigma\tilde\sigma\in C_3(\cover M)$.
We have 
$\partial\tilde c \in \ker\bigl(C_2(\cover M) \to C_2(M)\bigr)
=IC_2(\cover M)$.
Let $\lambda\in J^{k-1}$. Then $\lambda I \leq I^k$,
so $\lambda\partial\tilde c \in I^kC_2(M)\leq D_2(k)$,
and the cycle $\lambda\partial\tilde c$ represents an element $x$
of $H_2(D(k))$. Now 
$$\lambda\partial\tilde c 
=\lambda\sum_{\sigma,\tau}i_{\sigma,\tau}g_{\sigma,\tau}\tilde\tau
=\eta\Bigl(\lambda\sum_{\sigma,\tau}
i_{\sigma,\tau}g_{\sigma,\tau}\tau\Bigr),
$$
where $\eta$ is the map $C(M;\z[G]) \to C(\cover M)$
from the proof of Lemma \DKQuot. It follows that $\phi_k(x)$ is
the image in $C_1(M;A_k)$ of
$c'=\partial\bigl(\lambda\sum_{\sigma,\tau}
i_{\sigma,\tau}g_{\sigma,\tau}\tau\bigr)$. Now
$$c' = \lambda\sum_{\sigma,\tau,\upsilon}
i_{\sigma,\tau}i_{\tau,\upsilon}g_{\sigma,\tau}\upsilon
= \lambda\sum_{\upsilon}\mu_\upsilon\upsilon
\quad\hbox{where}\quad\mu_\upsilon =\sum_{\sigma,\tau}
i_{\sigma,\tau}i_{\tau,\upsilon}g_{\sigma,\tau}\in \z[G].
$$
Fix a 1--simplex $\upsilon$ of $M$. Let the 2--simplices of $M$
having $\upsilon$ as a face be $\tau_1,\ldots,\tau_n$, and
the 3--simplices $\sigma_1,\ldots,\sigma_n$. Let $\sigma_0 =\sigma_n$,
and choose the numbering so that $\tau_j$ is a face of
$\sigma_{j-1}$ and $\sigma_j$ for $1\leq j\leq n$. 
Let $i_j=i_{\sigma_{j-1},\tau_j}i_{\tau_j,\upsilon}=\pm1$.
Then $i_{\sigma_j,\tau_j}i_{\tau_j,\upsilon} = -i_j$, so
$\mu_\upsilon =
\sum_{j=1}^n i_j(g_{\sigma_{j-1},\tau_j}-g_{\sigma_j,\tau_j})\in I$.
By Lemma \AOneEqG, 
$[\mu_\upsilon]_1 = \omega\bigl(\prod_{j=1}^n
g_{\sigma_{j-1},\tau_j} g_{\sigma_j,\tau_j}\bigr)$.
Considering a lift to $\cover M$ of a meridian of $\upsilon$,
we see that $\prod_{j=1}^n
g_{\sigma_{j-1},\tau_j} g_{\sigma_j,\tau_j}$
is the color $g_\upsilon$ if $\upsilon$ is a 1--simplex
of $\Gamma$, and $1$ otherwise. Therefore
$$\phi_k(x) = \sum_{\upsilon\in S_1(\Gamma)}[\lambda]_{k-1}
\omega(g_\upsilon)\upsilon = \iota_k([\lambda]_{k-1}).
$$
Since $\Gamma$ is $k$--taut, this shows that $\phi_k$ is onto.
\qed

\proc{Lemma \checklabel\ResultSeqTaut} If $\Gamma$ is taut, 
there is a short exact sequence
$$0 \to \bigoplus_{H\in\indexset^\star}H_1(M_H)
\mapright\beta H_1(\cover M) \to \Lambda\to 0,$$
where $\Lambda$ satisfies $2^{d-1}\Lambda = 0$ and
$|\Lambda| = 2^m$ for $m=2^{d-2}(b_1(\Gamma)-5) + d+1$.
\endproc

\prf Since $D(d-1) = \im\beta$, $D(1)=\ker\gamma$
and $H_1(\ker\gamma)\cong H_1(\cover M)$, Lemmas
\ResultSeqOne\ and \ResultSeqTwo\ give an exact sequence
as claimed with $2^{d-1}\Lambda = 0$ and
$|\Lambda| = 2^m$ where $m$ is the sum of $b_1(\Gamma)-d$
and ${d-2 \choose k-1}(b_1(\Gamma)-1)-{d\choose k}+1$ for 
$1\leq k\leq d-2$. Since
${d-2 \choose k-1}(b_1(\Gamma)-1)-{d\choose k}+1$ is equal
to $b_1(\Gamma)-d$ when $k=d-1$,
$$m = \sum_{k=1}^{d-1}\left(
\smallchoose{d-2}{k-1}(b_1(\Gamma)-1)-\smallchoose{d}{k}+1\right)
=2^{d-2}(b_1(\Gamma)-1)-2^d + d+1=2^{d-2}(b_1(\Gamma)-5) + d+1.
$$
\qed

Even accepting the limitation to taut graphs,
Lemma \ResultSeqTaut\ is unsatisfactory in two respects.
First, it gives incomplete information about the group $\Lambda$.
(Theorem \ThmDThreeMCirc\ and 
Proposition \MCircExists\ show
that, at least for $d=3$, $\Lambda$ may be any group
satisfying the conditions of the lemma.)
Second, it gives no information
at all about the extension of
$\bigoplus_{H\in\indexset^\star}H_1(M_H)$ by $\Lambda$.
All the examples I know are consistent with the conjecture
that $\beta\bigl(\bigoplus_{H\in\indexset^\star}H_1(M_H)\bigr)
=2^{d-1}H_1(\cover M)$ whenever $\Gamma$ is taut, but I have been
unable to prove this. The following lemma suffices
in some cases.

\proc{Lemma \label\Extension} Let $\Gamma$ be taut and suppose
that, for every $H\in\indexset^\star$ such that $\Gamma_H$ is
disconnected, the cover $\maps{\pi_H}{\cover M}{M_H}$
can be factored through 2--fold covers
$\cover M=M_d \to\cdots\to M_2\to M_1=M_H$ so that
each transfer map $H_1(M_i;\zt)\to H_1(M_{i+1};\zt)$
is trivial. Then 
$\beta\bigl(\bigoplus_{H\in\indexset^\star}H_1(M_H)\bigr)
=2^{d-1}H_1(\cover M)$.
\endproc

\prf Lemma \ResultSeqTaut\ implies that
$2^{d-1}H_1(\cover M)\leq 
\beta\bigl(\bigoplus_{H\in\indexset^\star}H_1(M_H)\bigr)$,
so it is enough to show that
$\pi_H\shriek(H_1(M_H))\leq 2^{d-1}H_1(\cover M)$ for all
$H\in\indexset^\star$. If $\Gamma_H$ is connected, $H_1(M_H)$
has odd order and there is nothing to prove. 
If $\Gamma_H$ is disconnected and $\pi_H$ is factored as above
then the image of the transfer $H_1(M_i)\to H_1(M_{i+1})$
($1\leq i < d$) on integer homology is contained in 
$2H_1(M_{i+1})$, and the result follows.
\qed

\section{The mod 2 homology of 2-- and 4--fold branched covers}

In this section, the coefficients for homology will always 
be $\zt$, and will be omitted from the notation.
Let $L$ be a link in a connected, orientable 3--manifold
$N$. There is a double cover of $N$ with branch set $L$
iff $L$ represents zero in $H_1(N)$; suppose this is so.
Let $\maps{\theta}{H_1(N-L)}{\zt}$ be a homomorphism sending
each meridian of $L$ to $1$, and let $\maps{p}{\cover N}{N}$
be the corresponding branched cover. We wish to allow
the possibility that $L$ is empty (i.e., $p$ is unbranched);
in this case we insist that $\theta$ be onto, so that
$\cover N$ is connected.
There is an intersection pairing 
$H_1(N-L)\times H_2(N,L)\to \zt$ inducing an isomorphism
$H_2(N,L)\to \hom(H_1(N-L),\zt)$; we let
$\theta'\in H_2(N,L)$ correspond to $\theta$.
There is also an intersection pairing 
$H_2(N)\times H_2(N,L) \to H_1(N,L)$,
and we let $\maps{\theta''}{H_2(N)}{H_1(N,L)}$ 
be given by intersection with $\theta'$.

The transfer map with $\zt$ coefficients,
$\maps{p\shriek}{C(N)}{C(\cover N)}$, kills
$C(L)$, so there is an induced map
$\maps{p\shriek}{C(N,L)}{C(\cover N)}$. More generally,
if $X$ is any subcomplex of $N$ and $\cover X = p^{-1}(X)$, 
there is a map
$\maps{p\shriek}{C(N,L\cup X)}{C(\cover N,\cover X)}$
It was observed by Lee and Weintraub \cite{\LeeWein, Theorem 1}
that the sequence
$$0 \to C(N,L\cup X) \mapright{p\shriek} C(\cover N, \cover X)
\mapright{p} C(N,X) \to 0\eqno{(\label\LinkCoverSeq)}$$
is exact. When $N$ is a $\zt$ homology sphere,
it follows (taking $X=\emptyset$)
that $\maps{p\shriek}{H_1(N,L)}{H_1(\cover N)}$
is an isomorphism, which gives a different proof of
Sublemma 15.4 of \cite\Sakuma, that
$\dim H_1(\cover N) = b_0(L)-1$. The following
lemma generalizes this to other manifolds.

\proc{Lemma \label\LinkCover} In the above situation,
let $n=\dim H_1(N)$, let $r$ be the rank of the map
$H_1(L)\to H_1(N)$ induced by inclusion, and let $s$
be the rank of $\theta''$.
Then $r\leq s\leq n$, 
$\dim H_1(\cover N) = b_0(L)-1 +2n - r - s$,
and the rank of the map
$\maps{p\shriek}{H_1(N)}{H_1(\cover N)}$ equals $n-s$.
\endproc

Note a special case of this lemma: 
if $H_1(L)\to H_1(N)$ is onto then
$\dim H_1(\cover N)=b_0(L)-1$ and 
$\maps{p\shriek}{H_1(N)}{H_1(\cover N)}$ is the zero map.

\prf Certainly $s\leq \dim H_2(N) = n$. To see that $r\leq s$,
consider the composite of $\theta''$ and the connecting
homomorphism $H_1(N,L)\to H_0(L)$. 
This is the map $H_2(N)\to H_0(L)$ given by intersection with $L$,
which is dual to $H_1(L) \to H_1(N)$; therefore it has rank $r$.

From (\LinkCoverSeq) with $X=\emptyset$ we get an exact sequence
$$H_2(\cover N)\mapright{p} H_2(N)
\mapright{\partial} H_1(N,L) \to H_1(\cover N)
\to H_1(N) \to H_0(N,L) \to 0.$$
We claim that the connecting homomorphism labelled $\partial$ in this
sequence is equal to $\theta''$. We may take a 
(possibly non-orientable) surface $F$ in $N$ with boundary $L$
representing $\theta'\in H_2(N,L)$. Then $\cover N$ may be constructed
by gluing together two copies of $N$ cut open along $F$. Let 
$x \in H_2(N)$, and represent $x$ by a surface $F'$ transverse to
$F$. Then $p^{-1}(F')$ is the union of two copies of $F'$
cut open along $F\cap F'$. Either one of these carries
a 2--chain mapping to $F'$ under $p$, and their common
boundary is the image under $p\shriek$ of the element
of $C_1(N,L)$ carried by $F\cap F'$. Therefore $\partial(x)$
is represented by $F\cap F'$, which represents $\theta''(x)$
by the definition of $\theta''$, and the claim is proved. It
follows that $\dim H_1(\cover N) = \dim H_1(N,L)+n-s-\dim H_0(N,L)$,
and that the map $\maps{p}{H_2(\cover N)}{H_2(N)}$ has rank
$n-s$.
Now the exact sequence
$$ H_1(L)\to H_1(N) \to H_1(N,L) \to H_0(L) \to \zt \to H_0(N,L)\to 0$$
shows that $\dim H_1(N,L) = b_0(L)-1+n-r+\dim H_0(N,L)$,
so $\dim H_1(\cover N)$ is as claimed.
Also, the map $\maps{p\shriek}{H_1(N)}{H_1(\cover N)}$
is dual to $\maps{p}{H_2(\cover N)}{H_2(N)}$, and so has
rank $n-s$.
\qed

Now let $\Gamma$ be a $G(2)$--colored graph embedded in
a $\zt$ homology sphere $M$. Just as when $M$ is an
integral homology sphere, this determines a branched covering
$\maps{\pi}{\cover M}{M}$. We let the non-trivial elements of
$G$ be $g_1$, $g_2$ and $g_3$, and set 
$H_i=\genby{g_i}\in\indexset^\star$. Where $H_i$ would appear 
as a subscript, we just use $i$; thus we have 2--fold covers
$\maps{\rho_i}{M_i}{M}$ branched over $\Gamma_i$ and
$\maps{\pi_i}{\cover M}{M_i}$ branched over $\Delta_i$ 
for $1\leq i\leq 3$. If $\cover\Gamma=\pi^{-1}(\Gamma)$,
the map $\maps\pi{H_1(\cover M - \cover \Gamma)}{H_1(M - \Gamma)}$
kills each meridian of $\cover \Gamma$, so it induces a map
$\maps{\bar \pi}{H_1(\cover M)}{H_1(M - \Gamma)}$.

We wish to determine $\dim H_1(\cover M)$.
We deal first with the case where $\Gamma$ is connected,
since here we need some additional information.

\proc{Lemma \label\ModTwoGConn} When $\Gamma$ is connected,
$\dim H_1(\cover M)=b_1(\Gamma)-2$ and 
$\maps{\pi_i\shriek}{H_1(M_i)}{H_1(\cover M)}$ is the zero map for
$1\leq i\leq 3$. Further, the map
$\maps{\bar \pi}{H_1(\cover M)}{H_1(M - \Gamma)}$ is injective.
\endproc

\prf Let $1\leq i\leq 3$. Since $M$ is a $\zt$ homology sphere,
$\maps{\rho_i\shriek}{H_1(M,\Gamma_i)}{H_1(M_i)}$ is an isomorphism.
Since $\Gamma$ is connected,
every element of $H_1(M,\Gamma_i)$ is represented by a chain
of $\Gamma\setminus\Gamma_i$; since $\Delta_i$ is the inverse
image of $\Gamma\setminus\Gamma_i$ in $M_i$,
the map $H_1(\Delta_i)\to H_1(M_i)$ induced by inclusion
is onto. Since $\Delta_i$ is a link of $b_1(\Gamma)-1$
components, the first two claims follow from the special case of
Lemma \LinkCover\ noted above. The image of $\bar\pi$ is the
kernel of the homomorphism $H_1(M-\Gamma)\to G$ corresponding
to $\pi$; since this kernel has the same dimension as $H_1(\cover M)$,
it follows that $\bar \pi$ is injective.
\qed

Now let the components of $\Gamma$ be $\Gamma^k$ for
$1\leq k\leq b_0(\Gamma)$. We let $A=\{1,\ldots,b_0(\Gamma)\}$
be the index set for these components. For
$1\leq i\leq 3$, we partition $A$ into two sets $A_i$ and $A_i'$,
with $k\in A_i$ iff $\Gamma^k$ is a circular edge colored
$g_i$.
We also set $\Gamma^k_i = \Gamma^k\cap\Gamma_i$. 
(If $\Gamma^k$ is a circular edge, then $\Gamma^k_i$ is empty
if $\Gamma^k$ has color $g_i$, and equal to $\Gamma^k$ otherwise.)
If $\gamma$ is a 1--cycle of
$M-\Gamma^k$, we have $\sum_{i=1}^3\lk(\gamma,\Gamma^k_i)=0$,
where $\lk$ denotes mod 2 linking number. Hence, for $k\neq l$,
$$\eqalign{
\lk(\Gamma^k_1,\Gamma^l_2)+\lk(\Gamma^k_2,\Gamma^l_1)
&=\bigl(\lk(\Gamma^k_2,\Gamma^l_2)+\lk(\Gamma^k_3,\Gamma^l_2)\bigr)+
\bigl(\lk(\Gamma^k_2,\Gamma^l_2)+\lk(\Gamma^k_2,\Gamma^l_3)\bigr)\cr
&=\lk(\Gamma^k_2,\Gamma^l_3)+\lk(\Gamma^k_3,\Gamma^l_2),}$$
and similarly
$\lk(\Gamma^k_2,\Gamma^l_3)+\lk(\Gamma^k_3,\Gamma^l_2)
=\lk(\Gamma^k_3,\Gamma^l_1)+\lk(\Gamma^k_1,\Gamma^l_3)$;
we let $\lambda_{kl}\in \zt$ be this common value.
Note that if $k\in A_i$
then $\lambda_{kl} = \lk(\Gamma^k, \Gamma^l_i)$, and if also
$l\in A_j$ then
$\lambda_{kl}$ equals $\lk(\Gamma^k, \Gamma^l)$ if $i\neq j$
and $0$ if $i=j$. We also set 
$\lambda_{kk}=\sum_{l\in A, l\neq k}\lambda_{kl}$,
and let $\Lambda$ be the symmetric matrix 
$[\lambda_{kl}]_{k,l\in A}$.

\proc{Lemma \label\ModTwoGDisconn} We have
$\dim H_1(\cover M)=b_0(\Gamma)+b_1(\Gamma)-3 - \rank\Lambda$.
\endproc

\prf We shall prove this by applying Lemma \LinkCover\ to
the covering $\maps{\pi_1}{\cover M}{M_1}$. First we establish some 
notation.
\item{(a)} If $A'$ and $A''$ are subsets of $A$, we let
$\Lambda(A',A'')$ be the submatrix 
$[\lambda_{kl}]_{k\in A', l\in A''}$ of $\Lambda$.
Note that 
$\Lambda(A_1,A_1)$ is a diagonal matrix with diagonal
entries $\lambda_{kk}=\lk(\Gamma^k,\Gamma_1)$ for $k\in A_1$.

\item{(b)} We let $F$ be a surface in $M$ with 
$\partial F=\Gamma_1$.
Then $M_1$ can be constructed by gluing together
two copies of $M$ cut open along $F$.

\item{(c)} We denote the connecting homomorphisms in the exact 
sequences of the pairs $(M,\Gamma)$, $(M,\Gamma_1)$ and 
$(M,\Gamma\setminus\Gamma_1)$ by 
$\maps{\partial_i}{H_{i+1}(M,\Gamma)}{\tilde H_i(\Gamma)}$,
$\maps{\partial'_i}{H_{i+1}(M,\Gamma_1)}{\tilde H_i(\Gamma_1)}$
and
$\maps{\partial''_i}{H_{i+1}(M,\Gamma\setminus\Gamma_1)}
{\tilde H_i(\Gamma\setminus\Gamma_1)}$. (Here $\tilde H$
denotes reduced homology.) We need these maps only for
$i = 0$ or $1$, where they are isomorphisms.

\item{(d)} The case of (\LinkCoverSeq) for the cover $M_1 \to M$
with $X=\emptyset$ is
$$0 \to C(M,\Gamma_1) \mapright{\rho_1\shriek} C(M_1)
\mapright{\rho_1} C(M) \to 0.$$
The long exact sequence shows that
$$\maps{\alpha_i=\rho_1\shriek(\partial'_i)^{-1}}
{\tilde H_i(\Gamma_1)}{H_{i+1}(M_1)}$$
is an isomorphism for $i=0$ and an epimorphism for $i=1$.

\item{(e)} The case of (\LinkCoverSeq) for $M_1 \to M$
with $X=\Gamma\setminus\Gamma_1$ is
$$0\to C(M,\Gamma) \mapright{\rho_1\shriek} C(M_1,\Delta_1)
\mapright{\rho_1} C(M,\Gamma\setminus\Gamma_1) \to 0.$$
Denote the connecting homomorphisms in the long exact sequence
by $$\maps{\partial'''_i}{H_{i+1}(M,\Gamma\setminus\Gamma_1)}
{H_i(M, \Gamma)}.$$
We get an exact sequence
$$\eqalign{H_1(\Gamma)\mapright{\gamma_1}H_2(M_1,\Delta_1)
\mapright{\delta_1}H_1(\Gamma\setminus\Gamma_1)
&\mapright{\beta}\tilde H_0(\Gamma)\cr
&\mapright{\gamma_0}
H_1(M_1,\Delta_1)\mapright{\delta_0}
\tilde H_0(\Gamma\setminus\Gamma_1)\to 0,\cr}$$
where $\beta = \partial_0\partial'''_1(\partial''_1)^{-1}$,
$\gamma_i=\rho_1\shriek\partial_i^{-1}$ and
$\delta_i=\partial''_i\rho_1$.

\smallskip\noindent
The isomorphism 
$\maps{\alpha_0}{\tilde H_0(\Gamma_1)}{H_1(M_1)}$
gives
$$\dim H_1(M_1)=b_0(\Gamma_1)-1.\eqno{(\label\tempone)}$$
Next we determine $b_0(\Delta_1)$ 
and the rank of the map $H_1(\Delta_1)\to H_1(M_1)$.
Consider a non-circular edge $e$ of $\Gamma$ with color $g_1$;
the number of such edges is $b_1(\Gamma)-b_0(\Gamma)$,
and the inverse image $\rho_1^{-1}(e)$ is a single component of 
$\Delta_1$. The image under $\alpha_0^{-1}$ of the homology
class of $\rho_1^{-1}(e)$ is represented by $\partial e$,
and the subspace of $\tilde H_0(\Gamma_1)$ spanned by such
elements is 
$\hat H_0(\Gamma_1)= \bigoplus_{k\in A_1'}\tilde H_0(\Gamma^k_1)$,
which has dimension $b_0(\Gamma_1)-|A_1'|$. The remaining
components of $\Gamma\setminus\Gamma_1$ are the $\Gamma^k$
for $k\in A_1$, and such a $\Gamma^k$ is covered by a single
component of $\Delta_1$ if $\lambda_{kk}=1$,
and by two components if $\lambda_{kk}=0$. The number
of $k\in A_1$ with $\lambda_{kk}=1$ is the rank of the diagonal
matrix $\Lambda(A_1,A_1)$, and so
$$b_0(\Delta_1)=b_1(\Gamma)-b_0(\Gamma)+2|A_1|
-\rank \Lambda(A_1,A_1).\eqno{(\label\temptwo)}$$
For $k\in A_1$ with $\lambda_{kk}=1$, the component of
$\Delta_1$ covering $\Gamma^k$ is null-homologous.
Let $B$ be the set of $k\in A_1$ with $\lambda_{kk}=0$,
and $k\in B$. The two components
of $\Delta_1$ covering $\Gamma^k$ represent the same element
of $H_1(M_1)$, which we call $x^1_k$. We may assume that the
surface $F$ is disjoint from $\Gamma^k$, and
take a surface $F'$ with boundary
$\Gamma^k$ which is transverse to $F$.
Then $\rho_1^{-1}(F')$ is the union of two copies of
$F'$ cut open along $F\cap F'$, either one of which
shows that $x^1_k$ is the image under $\rho_1\shriek$ of the element
of $H_1(M,\Gamma_1)$ represented by $F\cap F'$. 
Thus $\alpha_0^{-1}(x^1_k)$ is represented by $\partial(F\cap F')$.
Now $H_0(\Gamma_1)/\hat H_0(\Gamma_1)$ has a basis with one element
$x^0_l$ for each $l\in A_1'$, and the image of $\alpha_0^{-1}(x^1_k)$
in this quotient is $\sum_{l\in A_1'}\lambda_{kl}x^0_l$.
Therefore the rank of $H_1(\Delta_1) \to H_1(M_1)$ is
$\dim \hat H_0(\Gamma_1) + \rank \Lambda(B,A_1')$.
But $\rank \Lambda(B,A_1') = \rank\Lambda(A_1,A)-\rank\Lambda(A_1,A_1)$, so
$$\rank\bigl(H_1(\Delta_1) \to H_1(M_1)\bigr)
=b_0(\Gamma_1)-|A_1'|+\rank\Lambda(A_1,A)-\rank\Lambda(A_1,A_1).
\eqno{(\label\tempthree)}$$

The 2--fold covering $\maps{\pi_1}{\cover M}{M_1}$ corresponds to
a homomorphism $\maps{\theta}{H_1(M_1-\Delta_1)}{\zt}$,
to which are associated $\theta'\in H_2(M_1,\Delta_1)$ and
$\maps{\theta''}{H_2(M_1)}{H_1(M_1,\Delta_1)}$; we must determine
the rank of $\theta''$. We first identify $\theta'$.
For $x\in H_1(M-\Gamma)$, $\lk(x,\Gamma_i)$ is well-defined for
$1\leq i\leq 3$, and $\sum_{i=1}^3 \lk(x,\Gamma_i)=0$.
Define a homomorphism $\maps{\phi}{H_1(M-\Gamma)}{G}$
by $\phi(x)=\prod_{i=1}^3 g_i^{\lk(x,\Gamma_i)}$.
Then $\phi$ sends the meridian of an edge of $\Gamma$ to the color
of that edge, so it is the homomorphism corresponding to the cover
$\cover M \to M$.
Let $\hat\Gamma = \rho_1^{-1}(\Gamma)$, and let
$\maps{\iota}{H_1(M_1-\hat\Gamma)}{H_1(M_1-\Delta_1)}$
be the surjection induced by inclusion. For $y\in H_1(M_1-\hat\Gamma)$,
we have $\rho_1(y)\in H_1(M-\Gamma)$ and 
$\phi\rho_1(y)=g_1^{\theta\iota(y)}$. 
It follows that
$\lk(\rho_1(y), \Gamma_2)= \lk(\rho_1(y), \Gamma_3)$,
$\lk(\rho_1(y), \Gamma_1)=0$, and
$\theta\iota(y)=\lk(\rho_1(y), \Gamma_2)$.
There are intersection pairings
$H_1(M_1-\Delta_1)\times H_2(M_1,\Delta_1)\to \zt$
and $H_1(M-\Gamma)\times H_2(M,\Gamma)\to \zt$
and a linking pairing
$H_1(M-\Gamma)\times H_1(\Gamma)\to \zt$, and they are related by
$$\iota(y)\cdot\gamma_1(z) 
= \iota(y)\cdot\rho_1\shriek\partial_1^{-1}(z)
= \rho_1(y)\cdot \partial_1^{-1}(z)
= \lk(\rho_1(y),z)$$
for $y\in H_1(M_1-\hat\Gamma)$ and $z\in H_1(\Gamma)$.
For $k\in A$, let $z_k\in H_1(\Gamma)$ be represented
by $\Gamma^k_2$. Then $\sum_{k\in A}z_k$ is represented by
$\Gamma_2$, and so 
$\iota(y)\cdot\gamma_1\bigl(\sum_{k\in A}z_k\bigr)=\theta\iota(y)$
for $y\in H_1(M_1-\hat\Gamma)$. Therefore
$\theta' = \gamma_1\bigl(\sum_{k\in A}z_k\bigr)$.

We have an epimorphism 
$\maps{\alpha_1}{H_1(\Gamma_1)}{H_2(M_1)}$.
For $k\in A_1'$, let $y^1_k\in H_1(\Gamma_1)$ be
represented by $\Gamma^k_1$, and let $\hat H_1(\Gamma_1)$
be the subspace of $H_1(\Gamma_1)$ generated by these elements.
Also let $\hat\alpha_1$ be the restriction of $\alpha_1$
to $\hat H_1(\Gamma_1)$. We shall show that 
$\ker(\theta''\alpha_1)\leq \hat H_1(\Gamma_1)$, from which
it will follow that
$$\rank \theta'' = \rank(\theta''\alpha_1)=
\rank(\theta''\hat\alpha_1) +\dim H_1(\Gamma_1) 
- \dim\hat H_1(\Gamma_1),$$ or
$$\rank \theta'' = \rank(\theta''\hat\alpha_1) +b_0(\Gamma_1)-|A_1'|.
\eqno{(\label\tempfour)}$$
Consider the composite 
$\maps{\delta_0\theta''\alpha_1}{H_1(\Gamma_1)}
{\tilde H_0(\Gamma\setminus\Gamma_1)}$. This may be described geometrically
as follows. If $x\in H_1(\Gamma_1)$ is represented by a circuit $C$,
take surfaces $F'$ and $F''$ with $\partial F'=C$ and 
$\partial F''=\Gamma_2$ that meet transversely except along the common
part of their boundaries, $C\cap\Gamma_2$. Then the closure of
$(F'\cap F'')-(C\cap\Gamma_2)$ represents an element $y$ of
$H_1(M,\Gamma)$. Now $\theta''\alpha_1(x)\in H_1(M_1,\Delta_1)$
is represented by $\rho_1^{-1}(F'\cap F'')$, and is therefore
the sum of $\rho_1\shriek(y)$ and the element represented by
$\rho_1^{-1}(C\cap\Gamma_2)$. Since $\rho_1\rho_1\shriek(y)=0$,
$\rho_1\theta''\alpha_1(x)\in H_1(M,\Gamma\setminus\Gamma_1)$
is represented by $C\cap\Gamma_2$. Hence $\delta_0\theta''\alpha_1(x)$
is represented by $\partial(C\cap\Gamma_2)$,
which is just the sum of the vertices of $\Gamma$ lying on $C$.
It follows that $\delta_0\theta''\alpha_1(x)=0$ iff 
$x\in \hat H_1(\Gamma_1)$, so 
$\ker(\theta''\alpha_1)\leq \hat H_1(\Gamma_1)$, as claimed.

We let the elements of
the natural basis for $H_0(\Gamma)$ be $y^0_k$ for $k\in A$,
and define 
$\maps{\hat \beta}{\hat H_1(\Gamma_1)}{\tilde H_0(\Gamma)}$
by $\hat\beta(y^1_k)=\sum_{l\in A}\lambda_{kl}y^0_l$
for $k\in A_1'$. We claim that 
$\maps{\gamma_0\hat \beta=\theta''\hat\alpha_1}
{\hat H_1(\Gamma_1)}{H_1(M_1,\Delta_1)}$.
Let $k\in A_1'$ and $l\in A$. We have 
$\hat\alpha_1(y^1_k)\in H_2(M_1)$ and $\gamma_1(z_l)\in H_2(M_1,\Delta_1)$,
with intersection 
$\hat\alpha_1(y^1_k)\cdot\gamma_1(z_l) \in H_1(M_1,\Delta_1)$.
Suppose $k\neq l$.
Then $(\partial_1')^{-1}(y^1_k)\in H_2(M,\Gamma_1)$
and $\partial_1^{-1}(z_l)\in H_2(M,\Gamma)$ may be represented
by transverse surfaces $F'$ and $F''$ with boundaries
$\Gamma^k_1$ and $\Gamma^l_2$, respectively, and
$\hat\alpha_1(y^1_k)\cdot\gamma_1(z_l)$
is the image under
$\maps{\rho_1\shriek}{H_1(M,\Gamma)}{H_1(M_1,\Delta_1)}$
of the class represented by $F'\cap F''$. Since the image of
this class under $\maps{\partial_0}{H_1(M,\Gamma)}{\tilde H_0(\Gamma)}$
is $\lk(\Gamma^k_1,\Gamma^l_2)(y^0_k+y^0_l)$, we have
$$\hat\alpha_1(y^1_k)\cdot\gamma_1(z_l)=
\gamma_0\bigl(\lk(\Gamma^k_1,\Gamma^l_2)(y^0_k+y^0_l)\bigr)
\quad\hbox{for $k\in A_1'$, $l\in A$, $k\neq l$.}$$
Now $\sum_{k\in A_1'}(\partial_1')^{-1}(y^1_k)$ is represented
by $F$, whose inverse image in $M_1$ is null homologous,
so $\sum_{k\in A_1'}\hat\alpha_1(y^1_k)=0$. Therefore, for $k\in A_1'$,
$$\hat\alpha_1(y^1_k)\cdot\gamma_1(z_k)
=\sum_{l\in A_1', l\neq k}\hat\alpha_1(y^1_l)\cdot\gamma_1(z_k)
=\sum_{l\in A, l\neq k}
\gamma_0\bigl(\lk(\Gamma^l_1,\Gamma^k_2)(y^0_k+y^0_l)\bigr),
$$
where in the last term we may sum over $A$ since $\Gamma^l_1$ is empty
for $l\notin A_1'$. Hence, again for $k\in A_1'$,
$$\eqalign{\theta''\hat\alpha_1(y^1_k)&=\hat\alpha_1(y^1_k)\cdot\theta'
=\sum_{l\in A}\hat\alpha_1(y^1_k)\cdot\gamma_1(z_l)\cr
&=\sum_{l\in A, l\neq k}\gamma_0\bigl(\lambda_{kl}(y^0_k+y^0_l)\bigr)
=\sum_{l\in A}\gamma_0(\lambda_{kl}y^0_l)
=\gamma_0\hat\beta(y^1_k),
}$$
and so indeed $\gamma_0\hat \beta=\theta''\hat\alpha_1$. Thus we
have a commutative diagram
{\abovedisplayskip=\diagdisplayskip  
\belowdisplayskip=\diagdisplayskip   
$$
\def\normalbaselines{\diagbaselines} 
\matrix{
&&\hat H_1(\Gamma_1)&\mapright{\hat\alpha_1}&H_2(M_1)\cr
&&\mapdown{\hat\beta}&&\mapdown{\theta''}&&\cr
H_1(\Gamma\setminus\Gamma_1)&\mapright{\beta}&\tilde H_0(\Gamma)
&\mapright{\gamma_0}&H_1(M_1,\Delta_1)&&
}$$}%
in which the bottom row is exact. Therefore
$$\rank(\theta''\hat\alpha_1) = \rank(\gamma_0\hat\beta)
=\dim(\im\beta+\im\hat\beta) - \rank\beta.$$
For $k\in A_1$, $\Gamma^k$ represents an element $y^1_k$ of
$H_1(\Gamma\setminus\Gamma_1)$, and these form a basis.
We claim that $\beta(y^1_k)=\sum_{l\in A}\lambda_{kl}y^0_l$
for $k\in A_1$, from which it will follow that
$$\rank(\theta''\hat\alpha_1) = \rank\Lambda - \rank\Lambda(A_1,A).
\eqno{(\label\tempfive)}$$
For $k\in A_1$, $(\partial''_1)^{-1}(y^1_k)$ is represented
by a surface $F'$ with boundary $\Gamma^k$, which we may
take to be transverse to $F$. 
Then $\rho_1^{-1}(F')$ is the union of two copies of
$F'$ cut open along $F\cap F'$. The boundary of either
one is the union of $\rho_1^{-1}(F\cap F')$ and part of
$\rho_1^{-1}(\Gamma^k)\subseteq \Delta_1$, so it represents
the same element of $C_1(M_1,\Delta_1)$ as
$\rho_1^{-1}(F\cap F')$. It follows that
$\partial'''_1(\partial''_1)^{-1}(y^1_k)$ is represented by
$F\cap F'$, and hence that
$$\beta(y^1_k)
=\partial_0\partial'''_1(\partial''_1)^{-1}(y^1_k)
=\sum_{l\in A_1'}\lk(\Gamma^k,\Gamma^l_1)(y^0_l+y^0_k)
=\sum_{l\in A}\lambda_{kl}y^0_l,$$
as claimed.

The proof of the lemma is completed by applying Lemma \LinkCover\ to
the covering $\cover M\to M_1$ and using the equations
(\tempone) -- (\tempfive).
\qed

In applying Lemma \ModTwoGDisconn, we compute the matrix $\Lambda$
using the following result, which is implicit in the proof
of Lemma 1 of Flapan \cite\Flapan.

\proc{Lemma \label\CalcLinking} Let $K$ be a knot in
a $\zt$ homology 3--sphere $N$, and let $A$ and $B$ be disjoint
arcs in $N$ meeting $K$ in their endpoints. Let $\cover N$
be the 2--fold cover of $N$ branched over $K$, and let
$\cover A$ and $\cover B$ be the inverse images of $A$ and $B$
in the $\zt$ homology sphere $\cover N$. Then 
$\lk(\cover A, \cover B)=1$ iff the endpoints of $A$ separate those of
$B$ on $K$.
\qed

\section{Proofs of theorems}

Recall that in the statement of each theorem,
$\Gamma$ is a $G(d)$--colored graph embedded in a homology 3--sphere 
$M$, with corresponding branched cover $\cover M$.

\proc{Theorem \checklabel\ThmDTwo}  If $d=2$ and $\Gamma$ is connected, then
there is a short exact sequence
$$0 \to \bigoplus_{H\in\indexset^\star}H_1(M_H)
\mapright\beta H_1(\cover M) \to \ztb^{b_1(\Gamma)-2}\to 0,$$
and 
$\beta\bigl(\bigoplus_{H\in\indexset^\star}H_1(M_H)\bigr)
=2H_1(\cover M)$.
\endproc

\prf Lemma \ResultSeqTaut\ gives the exact sequence,
while Lemma \ModTwoGConn\ shows that the mod 2 transfer 
$\maps{\pi_H\shriek}{H_1(M_H;\zt)}{H_1(\cover M;\zt)}$
is zero for $H\in\indexset^\star$, which implies
the second assertion by Lemma \Extension.
\qed

\proc{Theorem \checklabel\ThmDThreeMCirc} 
Let $d=3$ and let $\Gamma$ have an unsplittable coloring with 
a special $m$--circuit. Then $3\leq m\leq b_1(\Gamma)$,
there is a short exact sequence
$$0 \to \bigoplus_{H\in\indexset^\star}H_1(M_H)
\mapright\beta H_1(\cover M) \to 
\zmodb{4}^{m-3}\oplus\ztb^{2(b_1(\Gamma)-m)}\to 0,$$
and 
$\beta\bigl(\bigoplus_{H\in\indexset^\star}H_1(M_H)\bigr)
=4H_1(\cover M)$.
\endproc

\prf Let $H_0\in \indexset^\star$ be such that 
$\Gamma_0=\Gamma_{H_0}$ is a special $m$--circuit,
and let $M_0=M_{H_0}$ and $\Delta_0=\Delta_{H_0}$.
Since the coloring is unsplittable, $\Gamma$ is simple,
so any circuit has length at least 3. Further,
$\Gamma$ is connected, so $\chi(\Gamma\setminus\Gamma_0)
=1-b_1(\Gamma)+m$; since $\Gamma\setminus\Gamma_0$
is connected, this gives $m\leq b_1(\Gamma)$.
By Lemma \MCircIsTaut, $\Gamma$ is taut, so Lemma
\ResultSeqTaut\ gives an exact sequence
$$0 \to \bigoplus_{H\in\indexset^\star}H_1(M_H)
\mapright\beta H_1(\cover M) \to 
\zmodb{4}^a\oplus\ztb^b\to 0$$
for some $a$ and $b$ with $2a+b=2b_1(\Gamma)-6$.

Suppose $1\neq h\in H\in\indexset^\star$. There
is a cover $\cover M/\genby{h}\to M$ with group
$H/\genby{h}\cong \ztb^2$; its branch set is obtained
from $\Gamma$ by deleting all edges with color $h$.
Since $\Gamma$ is unsplittable, we may apply
Lemma \ModTwoGConn\ to this cover to show
that the transfer $H_1(M_H;\zt)\to H_1(\cover M/\genby{h};\zt)$
is zero. By Lemma \Extension, to show that
$\beta\bigl(\bigoplus_{H\in\indexset^\star}H_1(M_H)\bigr)
=4H_1(\cover M)$ it is then enough to show that, for each 
$H\in\indexset^\star$, there is some non-trivial $h$ in $H$
such that $H_1(\cover M/\genby{h};\zt)\to H_1(\cover M;\zt)$
is zero. Now consider the cover $\cover M \to M_0$,
with group $H_0\cong \ztb^2$ and branch set $\Delta_0$.
Since $\Gamma_0$ is a circuit, $M_0$ is a $\zt$ homology sphere.
Since $\Gamma\setminus\Gamma_0$ is connected, so is
$\Delta_0$, and Lemma \ModTwoGConn\ applies to this cover,
showing that $H_1(\cover M/\genby{h};\zt)\to H_1(\cover M;\zt)$
is zero whenever $1\neq h\in H_0$. Since $H\cap H_0$ contains
a non-trivial element for all $H\in\indexset^\star$,
the proof that $\beta\bigl(\bigoplus_{H\in\indexset^\star}H_1(M_H)\bigr)
=4H_1(\cover M)$ is complete. It follows that
$H_1(\cover M;\zt)\cong H_1(\cover M)/2H_1(\cover M)\cong
\zt^{a+b}$. On the other hand, Lemma \ModTwoGConn\ applied to
$\cover M \to M_0$ also shows that
$\dim H_1(\cover M;\zt)=b_1(\Delta_0)-2$. Since
$b_1(\Delta_0)=2b_1(\Gamma)-m-1$, we have
$a+b=2b_1(\Gamma)-m-3$. It follows that $a=m-3$ and
$b= 2(b_1(\Gamma)-m)$, and we are done.
\qed

\proc{Theorem \checklabel\ThmDThreeMob} 
Let $\Gamma$ be an $n$--rung M\"obius ladder ($n\geq 2$) with
a $G(3)$-coloring, and let $g_0$ be the product of the colors on the rungs.
Suppose that $g_0\neq 1$, and let $k$ be the number of rungs with
color $g_0$. If $k=0$, there is a short exact sequence
$$0 \to \bigoplus_{H\in\indexset^\star}H_1(M_H)
\mapright\beta H_1(\cover M) \to 
\zmodb{4}^{n-2}\to 0,$$
while if $k>0$ 
there is a short exact sequence
$$0 \to \bigoplus_{H\in\indexset^\star}H_1(M_H)
\mapright\beta H_1(\cover M) \to 
\zmodb{4}^{n-k-1}\oplus\ztb^{2(k-1)}\to 0.$$
In either case,
$\beta\bigl(\bigoplus_{H\in\indexset^\star}H_1(M_H)\bigr)
=4H_1(\cover M)$.
\endproc

\prf By Lemma \MobIsTaut, $\Gamma$ is taut, so Lemma
\ResultSeqTaut\ gives an exact sequence
$$0 \to \bigoplus_{H\in\indexset^\star}H_1(M_H)
\mapright\beta H_1(\cover M) \to 
\zmodb{4}^a\oplus\ztb^b\to 0$$
for some $a$ and $b$ with $2a+b=2n-4$.
Consider the cover $\maps{\pi'}{\cover M/\genby{g_0}}{M}$
with group $G'=G/\genby{g_0}\cong \ztb^2$. Its branch
set is the $(n-k)$--rung M\"obius ladder $\Gamma'$ obtained by 
deleting the rungs colored $g_0$, so Lemma
\ModTwoGConn\ shows that $H_1(M_H;\zt)\to H_1(\cover M/\genby{g_0};\zt)$
is zero whenever $g_0\in H\in\indexset^\star$.
By Lemma \EdgeParity, $\Gamma_H$ is connected if $g_0\notin H$,
so Lemma \Extension\ will imply that 
$\beta\bigl(\bigoplus_{H\in\indexset^\star}H_1(M_H)\bigr)
=4H_1(\cover M)$ provided that
$H_1(\cover M/\genby{g_0};\zt) \to H_1(\cover M;\zt)$
is also zero. Lemma \ModTwoGConn\ also gives
$\dim H_1(\cover M/\genby{g_0};\zt)= n-k-1$. The 2--fold cover
$\cover M\to \cover M/\genby{g_0}$ has as branch set
a link $L$, which is the inverse image of the rungs
of $\Gamma$ labelled $g_0$. Let $r$ be the rank of
$H_1(L;\zt)\to H_1(\cover M/\genby{g_0};\zt)$.
If $k=0$ then $L$ is empty and $r=0$.
Suppose $k>0$, and consider a rung $e$ labelled $g_0$. The endpoints
of $e$ lie on two edges of $\Gamma'$ with the same color in the
$G'$--labelling determining $\pi'$. Hence,
if $D\subset M$ is a 2--disk containing $e$ in its interior
and meeting $\Gamma'$ only in the endpoints of $e$, then
$(\pi')^{-1}(D)$ consists of two annuli. This shows first that
$(\pi')^{-1}(e)$ has two components, so $b_0(L)=2k$. It also
shows that under the map 
$\maps{\bar\pi'}{H_1(\cover M/\genby{g_0};\zt)}{H_1(M-\Gamma';\zt)}$,
each component of $(\pi')^{-1}(e)$ is sent to the element
of $H_1(M-\Gamma';\zt)$ represented by $\partial D$. This element
is non-trivial and independent of the choice of $e$. By Lemma 
\ModTwoGConn, $\bar\pi'$ is injective, and it follows that $r=1$.
Using the Kronecker delta, we may say that in all cases
$b_0(L)=2k$ and $r=1-\delta_{k0}$. 
It now follows from Lemma \LinkCover\ applied to the 
cover $\cover M\to \cover M/\genby{g_0}$ that
$$\dim H_1(\cover M;\zt)-
\rank\bigl(H_1(\cover M/\genby{g_0};\zt) \to H_1(\cover M;\zt)\bigr)
=n+k-3+\delta_{k0}\eqno{(\label\temp)}.$$

Now choose $H\in\indexset^\star$ with $g_0\notin H$. Then $M_H$
is a $\zt$ homology sphere, and we may compute $\dim H_1(\cover M;\zt)$
by applying Lemma \ModTwoGDisconn\ to the cover
$\maps{\pi_H}{\cover M}{M_H}$. We must compute the matrix $\Lambda$
of that lemma. Suppose that $\Gamma_H$ contains the $m$ rungs
$\tau_{i_j}$ for $0\leq j<m$, where $0\leq i_0<\ldots<i_{m-1}<n$,
and let the color of $\tau_{i_j}$ be $h_j\in G-H$. Then 
$\Gamma\setminus\Gamma_H$ has $m$ components $C_0,\ldots,C_{m-1}$,
and the components of $\Delta_H$ are 
$\rho_H^{-1}(C_0),\ldots,\rho_H^{-1}(C_{m-1})$. We may choose
the numbering of the $C_j$ so that $\tau_{i_j}$ has one vertex
on $C_j$ and the other on $C_{j+1}$. (The subscripts on the $C_j$
are to be taken modulo $m$.) By Lemma \CalcLinking, all the off-diagonal
elements of $\Lambda$ except $\lambda_{j,j\pm 1}$ are zero.
Also, if the edges of $C_j$ and $C_{j+1}$ that meet $\tau_j$
have colors $h_j'$ and $h_j''$, then $\lambda_{j,j+ 1}=1$ iff
$h_j'\neq h_j''$. However, $h_j'h_j''=g_0h_j$, 
so $\lambda_{j,j+ 1}=1$ iff $h_j\neq g_0$. Since exactly $k$
of the $h_j$ are equal to $g_0$, it follows that
$\rank\Lambda = m-k-\delta_{k0}$.

The trivalent graph $\Delta_H$ has $2(n-m)$ vertices,
so $\chi(\Delta_H)=m-n$, and since $b_0(\Delta_H)=m$
we have $b_1(\Delta_H)=n$. Now Lemma \ModTwoGDisconn\ gives
$\dim H_1(\cover M;\zt) = n+k-3+\delta_{k0}$.
Comparing this to (\temp), we see that 
$H_1(\cover M/\genby{g_0};\zt) \to H_1(\cover M;\zt)$
is the zero map, and hence
$\beta\bigl(\bigoplus_{H\in\indexset^\star}H_1(M_H)\bigr)
=4H_1(\cover M)$.
It then follows that 
$a+b=\dim H_1(\cover M;\zt) = n+k-3+\delta_{k0}$,
giving $a=n-k-1-\delta_{k0}$ and $b=2(k-1+\delta_{k0})$,
completing the proof.
\qed

Suppose that $\Gamma$ is taut. By Lemmas
\ResultSeqOne\ and \ResultSeqTwo,
we may identify $\bigoplus_{H\in\indexset^\star}H_1(M_H)$
and $H_1(\im\beta)$ with their images in $H_1(\cover M)$;
thus 
$\bigoplus_{H\in\indexset^\star}H_1(M_H) \leq
H_1(\im\beta) \leq H_1(\cover M)$. For any $x\in H_1(\cover M)$,
$2^{d-1}x \in \bigoplus_{H\in\indexset^\star}H_1(M_H)$.
In the proofs of the remaining theorems, we need to show
that we may choose $x$ so that
$2^{d-2}x \notin \bigoplus_{H\in\indexset^\star}H_1(M_H)$.
Now $2^{d-2}x$ is in $H_1(\im\beta)$, and so it is in
$\bigoplus_{H\in\indexset^\star}H_1(M_H)$
iff it is in the kernel
of the map $H_1(\im\beta) \to \ztb^{b_1(\Gamma)-d}$ from Lemma 
\ResultSeqOne.
From the proof of that lemma, the kernel
of this map is equal to the kernel of the composite
of the maps $H_1(\im\beta)\to H_0(\ker\beta)$ from the long
exact sequence of (\ShortSeqTwo), and 
$H_0(\ker\beta)\to H_0(\Gamma\sep d-1)$ from Lemma \HKerB.

\proc{Lemma \label\HighOrder} Suppose that
$\Gamma$ is taut, and let $e_1,\ldots,e_n$ be edges of 
$\Gamma$ with colors $g_1,\ldots,g_n$ such that $g_1\cdots g_n=1$.
For $1\leq i \leq n$, pick a vertex $v_i$ of $e_i$. Then
there is an element $x$ of $H_1(\cover M)$ such that the image
of $2^{d-2}x$ in $H_0(\Gamma\sep d-1)$ is represented by
$\sum_{i=1}^n\sum_{H\in\indexset^\star}\achar_H(g_i)v_i H
\in C'_0(\Gamma\sep d-1)$.
\endproc

\prf Consider an element $x$ of $H_1(\cover M)$
represented by a cycle of the form
$z = \sum_{\sigma\in S_1(M)}(1-h_\sigma)\tilde\sigma$
for some $h_\sigma \in G$. Let $S'$ be the set of those $\sigma$
for which $h_\sigma \neq 1$, and
for each $\sigma\in S'$, 
define an element $c_\sigma$ of 
$\sum_{H\in \indexset}C_1(M_H)$ by
$$c_\sigma = -2^{d-2}\sigma G 
+\sum_{H\in\indexset}\half\bigl(1-\mchar_H(h_\sigma)\bigr)\sigma_H H.$$
Then
$$\eqalign{
\beta(c_\sigma)&=
-2^{d-2}\sum_{g\in G}g\tilde\sigma 
+\sum_{H\in\indexset}\half\bigl(1-\mchar_H(h_\sigma)\bigr)
\sum_{h\in H} h\tilde\sigma\cr
&=\sum_{g\in G}\Bigl(-2^{d-2} + 
\sum_{H\in\indexset}\half\bigl(1-\mchar_H(h_\sigma)\bigr)
\half\bigl(1+\mchar_H(g)\bigr)\Bigr)g\tilde\sigma\cr
&=\sum_{g\in G}\sum_{H\in\indexset}\smallfrac{1}{4}
\bigl(\mchar_H(g)-\mchar_H(h_\sigma)-\mchar_H(gh_\sigma)\bigr)
g\tilde\sigma\cr
&= 2^{d-2}(1-h_\sigma)\tilde\sigma.
}$$
Therefore $\beta\bigl(\sum_{\sigma\in S'}c_\sigma\bigr) = 2^{d-2}z$,
and so the image of $2^{d-2}x$ in $H_0(\ker\beta)$ is represented
by $\sum_{\sigma\in S'}\partial c_\sigma$.  From the proofs
of Lemmas \HKerB\ and \KerBImA, the image
of $2^{d-2}x$ in $H_0(\Gamma\sep d-1)$ is represented by
$$z'=\sum_{\sigma\in S', H\in\indexset^\star}
\achar_H(h_\sigma)(\partial\sigma)H
= \sum_{\sigma\in S_1(M), H\in\indexset^\star}
\achar_H(h_\sigma)(\partial\sigma)H$$
(since $\half(1-\mchar_H(g)) \bmod 2=\achar_H(g)$).

We now construct a specific 1--cycle. Take a disc $D$ in $M$ meeting
$\Gamma$ transversely in $n$ points $p_1,\ldots,p_n$, where
$p_i$ lies on the edge $e_i$. Take disjoint arcs
$A_1,\ldots A_n$ on $D$, where $A_i$ joins $p_i$ to a point $q_i$
of $\partial D$ and $q_i$ is adjacent to $q_{i+1}$ on $\partial D$.
(Here and in the rest of the proof, subscripts are to be
taken modulo $n$.) We may assume that $D$ and each $A_i$ are triangulated
by subcomplexes of $M$ (and hence the $p_i$ and $q_i$ are 0--simplices
of $M$). Let $c_i\in C_1(M)$ be a 1--chain carried by $A_i$
with $\partial c_i = q_i-p_i$. Also let $d_i\in C_1(M)$ be carried
by one of the arcs into which the $q_i$ divide $\partial D$,
with $\partial d_i= q_{i+1}-q_i$. Let $\tilde c_i$ and $\tilde d_i$
be the images of $c_i$ and $d_i$ under the $\z$--module homomorphism
$C(M)\to C(\cover M)$ taking $\sigma$ to $\tilde\sigma$ ($\sigma
\in S(M)$). Now $\pi^{-1}(D)$ is the union of $2^d$ copies of $D$
cut open along the $A_i$; let $\cover D$ be one copy. If $\sigma$
is either $p_i$ or a 1--simplex of $\partial D$, there is just
one lift of $\sigma$ lying in $\partial \cover D$; we take this to be
$\tilde\sigma$. If $\sigma$
is either $q_i$ or a 1--simplex of $A_i$, there are
two lifts of $\sigma$ lying in $\partial \cover D$, and $g_i$
takes one to the other. We may choose $\tilde\sigma$ to be one of
these lifts in such a way that 
$\partial\tilde d_i = g_{i+1}\tilde q_{i+1}-\tilde q_i$
and $\partial\tilde c_i =\tilde q_i-\tilde p_i$. With these choices,
$z_1 = \sum_{i=1}^n(\tilde c_i -g_i\tilde c_i +\tilde d_i)$
is a 1--cycle of $\cover M$ carried by $\partial\cover D$.
Set $g_i'=\prod_{j=1}^i g_j$; $g_i'$ depends only on 
$i \bmod n$ since $g_1\cdots g_n=1$, and so
$z_2 = \sum_{i=1}^n g_i'\tilde d_i$ is another 1--cycle of $\cover M$.
(It is carried by a single lift of $\partial D$.) Let 
$x\in H_1(\cover M)$ be represented by
$z = z_1-z_2 = \sum_{i=1}^n\bigl((1-g_i)\tilde c_i 
+ (1-g_i')\tilde d_i\bigr)$. By the previous paragraph,
the image of $2^{d-2}x$ in $H_0(\Gamma\sep d-1)$ is represented by
$$\eqalignno{
z' &= \sum_{i=1}^n\sum_{H\in\indexset^\star}
\bigl(\achar_H(g_i)(q_i+p_i)+\achar_H(g_i')(q_{i+1}+q_i)\bigr)H\cr
&= \sum_{i=1}^n\sum_{H\in\indexset^\star}
\bigl(\achar_H(g_i)p_i 
+(\achar_H(g_i)+\achar_H(g_i')+\achar_H(g_{i-1}'))q_i\bigr)H\cr
&= \sum_{i=1}^n\sum_{H\in\indexset^\star} \achar_H(g_i)p_i H
\quad\hbox{(because $g_ig_i'g_{i-1}'=1$).}
}$$
Now
$\sum_{i=1}^n\sum_{H\in\indexset^\star} \achar_H(g_i)p_i H$
is homologous to
$\sum_{i=1}^n\sum_{H\in\indexset^\star} \achar_H(g_i)v_i H$,
and the proof is complete.
\qed

The next lemma will be used in the proof of Theorem \ThmDFour\ to
show that the 0--chain of
the previous one is not a boundary.

\proc{Lemma \label\MobParity} Let $\Gamma$ be an $m$--rung
M\"obius ladder, and let $0\leq i_1<i_2<\cdots<i_k<m$,
where either $m$ is odd and $k$ is even, or $m=k=2$.
Let the two circuits of $\Gamma$ that contain all the rungs
be $\Gamma_1$ and $\Gamma_2$, and for $\alpha=1$ or $2$, let
$c_\alpha\in C_1(\Gamma_\alpha;\zt)$ be such that
$\partial c_\alpha=\sum_{j=1}^k v_{i_j}$. Let $a\in\zt$
be the sum of the coefficients of the rungs $\tau_{i_j}$
in $c_\alpha$ for $1\leq j\leq k$ and $\alpha =1$ or $2$.
Then $a=1$ iff $m=2$.
\endproc

\prf Since each $\Gamma_\alpha$ contains all the rungs
$\tau_{i_j}$, $a$ is independent of the choice of
$c_1$ and $c_2$. Suppose first that $m=k=2$, and so
$i_1=0$ and $i_2=1$. If $\Gamma_1$ is taken to be the circuit
containing $\sigma_0$, we may take $c_1=\sigma_0$ and
$c_2=\tau_0+\sigma_1$, and so $a=1$.

Now suppose that $m$ is odd. We show that for $1\leq j\leq\half k$,
there are chains $c_{\alpha j}\in C_1(\Gamma_\alpha;\zt)$
with $\partial c_{\alpha j}=v_{i_{2j-1}}+v_{i_{2j}}$
such that 
$$c_{1j}+c_{2j} = \tau_{i_{2j-1}}+\tau_{i_{2j}}
+\sum_{i=i_{2j-1}}^{i_{2j}-1}(\sigma_i+\sigma_{i+n}),$$
the sum being taken in $C_1(\Gamma;\zt)$.
Then we may set $c_\alpha=\sum_{j=1}^{k/2}c_{\alpha j}$
and conclude that $a=0$. Given $j$, let $\Gamma_\alpha$
be that one of $\Gamma_1$ and $\Gamma_2$ that contains 
$\sigma_{i_{2j-1}}$,
and $\Gamma_\beta$ the other. If $i_{2j}-i_{2j-1}$ is odd, we may set
$$\eqalign{
c_{\alpha j} &= \sigma_{i_{2j-1}}+\tau_{i_{2j-1}+1}
+\sigma_{i_{2j-1}+1+n}+\tau_{i_{2j-1}+2}+\cdots
+\tau_{i_{2j}-1}+\sigma_{i_{2j}-1}\cr
\hbox{and\quad}
c_{\beta j} &= \tau_{i_{2j-1}}
+\sigma_{i_{2j-1}+n}+\tau_{i_{2j-1}+1}+\sigma_{i_{2j-1}+1}\cdots
+\sigma_{i_{2j}-1+n}+\tau_{i_{2j}},
}$$
while if $i_{2j}-i_{2j-1}$ is even, we may set
$$\eqalign{
c_{\alpha j} &= \sigma_{i_{2j-1}}+\tau_{i_{2j-1}+1}
+\sigma_{i_{2j-1}+1+n}+\tau_{i_{2j-1}+2}+\cdots
+\sigma_{i_{2j}-1+n}+\tau_{i_{2j}}\cr
\hbox{and\quad}
c_{\beta j} &= \tau_{i_{2j-1}}
+\sigma_{i_{2j-1}+n}+\tau_{i_{2j-1}+1}+\sigma_{i_{2j-1}+1}\cdots
+\tau_{i_{2j}-1}+\sigma_{i_{2j}-1}.
}$$
\qed

\proc{Theorem \checklabel\ThmDFour} Let $d=4$ and
let $\Gamma$ be an $n$--rung
M\"obius ladder with $n\geq 3$. Give $\Gamma$
the $G(4)$--coloring of Example \ExMobDFour. Then
$$H_1(\cover M) \cong \cases{
\bigoplus_{H\in\indexset^\star}H_1(M_H)\oplus \zt,
& if $n=3$;\cr
\bigoplus_{H\in\indexset^\star}H_1(M_H)
\oplus \zmod{8} \oplus \ztb^{4n-14},
&if $n\geq 4$.
}$$
\endproc

\prf This coloring is $4$--taut,
so Lemma \ResultSeqTaut\ applies and 
$\bigoplus_{H\in\indexset^\star}H_1(M_H)$ has odd order.
Therefore
$$H_1(\cover M) \cong \bigoplus_{H\in\indexset^\star}H_1(M_H)
\oplus \zmodb{8}^a \oplus \zmodb{4}^b \oplus \ztb^c$$
for some $a$, $b$ and $c$ with $3a+2b+c = 4n-11$, and 
$a+b+c=\dim H_1(\cover M;\zt)$.
If $n=3$, we must have $a=b=0$ and $c=1$, which proves this case of the 
theorem. From now on we assume that $n \geq 4$.

Let $H_0=\genby{x_1,x_2,x_3}\in \indexset^\star$, and
set $\Gamma_0=\Gamma_{H_0}$, $M_0=M_{H_0}$, and 
$\Delta_0=\Delta_{H_0}$; note that $\Gamma_0$ is the rim.
The chain  of subgroups 
$1\leq \genby{x_3} \leq H_0\leq G$
determines a chain of coverings
$\cover M \to M_1 \to M_0 \to M$, of which the middle one
has group $H_0/\genby{x_3}\cong \ztb^2$ and the others are 2--fold.
The branch set $\Delta_0$ of $\cover M\to M_0$ is a link
of $n$ components, any two of which have linking number 1
by Lemma \CalcLinking. The branch set of $M_1\to M_0$ is
a 3--component sublink $L_0$ of $\Delta_0$, lying
over the three rungs whose color is not $x_3$,
and it follows from Lemma \ModTwoGDisconn\ that
$H_1(M_1;\zt)\cong \zt$. Also, each component of $\Delta_0-L_0$
is covered by four simple closed curves in $M_1$,
so $\cover M\to M_1$ is branched over a link $\Delta_1$
of $4n-12$ components. Each element of $H_1(M_1;\zt)$
represented by a 
component of $\Delta_1$ has non-trivial image under the map
$H_1(M_1;\zt)\to H_1(M_0-L_0)$ defined just before Lemma 
\ModTwoGConn, and is 
therefore non-trivial. Hence $H_1(\Delta_1;\zt)\to
H_1(M_1;\zt)$ is onto, and the special case of Lemma \LinkCover\ shows
that $\dim H_1(\cover M;\zt)=4n-13$. It follows that
$2a+b = 2$. If we show that $a>0$, it will follow that
$a=1$, $b=0$ and $c=4n-14$, completing the proof.

Let $g_i\in G$ be the color of the rung $\tau_i$. Since
$n\geq 4$, there is at least one rung with color $x_3$,
which we may take to be $\tau_0$. Then $g_1\cdots g_{n-1}=1$,
and applying Lemma \HighOrder\ to the rungs 
$\tau_1,\ldots,\tau_{n-1}$ we see that it is enough to show
that 
$$z = \sum_{i=1}^{n-1}\sum_{H\in\indexset^\star}\achar_H(g_i)v_iH
\in C'_0(\Gamma\sep 3)$$
represents a non-zero element
of $H_0(\Gamma\sep 3)$.
Recall that  $C'(\Gamma\sep 3)$ is a subcomplex of
$C'(\Gamma\sep 4) = \bigoplus_{H\in\indexset^\star}C(\Gamma_H;\zt)$.
Since, for each $H\in\indexset^\star$,
$\Gamma_H$ is a circuit and 
there are an even number of $i$ 
($1\leq i\leq n-1$) with $\achar_H(g_i)=1$, $z$ is a boundary in
$C'(\Gamma\sep 4)$.
Let $c \in C'(\Gamma\sep 4)$,
with $c= \sum_{\sigma\in S_1(\Gamma),H\in\indexset^\star}
c(\sigma,H)\sigma H$,
and set 
$\phi(c) = \sum_{i=1}^{n-1}\sum_{H\in\indexset^\star} c(\tau_i,H)
\in\zt$.
If $c$ is a cycle, then 
$\phi(c) = 0$, so if $c_1$ and $c_2$ both have boundary $z$,
then $\phi(c_1)=\phi(c_2)$. On the other hand, if $c$ lies in
$C'(\Gamma\sep 3)$, then $\phi(c)=0$. Thus if we can find 
$c\in C'(\Gamma\sep 4)$
with $\partial c=z$ and $\phi(c)=1$, it will follow that
$z$ represents a non-zero element
of $H_0(\Gamma\sep 3)$. 

For $H\in \indexset^\star$, let 
$z_H = \sum_{i=1}^{n-1}\achar_H(g_i)v_i\in C_0(\Gamma_H;\zt)$,
so $z=\sum_{H\in\indexset^\star}z_H H$. A chain $c\in C'(\Gamma\sep 4)$
with $\partial c=z$ has the form $c=\sum_{H\in\indexset^\star}c_H H$
with $c_H\in C_1(\Gamma_H;\zt)$ and $\partial c_H =z_H$
for $H\in \indexset^\star$. Now $z_{H_0}=0$ and we may take
$c_{H_0}=0$. The remaining elements of $\indexset^\star(G)$
are in 2--1 correspondence with the elements of $\indexset^\star(H_0)$.
For $H\in\indexset^\star(H_0)$ let $H_1$ and $H_2$ be the two elements
of $\indexset^\star(G)$ with $H_1\cap H_0 =H = H_2\cap H_0$.
Then $\Gamma_{H_1}$ and $\Gamma_{H_2}$ contain the same rungs;
let $m_H$ be the number of these rungs, and $k_H$ the number of
them distinct from $\tau_0$. The union of $\Gamma_{H_1}$ and 
$\Gamma_{H_2}$ is an $m_H$--rung M\"obius ladder. If
$c_{H_1}\in C_1(\Gamma_{H_1};\zt)$ and 
$c_{H_2}\in C_1(\Gamma_{H_2};\zt)$ both have boundary
$z_{H_1}=z_{H_2}$, we may compute the sum $a_H$ of the coefficients
of the $\tau_i$ for $1\leq i\leq n-1$ in $c_{H_1}$ and $c_{H_2}$
using Lemma \MobParity\ (provided that $m_H$ and $k_H$ satisfy
the hypotheses of that lemma, as we shall see they do),
and then $c=\sum_{H\in\indexset^\star(H_0)}(c_{H_1}H_1+c_{H_2}H_2)$
is an element of $C_1'(\Gamma\sep 4)$ with $\partial c = z$
and $\phi(c)=\sum_{H\in\indexset^\star(H_0)}a_H$.
Now, for any $n$ and $H=\genby{x_1,x_3}$, $\genby{x_2,x_3}$
or $\genby{x_1x_2,x_3}$ we have $m_H=k_H=2$, and so $a_H=1$. For $n$
even we have $m_H=n-1$ and $k_H=n-2$ for $H=\genby{x_1,x_2x_3}$,
$\genby{x_2,x_1x_3}$ or $\genby{x_1x_2,x_1x_3}$,
and $m_H=n-3$ and $k_H=n-4$ for $H=\genby{x_1,x_2}$;
while for $n$ odd
we have $m_H=n-2$ and $k_H=n-3$ for $H=\genby{x_1,x_2}$,
$\genby{x_1,x_2x_3}$ or $\genby{x_2,x_1x_3}$,
and $m_H=n$ and $k_H=n-1$ for $H=\genby{x_1x_3,x_2x_3}$;
in all these cases, $a_H=0$.
This gives $\phi(c)=1$, completing the proof.
\qed

\proc{Theorem \checklabel\ThmDFive} Let $d=5$, and let $\Gamma$ 
be the Petersen graph with
the $G(5)$--coloring of Example \ExPetersen. Then
$$H_1(\cover M) \cong \bigoplus_{H\in\indexset^\star}H_1(M_H)
\oplus \zmod{16} \oplus \zmodb{4}^4 \oplus \ztb^{2}.$$
\endproc

\prf Let $S = \bigoplus_{H\in\indexset^\star}H_1(M_H)$.
This coloring is 5--taut, and therefore taut.
By Lemmas \ResultSeqOne\ and \ResultSeqTwo,
we may identify $S$ and $H_1(D(k))$ ($1\leq k \leq 4$)
with their images in $H_1(\cover M)$, so we have a filtration
$$S\leq H_1(\im\beta)=H_1(D(4))\leq H_1(D(3))\leq H_1(D(2))\leq 
H_1(D(1)) = H_1(\cover M).$$
Moreover, $H_1(D(4))/S \cong \zt$, and there are exact sequences
$$\eqalignno{
&0 \to H_1(D(4))/S \to H_1(D(3))/S \to \ztb^6 \to 0, 
&(\label\tempone)\cr
&0 \to H_1(D(3))/S \to H_1(D(2))/S \to \ztb^6 \to 0, 
&(\label\temptwo)\cr
\hbox{and}\quad &0 \to H_1(D(2))/S \to H_1(\cover M)/S \to \ztb \to 0.
&(\label\tempthree)
}$$
We show first that $H_1(\cover M)/S$ has an element of order 16.
Applying Lemma \HighOrder\ to
the edges $\tau_0,\ldots,\tau_4$, we see that it
is enough to show that
$$z = \sum_{i=0}^4\sum_{H\in\indexset^\star}
\achar_H(x_{i-1}x_{i+2})v_i H
\in C'_0(\Gamma\sep 4)$$ 
represents a non-zero element
of $H_0(\Gamma\sep 4)$.
Now  $C'(\Gamma\sep 4)$ is a subcomplex of
$C'(\Gamma\sep 5) = \bigoplus_{H\in\indexset^\star}C(\Gamma_H;\zt)$.
For each $H\in\indexset^\star$, 
$\achar_H(x_{i-1}x_{i+2})$ is non-zero
for an even number of $i$ ($0\leq i\leq 4$), 
and so $z$ is a boundary in
$C'(\Gamma\sep 5)$.
Let $c \in C'(\Gamma\sep 5)$,
with $c= \sum_{\sigma\in S_1(\Gamma),H\in\indexset^\star}
c(\sigma,H)\sigma H$,
and set 
$\phi(c) = \sum_{i=0}^4\sum_{H\in\indexset^\star} c(\tau_i,H)
\in\zt$.
If $c$ is a cycle, or if $c\in C'(\Gamma\sep 4)$, then 
$\phi(c) = 0$. Thus if we can find 
$c\in C'(\Gamma\sep 5)$
with $\partial c=z$ and $\phi(c)=1$, it will follow that
$z$ represents a non-zero element
of $H_0(\Gamma\sep 4)$. 

For $H\in \indexset^\star$, let 
$z_H = \sum_{i=0}^4\achar_H(x_{i-1}x_{i+2})v_i \in C_0(\Gamma_H;\zt)$.
If $c_H \in  C_1(\Gamma_H;\zt)$ has $\partial c_H = z_H$,
then $c = \sum_{H\in\indexset^\star}c_H H \in C'(\Gamma\sep 5)$
has $\partial c = z$. Let $H_0\in\indexset^\star$ have 
$\achar_{H_0}(x_i)=1$ for all $i$. Then $\Gamma_{H_0}$ is 
the outer rim
and $z_{H_0} = 0$, so we may take $c_{H_0}=0$. The remaining elements
of $\indexset^\star$ may be numbered as $H_{ji}$, $1\leq j\leq 6$
and $0\leq i \leq 4$. In Table 1 below, we list for
$H=H_{ji}$ the values of $\achar_H$ on the basis $x_0,\ldots, x_4$;
it will be apparent from the table that we have listed
every $H\neq H_0$.
$$\displaylines{\matrix{
&x_i&x_{i+1}&x_{i+2}&x_{i+3}&x_{i+4}\cr
H_{1i}:&1&0&0&0&0\cr
H_{2i}:&0&1&1&1&1\cr
H_{3i}:&1&1&0&0&0\cr
H_{4i}:&0&0&1&1&1\cr
H_{5i}:&1&0&1&0&0\cr
H_{6i}:&0&1&0&1&1\cr
}\cr
\hbox{Table 1.}
}$$
The $x_i$ are the colors on the $\sigma_i$; in Table 2 we list
the values of the $\achar_H$ on the colors of the other edges.
$$\displaylines{\matrix{
&\tau_i&\tau_{i+1}&\tau_{i+2}&\tau_{i+3}&\tau_{i+4}
&\rho_i&\rho_{i+1}&\rho_{i+2}&\rho_{i+3}&\rho_{i+4}\cr
H_{1i},H_{2i}:&0&1&0&1&0&1&1&0&0&0\cr
H_{3i},H_{4i}:&0&1&1&1&1&1&0&1&0&0\cr
H_{5i},H_{6i}:&1&1&0&0&0&1&1&1&1&0\cr
}\cr
\hbox{Table 2.}
}$$
We can read off the 0--chains $z_H$ from these tables;
we list these below, together with $c_H\in C_1(\Gamma;\zt)$
with $\partial c_H=z_H$; reference to the tables will show that
in fact $c_H\in C_1(\Gamma_H;\zt)$.
$$\eqalign{
H=H_{1i}:\quad z_H&=v_{i+1}+v_{i+3},\cr
H=H_{2i}:\quad z_H&=v_{i+1}+v_{i+3},\cr
H=H_{3i}:\quad z_H &= v_{i+1} + v_{i+2} + v_{i+3} + v_{i+4},\cr
H=H_{4i}:\quad z_H &= v_{i+1} + v_{i+2} + v_{i+3} + v_{i+4},\cr
H=H_{5i}:\quad z_H &= v_i + v_{i+1},\cr
H=H_{6i}:\quad z_H &= v_i + v_{i+1},\cr
}\quad\eqalign{
c_H&=\tau_{i+1};\cr
c_H&=\tau_{i+1};\cr
c_H &= \tau_{i+1}+\tau_{i+2};\cr
c_H &= \tau_{i+1}+\tau_{i+2};\cr
c_H &= \rho_i+\sigma_i+\rho_{i+1};\cr
c_H &= \tau_i+\rho_{i+2}+\sigma_{i+1}+\rho_{i+1}.\cr
}$$
Now $\phi\bigl(\sum_{H\in\indexset^\star}c_H H\bigr) = 1$,
and the proof that $H_1(\cover M)/S$ has an element of order 16
is complete.
It follows that $H_1(D(3))/S$ has an element of order 4.
Since $H_1(D(4))/S\cong\zt$, the sequence (\tempone) gives
$H_1(D(3))/S\cong \zmod{4}\oplus \ztb^5$. Also,
$H_1(D(2))/S$ has an element of order 8, so the sequence
(\temptwo) gives
$H_1(D(2))/S\cong \zmod{8}\oplus\zmodb{4}^a\oplus \ztb^b$
for some $a$ and $b$ with $2a+b=10$, and then (\tempthree)
gives $H_1(\cover M)/S\cong \zmod{16}\oplus\zmodb{4}^a\oplus \ztb^b$.
Since $S$ has odd order, we have
$H_1(\cover M) \cong S\oplus \zmod{16}\oplus\zmodb{4}^a\oplus \ztb^b$.

Consider the tower of coverings 
$\cover M \to M_2 \to M_1 \to M_0 \to M$
corresponding to the chain of subgroups 
$1\leq \genby{x_0x_1, x_4x_0} 
\leq \genby{x_0x_1, x_2x_3, x_4x_0} \leq H_0\leq G$.
Here $\cover M\to M_2$ has group $\genby{x_0x_1, x_4x_0}\cong\ztb^2$
and the others
are 2--fold. Now $M_0$ is a $\zt$ homology sphere,
and the branch set $\Delta_{H_0}$ of $\cover M \to M_0$
depends on the mod 2 linking number of the inner and outer rims
of $\Gamma$ in $M$. If this number is $0$ then $\Delta_{H_0}$
is a graph with vertices $w_0,\ldots,w_4,w'_0,\ldots,w'_4$
and edges $\{w_i,w_{i+2}\}$ colored $x_{i-1}x_{i+2}$,
$\{w'_i,w'_{i+2}\}$ also colored $x_{i-1}x_{i+2}$,
and $\{w_i,w'_i\}$ colored $x_{i-1}x_i$. 
If the linking number is $1$ then $\Delta_{H_0}$
is obtained from that graph by replacing, say, the edges
$\{w_0,w_2\}$ and $\{w'_0,w'_2\}$ by
$\{w_0,w'_2\}$ and $\{w'_0,w_2\}$.
In either case, the branch
set of $M_1\to M_0$ is a Hamiltonian circuit in $\Delta_{H_0}$,
and the edges not on this circuit are
$\{w_3,w'_3\}$, $\{w_1,w'_1\}$, $\{w_0,w'_0\}$,
$\{w_2,w_4\}$ and $\{w'_2,w'_4\}$.
Therefore $M_1$ is a $\zt$ homology sphere and the branch set
of $\cover M\to M_1$ is a link $L^1$ of 5 components.
We number the components as $L^1_0$, $L^1_1$, $L^1_2$, 
$L^1_{31}$ and $L^1_{32}$, where $L^1_0$ has color $x_2x_3$, 
$L^1_1$ has color $ x_0x_1$, $L^1_2$ has color $x_4x_0$, 
and $L^1_{31}$ and $L^1_{32}$ have color $x_1x_4$.
By Lemma \CalcLinking,
we have $\lk(L^1_i,L^1_{3j})=1$ for $0\leq i\leq 2$
and $j=1$ or $2$, and the linking number of any other pair of components 
except for $\{L^1_{31},L^1_{32}\}$ is $0$. 
Now the branch set of $M_2\to M_1$
is $L^1_0$, so $M_2$ is a $\zt$ homology sphere.
For $i=1$ or $2$, $L^1_i$ is covered by two simple
closed curves $L^2_{i1}$ and $L^2_{i2}$ in $M_2$,
while $L^1_{3i}$ is covered by a single curve
$L^2_{3i}$. The branch set of $\cover M\to M_2$
is the link with these six components.
There is a surface $F$ in $M_1$ with $\partial F= L^1_1$
and disjoint from $L^1_2$. Its inverse image in $M_2$
shows that 
$\lk(L^2_{11}, L^2_{21}) = \lk(L^2_{12}, L^2_{21})$ and
$\lk(L^2_{11}, L^2_{22}) = \lk(L^2_{12}, L^2_{22})$.
Switching the roles of $L^1_1$ and $L^1_2$ shows that
all four of these linking numbers are equal. Now,
for $i,j=1$ or $2$, there is a surface $F'$ in $M_1$
with $\partial F' = L^1_{3j}$ that meets
$L^1_i$ in a single point. Its inverse image shows that
$\lk(L^2_{i1},L^2_{3j})=\lk(L^2_{i2},L^2_{3j})=1$.
These linking numbers determine the matrix $\Lambda$ of
Lemma \ModTwoGDisconn\ for the covering
$\cover M \to M_2$; it 
has rank 2, and it follows that
$H_1(\cover M;\zt)\cong \ztb^7$.
Hence $a+b = 6$, so $a=4$ and $b=2$, and we are done.
\qed

\sect{References}
\item{{\bf \Flapan.}} Erica Flapan,
{\it Symmetries of M\"obius ladders},
Math. Ann. {\bf 283} (1989), 271--283.

\item{{\bf \LeeWein.}} Ronnie Lee and Steven H.~Weintraub,
{\it On the homology of double branched covers},
Proc. Amer. Math. Soc. {\bf 123} (1995), 1263--1266.

\item{{\bf \Massey.}} W.~S.~Massey,
{\it Completion of link modules},
Duke Math. J. {\bf 47} (1980), 399--420.

\item{{\bf \Sakuma.}} Makoto Sakuma,
{\it Homology of abelian coverings of links and spatial graphs}, 
Canad. J. Math. {\bf 47} (1995) no. 1, 201--224. 

\item{{\bf \Watkins.}} Mark E. Watkins,
{\it A theorem on Tait colorings with an application
to the generalized \it Petersen graphs}, 
J. Combin. Theory {\bf 6} (1969), 152--164. 

\item{{\bf \ZS.}} Oscar Zariski and Pierre Samuel,
{\it Commutative Algebra, Vol. II}, Springer-Verlag, New York,
1975; originally published by Van Nostrand, Princeton N.J., 1960.

\bigskip
\noindent Department of Mathematics
\par\noindent Louisiana State University
\par\noindent Baton Rouge, LA 70803, USA
\par\noindent E-mail: {\tt lither@marais.math.lsu.edu}
\bye